\documentclass[a4paper,10pt,twoside]{article}
\usepackage{german,amssymb,amsmath,latexsym,epsfig,graphics,graphicx,mathrsfs,color,amsthm,cite}
\usepackage[latin1]{inputenc}
\usepackage{fancyhdr} 
\newtheorem{theorem}{Theorem}[section]

\newtheorem{definition}[theorem]{Definition}
\newtheorem{lemma}[theorem]{Lemma}

\newtheorem{folgerung}[theorem]{Folgerung}
\newtheorem{beispiel}[theorem]{Beispiel}
\newtheorem{satz}[theorem]{Satz}
\newtheorem{bemerkung}[theorem]{Bemerkung}

\numberwithin{equation}{section}
\newcommand{\R}{{\mathbb R}}
\newcommand{\N}{{\mathbb N}}

\begin{document}
\thispagestyle{empty}
\rule{1.0\textwidth}{1.0pt}

\vspace*{5mm}
{\bf{\large Schwache lokale Extrema in Steuerungsproblemen mit \\[2mm] unendlichem Zeithorizont} \\[10mm]
Nico Tauchnitz} \\[25mm]
{\bf Vorwort} \\[2mm]
In der vorliegenden Ausarbeitung stelle ich notwendige und hinreichende Optimalit"atsbedingungen f"ur
schwache lokale Extrema in Aufgaben mit unendlichem Zeithorizont vor.
Im Vergleich zu den
klassischen Steuerungsproblemen besitzen Aufgaben mit unendlichem Zeithorizont ihren eigenen Charakter,
da durch das unbeschr"ankte Zeitintervall die Aufgabenstellung eine Singularit"at beinhaltet.
Die L"osungsmethoden, die f"ur die klassischen Aufgaben entwickelt wurden,
k"onnen f"ur die Situation des unendlichen Horizontes nicht einfach "ubernommen werden. \\[2mm]
Die Herleitung notwendiger Optimalit"atsbedingungen f"ur Aufgaben mit unendlichem Zeithorizont
auf der Grundlage der Ideen von Pontrjagin et. al., Dubovitskii \& Milyutin und Ioffe \& Tichomirov
ist ziemlich aufwendig und setzt ein tiefgreifendes Vorwissen "uber diese Verfahren voraus.
Dadurch motiviert wurde diese umfassende Ausarbeitung "uber schwache lokale Extrema,
das einfacher [nachzuvollziehen] ist als der Nachweis von Optimalit"atsbedingungen f"ur ein starkes lokales Minimum in \cite{TauchnitzPMPIHOC}.
Die Resultate "uber schwache lokale Extrema ersetzen kein Pontrjaginsches Maximumprinzip,
da Bedingungen f"ur die starke lokale Optimlit"at bei einer weitaus allgemeineren Klasse von
optimalen L"osungen greifen. \\[2mm]
Die vorliegende "uberarbeitete Version behandelt zus"atzlich die Aufgabe mit Zustandsbeschr"ankungen "uber dem unendlichen Zeithorizont.
Bei n"aherer Betrachtung der mathematischen Theorie zeigen sich aber die Schwierigkeiten und Pathologien,
die sich im Rahmen gewichteter Funktionenr"aume ergeben k"onnen.
Diesbez"uglich geben wir zum Abschluss dieser Arbeit ausf"uhrliche Bemerkungen an.
Ein alternativer Zugang zu starken lokalen Minimalstellen in Steuerungsproblemen mit unendlichem Zeithorizont
wird in der Arbeit \cite{TauchnitzPMPIHOC} angegeben. 
Eine umfassendere Darstellung kann man \cite{TauchnitzOC} entnehmen. \\[5mm]
Juli 2018

\newpage
\lhead[\thepage \, Inhaltsverzeichnis]{Optimale Steuerung mit unendlichem Zeithorizont}
\rhead[Optimale Steuerung mit unendlichem Zeithorizont]{Inhaltsverzeichnis \thepage}
\tableofcontents

\newpage
\lhead[\thepage \hspace*{1mm} Schwaches lokales Minimum]{}
\rhead[]{Schwaches lokales Minimum \hspace*{1mm} \thepage}
\section{Schwaches lokales Minimum "uber unendlichem Zeithorizont} \label{KapitelWeak}
Die Optimale Steuerung mit unendlichem Zeithorizont liefert die wesentliche Grundlage zur Formulierung und Untersuchung
von Aufgaben in der "Okonomischen Wachstumstheorie.
Im Rahmen der "Okonomischen Wachstumstheorie werden z.\,B. die Interaktionen sich "uberschneidender Generationen oder die Determinanten 
des wirtschaftlichen Wachstums,
insbesondere unter sich "andernden Umweltbedingungen wie globaler Erw"armung oder ersch"opfenden nat"urlichen Ressourcen, untersucht.
Aufgrund der Langlebigkeit der wirtschafts- und sozialpolitischen Entscheidungen muss dabei die Frage nach einem geeigneten 
Planungszeitraum aufgeworfen werden:
Jeder endliche Zeithorizont stellt die Forderung nach einer ad"aquaten Ausgangslage f"ur die nachfolgenden Generationen.
Um die Beachtung aller nachfolgenden Generationen zu gew"ahrleisten,
wird der Zeitrahmen in Form des unendlichen Zeithorizontes idealisiert (Arrow \& Kurz \cite{Arrow}). \\[2mm]
Der erste mathematische Beitrag zu einem Problem mit unbeschr"anktem Zeitintervall besteht in einer Aufgabe der Variationsrechnung,
in der die Frage nach der optimalen Sparquote einer Gesellschaft behandelt wird (Ramsey \cite{Ramsey}).
Bei n"aherer Betrachtung entsteht in diesem Problem nicht nur die Aufgabe die Optimalit"atsbedingungen "uber dem unendlichen Zeithorizont zu formulieren,
sondern insbesondere die Gestalt der Transversalit"atsbedingungen im Unendlichen zu charakterisieren. 
Diese Fragestellung stellt allerdings eine schwerwiegende Herausforderung dar,
denn die bekannten Resultate k"onnen nicht einfach "ubernommen und die bekannten Methoden k"onnen nicht
unmittelbar an das unbeschr"ankte Intervall angepasst werden. \\[2mm]
Von den wenigen uns bekannten Resultaten,
die den unendlichen Zeithorizont nicht auf ein endliches Intervall reduzieren,
z"ahlen wir die Arbeiten von Brodskii \cite{Brodskii} und Pickenhain \cite{Pickenhain} auf.
In diesen Beitr"agen werden auf der Basis von schwachen lokalen Variationen und der Anwendung geeigneter funktionalanalytischer Methoden
notwendige Optimalit"atsbedingungen f"ur Steuerungsprobleme mit unendlichem Zeithorizont erzielt. \\
In \cite{Brodskii} wird eine sehr allgemeine Aufgabenklasse mit Zustandsbeschr"ankungen und Randbedingungen im Unendlichen betrachtet.
Aufgrund der Wahl des Raumes der messbaren und beschr"ankten Funktionen f"uhrt der funktionalanalytische Rahmen allerdings
zu keiner ``"asthetischen'' Darstellung der Adjungierten. \\
Demgegen"uber bezieht sich die Arbeit \cite{Pickenhain} auf Aufgaben mit eindimensionalen linearen Nebenbedingungen und mit freiem rechten Endpunkt.
Der innovative Beitrag in \cite{Pickenhain} ist die Wahl des gewichteten Sobolev-Raumes,
in dessen Rahmen die ``nat"urlichen'' Transversalit"atsbedingungen in gewisser Weise Eigenschaften der Elemente dieser R"aume sind.
An der Arbeit \cite{Pickenhain} ist jedoch anzumerken,
dass sich die angewandte Beweismethode auf linear-quadratische Aufgaben fokussiert.
Die Anwendbarkeit der Beweisstrategie im Fall von Aufgaben mit einer nichtlinearen Dynamik ist offen. \\[2mm]
In dieser Arbeit stellen wir einen Zugang zu Steuerungsproblemen mit unendlichem Zeithorizont im Rahmen gewichteter Funktionenr"aume vor.
Das gestellte Steuerungsproblem besitzt einen weniger allgemeinen Charakter als die Aufgabe bei Brodskii \cite{Brodskii}.
Aber auf der Grundlage der Wahl gewichteter Sobolev-R"aume f"ur die Zustandstrajektorie nach Pickenhain \cite{Pickenhain}
erhalten wir in den notwendigen Optimalit"atsbedingungen qualitativ bessere Informationen "uber die Adjungierte.
Aufgrund der G"ultigkeit der ``nat"urlichen'' Transversalit"atsbedingungen lassen sich hinreichende Bedingungen
unmittelbar in die Ausf"uhrungen einbinden. \\[2mm]
Der wesentliche Vorteil der Wahl gewichteter Funktionenr"aume ist,
dass eine recht umfassende Menge an zul"assigen Steuerungsprozessen betrachtet werden darf.
Allerding treten dadurch Schwierigkeiten auf,
die einerseits mit dem expandierenden Umgebungsbegriff verbunden sind.
Andererseits ergeben sich bei der Verallgemeinerung der gundlegenden Aufgabenklasse Probleme,
insbesondere in der Aufgabe mit Zustandsbschr"ankungen oder bei Aufgaben mit Randwerten im Unendlichen,
tiefgreifende H"urden, zu denen wir im letzten Abschnitt ausf"uhrliche Bemerkungen geben.


\section{Die Aufgabenstellung}
Wir untersuchen schwache lokale Minimalstellen der Aufgabe
\begin{eqnarray}
&& \label{WMPAufgabe1} J\big(x(\cdot),u(\cdot)\big) = \int_0^\infty \omega(t) f\big(t,x(t),u(t)\big) \, dt \to \inf, \\
&& \label{WMPAufgabe2} \dot{x}(t) = \varphi\big(t,x(t),u(t)\big), \qquad x(0)=x_0, \\
&& \label{WMPAufgabe3} u(t) \in U \subseteq \R^m, \quad U \not= \emptyset \mbox{ und konvex}, \\
&& \label{WMPAufgabe4} g_j\big(t,x(t)\big) \leq 0 \quad \mbox{f"ur alle } t \in \R_+, \quad j=1,...,l.
\end{eqnarray}
Dabei gelten $f:\R \times \R^n \times \R^m \to \R$, $\varphi:\R \times \R^n \times \R^m \to \R^n$ und $g_j(t,x):\R \times \R^n \to \R$. \\
Wir nennen die Trajektorie $x(\cdot)$ eine L"osung des dynamischen Systems (\ref{WMPAufgabe2}) zur Anfangsbedingung $x(0)=x_0$,
falls $x(\cdot)$ auf $\R_+$ definiert ist und auf jedem endlichen Intervall die Dynamik mit Steuerung $u(\cdot)$
im Sinn von Carath\'eodory l"ost. \\[2mm]
In der Aufgabe (\ref{WMPAufgabe1})--(\ref{WMPAufgabe4}) gilt stets die Annahme:
\begin{enumerate}
\item[(A$_0$)] Es sei $\omega(\cdot) \in L_1(\R_+,\R_+)$ und es sei $\nu(t) = e^{-at}$ ein Gewicht mit $a>0$.
\end{enumerate}
F"ur die Aufgabe (\ref{WMPAufgabe1})--(\ref{WMPAufgabe4}) betrachten wir Variationen im Raum gewichteter stetiger Funktionen,
die im Unendlichen verschwindenn.
Da unter den Eigenschaften der Gewichtsfunktion $\nu(t) = e^{-at}$ nach Lemma \ref{LemmaGrenzwert2} die Implikationen
$$x(\cdot) \in W^1_2(\R_+,\R^n;\nu) \quad\Rightarrow\quad \nu(\cdot)x(\cdot) \in W^1_1(\R_+,\R^n) \quad\Rightarrow\quad x(\cdot) \in C_0(\R_+,\R^n;\nu)$$
gelten,
formulieren wir wie Pickenhain \cite{Pickenhain} die Aufgabe (\ref{WMPAufgabe1})--(\ref{WMPAufgabe4}) im Rahmen gewichteter Sobolev-R"aume. \\[2mm]
Wir definieren zu $\big(x(\cdot),u(\cdot)\big) \in W^1_2(\R_+,\R^n;\nu) \times L_\infty(\R_+,U)$ die Menge $U_{\gamma,\nu}$ wie folgt:
$$U_{\gamma,\nu}= \{ (t,x,u) \in \R_+ \times \R^n \times \R^m\,|\, \nu(t)\|x-x(t)\| \leq \gamma, \|u-u(t)\| \leq \gamma\}.$$
Dann geh"oren zur Menge $\mathscr{A}_{\rm Lip}$ diejenigen $\big(x(\cdot),u(\cdot)\big) \in W^1_2(\R_+,\R^n;\nu) \times L_\infty(\R_+,U)$,
f"ur die es eine Zahl $\gamma>0$ derart gibt, dass auf dem Abschluss der Menge $U_{\gamma,\nu}$ gelten:
\begin{enumerate}
\item[(A$_1$)] Die Abbildungen $f(t,x,u)$, $\varphi(t,x,u)$ und $g_j(t,x)$ sind stetig differenzierbar.
\item[(A$_2$)] F"ur alle $(t,x,u) \in U_{\gamma,\nu}$ ein $L(\cdot) \in L_1(\R_+,\R;\omega)$ und ein $C_0>0$ mit
               \begin{eqnarray*}
               && \big\|\big(f(t,x,u),f_u(t,x,u)\big) \big\| \leq L(t), \quad \|f_x(t,x,u)\| \leq L(t)\nu(t) \\
               && \|\varphi(t,x,u)\| \leq C_0 (1+\|x\|+\|u\|), \quad
                  \big\|\big(\varphi_x(t,x,u),\varphi_u(t,x,u)\big)\big\| \leq C_0.
               \end{eqnarray*}    
               Au"serdem gilt f"ur alle $(t,x,u),(t,x',u') \in U_{\gamma,\nu}$:
               \begin{eqnarray*}
               && \|\varphi_x(t,x',u') - \varphi_x(t,x,u) \| \leq C_0 \big(e^{-at}\|x'-x\| + \|u'-u\|\big), \\
               && \|\varphi_u(t,x',u') - \varphi_u(t,x,u) \| \leq C_0 \big(\|x'-x\| + e^{at} \|u'-u\|\big).
               \end{eqnarray*}
\item[(A$_3$)] F"ur alle $(t,x), (t,x') \in U_{\gamma,\nu}$ existiert ein $C_0>0$ mit
               $$|g_j(t,x)| \leq C_0 (1+\|x\|), \quad \|g_{jx}(t,x)\| \leq C_0, \quad \|g_{jx}(t,x)-g_{jx}(t,x')\| \leq C_0 e^{-at}\|x - x'\|.$$
\end{enumerate}
Die Terme $e^{-at}, e^{at}$ sind eine Konsequenz der Variation in der expandierenden Umgebung $U_{\gamma,\nu}$. \\[2mm]
Wir nennen $\big(x(\cdot),u(\cdot)\big) \in W^1_2(\R_+,\R^n;\nu) \times L_\infty(\R_+,U)$
einen zul"assigen Steuerungsprozess in der Aufgabe (\ref{WMPAufgabe1})--(\ref{WMPAufgabe4}),
falls $\big(x(\cdot),u(\cdot)\big)$ dem System (\ref{WMPAufgabe2}) gen"ugt, die Zustandsbeschr"ankungen (\ref{WMPAufgabe4}) erf"ullt
und das Lebesgue-Integral im Zielfunktional in (\ref{WMPAufgabe1}) endlich ist.
Die Menge $\mathscr{A}_{\rm adm}$ bezeichnet die Menge der zul"assigen Steuerungsprozesse $\big(x(\cdot),u(\cdot)\big)$. \\[2mm]
Ein zul"assiger Steuerungsprozess $\big(x_*(\cdot),u_*(\cdot)\big)$ ist eine schwache lokale Minimalstelle\index{Minimum, schwaches lokales}
der Aufgabe (\ref{WMPAufgabe1})--(\ref{WMPAufgabe4}),
falls eine Zahl $\varepsilon > 0$ derart existiert, dass die Ungleichung 
$$J\big(x(\cdot),u(\cdot)\big) \geq J\big(x_*(\cdot),u_*(\cdot)\big)$$
f"ur alle $\big(x(\cdot),u(\cdot)\big) \in \mathscr{A}_{\rm adm}$ mit 
$\|x(\cdot)-x_*(\cdot)\|_\infty \leq \varepsilon$, $\|u(\cdot)-u_*(\cdot)\|_{L_\infty} \leq \varepsilon$ gilt. \\[2mm]
Als abschlie"sende Bemerkung weisen wir darauf hin,
dass die Variationen im gewichteten Rahmen streng genommen auf einen $\nu$-lokalen Optimalit"atsbegriff f"uhren.
D.\,h., dass die Ungleichung
$$J\big(x(\cdot),u(\cdot)\big) \geq J\big(x_*(\cdot),u_*(\cdot)\big)$$
f"ur alle $\big(x(\cdot),u(\cdot)\big) \in \mathscr{A}_{\rm adm}$ mit 
$\|x(\cdot)-x_*(\cdot)\|_{\infty,\nu} \leq \varepsilon$, $\|u(\cdot)-u_*(\cdot)\|_{L_\infty} \leq \varepsilon$ gilt.
Da aber in der Ungleichung $\|x(\cdot)-x_*(\cdot)\|_{\infty,\nu} \leq \varepsilon$ die Elemente $x(\cdot)$ mit der Eigenschaft
$\|x(\cdot)-x_*(\cdot)\|_\infty \leq \varepsilon$ inbegriffen sind,
entsteht kein Widerspruch.
       
       \rhead[]{Schwaches Optimalit"atsprinzip \hspace*{1mm} \thepage}
       \section{Ein Schwaches Optimalit"atsprinzip} \label{AbschnittWeak}
\subsection{Notwendige Optimalit"atsbedingungen}
Im Weiteren bezeichnet $H: \R \times \R^n \times \R^m \times \R^n \times \R \to \R$ die Pontrjagin-Funktion
$$H(t,x,u,p,\lambda_0) = \langle p, \varphi(t,x,u) \rangle-\lambda_0 \omega(t)f(t,x,u).$$

\begin{theorem} \label{SatzWMP} \index{Schwaches Optimalit"atsprinzip}
Es sei $\big(x_*(\cdot),u_*(\cdot)\big) \in \mathscr{A}_{\rm adm} \cap \mathscr{A}_{\rm Lip}$ und $x_*(\cdot) \in W^1_2(\R_+,\R^n;\nu)$.
Ferner sei die folgende Kontraktionsbedingung erf"ullt:
\begin{equation} \label{WMPBedingung}
\sup_{t \in \R_+} \int_0^t \frac{\nu(t)}{\nu(s)} \big\|\varphi_x\big(s,x_*(s),u_*(s)\big)\big\| \, ds < 1.
\end{equation}
Ist $\big(x_*(\cdot),u_*(\cdot)\big)$ ein schwaches lokales Minimum der Aufgabe (\ref{WMPAufgabe1})--(\ref{WMPAufgabe3}),
dann existieren nicht gleichzeitig verschwindende Multiplikatoren $\lambda_0 \geq 0$ und $p(\cdot) \in L_2(\R_+,\R^n;\nu^{-1})$
derart, dass
\begin{enumerate}
\item[(a)] die Funktion $p(\cdot)$ fast "uberall der adjungierten Gleichung\index{adjungierte Gleichung}
           \begin{equation}\label{WMP1}
           \dot{p}(t) = -\varphi_x^T\big(t,x_*(t),u_*(t)\big) p(t) + \lambda_0 \omega(t)f_x\big(t,x_*(t),u_*(t)\big)
           \end{equation}
           gen"ugt und die ``nat"urlichen'' Transversalit"atsbedingungen\index{Transversalit"atsbedingungen!nat@--, nat"urliche}
           \begin{equation}\label{WMP2}
           \lim_{t \to \infty} \|p(t)\|^2\nu^{-1}(t) = 0, \qquad
           \lim_{t \to \infty} \langle p(t),x(t) \rangle = 0 \quad \forall \; x(\cdot) \in W^1_2(\R_+,\R^n;\nu)
           \end{equation}
           erf"ullt;
\item[(b)] in fast allen Punkten $t \in \R_+$ und f"ur alle $u \in U$ die Variationsungleichung
           \begin{equation}\label{WMP3}
           \big\langle H_u\big(t,x_*(t),u_*(t),p(t),\lambda_0\big),\big(u-u_*(t)\big) \big\rangle\leq 0
           \end{equation}
           gilt.
\end{enumerate}
\end{theorem}

\begin{bemerkung}{\rm
Die Kontraktionsbedingung (\ref{WMPBedingung}) stellt eine wesentliche Voraussetzung dar,
auf die wir in der Beweisf"uhrung zur"uckgreifen.
Diese Bedingung bezieht sich auf die Frage,
f"ur welche $z(\cdot)$ die lineare Intregralgleichung
$$x(t)+ \int_0^t A(s) x(s) \, ds = z(t), \qquad t \in \R_+,$$
eine L"osung $x(\cdot)$ besitzt (Lemma \ref{LemmaDGL2}).
Dabei erweist sich der korrespondierende Integraloperator im Rahmen gewichteter stetiger Funktionen als kontraktiv,
wenn (\ref{WMPBedingung}) f"ur $A(t) = \varphi_x\big(t,x_*(t),u_*(t)\big)$ erf"ullt ist.
Beachten wir nun f"ur die Gewichte $\nu(t)=e^{-at}$ mit $a>0$ die Gleichung
$$\sup_{t \in \R_+} \int_0^t \frac{\nu(t)}{\nu(s)} \, ds = \frac{1}{a},$$
so lassen sich stets Gewichtsfunktionen $\nu(\cdot)$ angeben,
die bei der messbar und beschr"ankten Funktion $t \to \varphi_x\big(t,x_*(t),u_*(t)\big)$ die Bedingung (\ref{WMPBedingung}) erf"ullen. \\
Allerdings ist im Zugang "uber gewichtete R"aume stetiger Funktionen der expandierende Umgebungsradius zu beachten:
$$\|x(\cdot)\|_{\infty,\nu} \leq \varepsilon \quad\Leftrightarrow\quad \|x(t)\| \leq \varepsilon e^{at} \mbox{ f"ur alle } t \in \R_+.$$
Je gr"o"ser der Parameter $a$ gew"ahlt wurde,
desto schwieriger gestaltet sich der Umgang mit den Daten der Aufgabe. \hfill $\square$}
\end{bemerkung}

\begin{lemma} \label{LemmaMichel}
Zus"atzlich zu Theorem \ref{SatzWMP} gelte f"ur $L(\cdot) \in L_1(\R_+,\R_+;\omega)$ in (A$_2$):
$$\lim_{t \to \infty} \omega(t)L(t)=0.$$
Dann ist die Bedingung von Michel \cite{Michel}\index{Transversalit"atsbedingungen!von@-- von Michel} erf"ullt:
\begin{equation} \label{NaturalHamilton} \lim_{t \to \infty} H\big(t,x_*(t),u_*(t),p(t),\lambda_0\big)=0. \end{equation}
\end{lemma}
 
{\bf Beweis} Da $\big(x_*(\cdot),u_*(\cdot)\big)$ zur Menge $\mathscr{A}_{\rm Lip}$ geh"oren,
gilt nach Voraussetzung (A$_2$):
\begin{eqnarray*}
&& \hspace*{-10mm}\big|H\big(t,x_*(t),u_*(t),p(t),\lambda_0\big)\big|^2
         \leq \big[\|p(t)\| \|\varphi\big(t,x_*(t),u_*(t)\big)\| + \omega(t)\|f\big(t,x_*(t),u_*(t)\big)\| \big]^2 \\
&& \leq C \cdot \Big[\big(\|p(t)\|^2 \nu^{-1}(t) \big) \cdot \big((1+\|x_*(t)\|^2+\|u_*(t)\|^2)\nu(t)\big)
    + \big(\omega(t)L(t)\big)^2\Big].
\end{eqnarray*}
Der letzte Ausdruck verschwindet im Unendlichen, da nach Lemma \ref{LemmaGrenzwert2}
$$\lim_{t \to \infty} \|x_*(t)\|^2\nu(t) = 0$$ 
und nach (\ref{WMP2}) die ``nat"urliche'' Transversalit"atsbedingung
$$\lim_{t \to \infty} \|p(t)\|^2 \nu^{-1}(t) = 0$$
gelten.
Weiterhin sind nach Voraussetzung und wegen der Beschr"anktheit von $u_*(\cdot)$
$$\lim_{t \to \infty} \omega(t)L(t)=0,\qquad \lim_{t \to \infty} \big(1+\|u_*(t)\|^2\big)\nu(t) = 0$$
erf"ullt. \hfill $\blacksquare$

\begin{beispiel}\label{BeispielDockner} {\rm Wir betrachten nach Dockner et\,al. \cite{Dockner} das Differentialspiel
\begin{eqnarray*}
&& \tilde{J}_i\big(x(\cdot),u_1(\cdot),u_2(\cdot)\big) =\int_0^\infty \omega(t)\big(p x(t)-c_i\big)u_i(t) \, dt \to \sup, \\
&& \dot{x}(t)=x(t) \big(\alpha-r\ln x(t) \big) -u_1(t)x(t)-u_2(t)x(t), \quad x(0)=x_0>0, \\
&& u_i > 0, \quad \alpha,c_i,p,r,\varrho >0, \quad \alpha> \frac{1}{c_1+c_2}, \quad i=1,2.
\end{eqnarray*}
Im Gegensatz zu \cite{Dockner} sei $\omega(\cdot) \in L_1(\R_+,\R_+)$ nicht auf $\omega(t)=e^{-\varrho t}$ festgelegt und
kann durch eine Weibull-Verteilung $\omega(t)=t^{k-1} e^{-t^k}$,
insbesondere mit einem Formparameter $k \in (0,1)$, beschrieben werden.
Es sei der Preis $p$ nicht konstant, sondern umgekehrt proportional zum Angebot:
$$p=p(u_1x+u_2x)= \frac{1}{u_1x+u_2x}.$$
Dieser Ansatz spiegelt die "Okonomie einer ``Eskimo''-Gesellschaft wider,
in der der Fischbestand die wichtigste Nahrungsgrundlage darstellt und kein echtes Ersatzprodukt existiert.
Unter diesen Umst"anden f"uhrt eine prozentuale Preissteigerung zu einem Umsatzr"uckgang in gleicher Relation. \\[2mm]
Nach Anwendung der Transformation $z=\ln x$ ergibt sich das Spielproblem
\begin{eqnarray*}
&& J_i\big(z(\cdot),u_1(\cdot),u_2(\cdot)\big)
        =\int_0^\infty \omega(t)\bigg(\frac{1}{u_1(t)+u_2(t)}-c_i\bigg)u_i(t) \, dt \to \sup, \\
&& \dot{z}(t)=-r z(t)+\alpha-u_1(t)-u_2(t), \quad z(0)=\ln x_0>0, \\
&& u_i>0, \quad \alpha,c_i,p,r,\varrho >0, \quad \alpha> \frac{1}{c_1+c_2}, \quad i=1,2.
\end{eqnarray*}
In diesem Spielproblem sind die Dynamiken linear und au"serdem flie"st die Zustandsvariable nicht im Integranden ein.
W"ahlen wir nun $\nu(t) = e^{-at}$ mit $a >r$,
so sind f"ur jedes $\omega(\cdot) \in L_1(\R_+,\R_+)$ s"amtliche Annahmen in Theorem \ref{SatzWMP},
insbesondere die Kontraktionsbedingung (\ref{WMPBedingung}), erf"ullt. \\
Der L"osungsansatz "uber ein Nash-Gleichgewicht liefert die Steuerungen
$$u_1^*(t)\equiv \frac{c_2}{(c_1+c_2)^2}, \qquad u_2^*(t)\equiv \frac{c_1}{(c_1+c_2)^2},$$
die optimale Trajektorie
$$z_*(t)=(z_0-c_0)e^{-rt} + c_0, \qquad c_0=  \frac{1}{r}\bigg(\alpha-\frac{1}{c_1+c_2}\bigg)$$
und die Adjungierten
$$p_i(t) \equiv 0, \qquad i=1,2.$$
Die Funktion $z_*(\cdot)$ ist streng monoton und nimmt nur Werte des Segments $[z_0,c_0]$ an.
Da $c_0$ und $z_0$ positiv sind, ist $x_*(t)=\exp\big(z_*(t)\big)$ "uber $\R_+$ wohldefiniert. \hfill $\square$}
\end{beispiel}
       \subsection{Der Nachweis der notwendigen Optimalit"atsbedingungen} \label{AbschnittBeweisWMP}
Es sei $\big(x_*(\cdot),u_*(\cdot)\big) \in \mathscr{A}_{\rm Lip}$.
Da die Funktion $y(\cdot) = \nu(\cdot) x(\cdot)$ nach Lemma \ref{LemmaGrenzwert2} dem Raum $W^1_1(\R_+,\R^n)$ angeh"ort,
verschwindet diese Funktion im Unendlichen.
Also ist es gerechtfertigt,
die Extremalaufgabe in den gewichteten Raum $C_0(\R_+,\R^n;\nu)$ einzubetten. \\[2mm]
Die expandierende Menge $U_{\gamma,\nu}$ ist wie folgt definiert:
$$U_{\gamma,\nu}= \{ (t,x,u) \in \R_+ \times \R^n \times \R^m\,|\, e^{-at}\|x-x(t)\| \leq \gamma, \|u-u(t)\| \leq \gamma\}.$$
Wir betrachten f"ur $\big(x(\cdot),u(\cdot)\big) \in C_0(\R_+,\R^n;\nu) \times L_\infty(\R_+,\R^m)$ die Abbildungen
\begin{eqnarray*}
J\big(x(\cdot),u(\cdot)\big) &=& \int_{t_0}^{t_1} \omega(t) f\big(t,x(t),u(t)\big) \, dt, \\
F\big(x(\cdot),u(\cdot)\big)(t) &=& x(t) -x(t_0) -\int_{t_0}^t \varphi\big(s,x(s),u(s)\big) \, ds, \quad t \in \R_+,\\
H_0\big(x(\cdot)\big) &=& x(0).
\end{eqnarray*}
Dabei fassen wir sie als Abbildungen zwischen folgenden Funktionenr"aumen auf:
\begin{eqnarray*}
J &:& C_0(\R_+,\R^n;\nu) \times L_\infty(\R_+,\R^m) \to \R, \\
F &:& C_0(\R_+,\R^n;\nu) \times L_\infty(\R_+,\R^m) \to C_0(\R_+,\R^n;\nu), \\
H_0 &:& C_0(\R_+,\R^n;\nu) \to \R^n.
\end{eqnarray*}

Wir setzen $\mathscr{F}=(F,H_0)$ und pr"ufen f"ur die Extremalaufgabe
\begin{equation} \label{ExtremalaufgabeWMP}
J\big(x(\cdot),u(\cdot)\big) \to \inf, \qquad \mathscr{F}\big(x(\cdot),u(\cdot)\big)=0, \qquad u(\cdot) \in L_\infty(\R_+,U)
\end{equation}
die Voraussetzungen von Theorem \ref{SatzExtremalprinzipSchwach} im Punkt $\big(x_*(\cdot),u_*(\cdot)\big)$:

\begin{enumerate}
\item[(A$_1$)] Das Funktional $J$ ist nach Beispiel \ref{DiffZielfunktionalW}
               im Punkt $\big(x_*(\cdot),u_*(\cdot)\big)$ Fr\'echet-differenzierbar.
\item[(A$_2$)] Die Abbildung $F$ ist die Summe der Abbildung $x(\cdot) \to x(t)$ und der Abbildung
               $$\big(x(\cdot),u(\cdot)\big) \to -\int_{t_0}^t \varphi\big(s,x(s),u(s)\big) \, ds.$$
               Im Beispiel \ref{DiffDynamikW1} ist nachgewiesen, dass die Abbildung $F$ in den Raum $C_0(\R_+,\R^n;\nu)$ abbildet
               und im Punkt $\big(x_*(\cdot),u_*(\cdot)\big)$ stetig Fr\'echet-differenzierbar ist.
               F"ur die Abbildung $H_0$ ist die stetige Fr\'echet-Differenzierbarkeit offensichtlich.
\item[(B)] Wir setzen $A(t)=\varphi_x\big(t,x_*(t),u_*(t)\big)$.
           Dann ist die Surjektivit"at des Operators $F_x\big(x_*(\cdot),u_*(\cdot)\big)$ "aquivalent zur Aussage,
           dass f"ur jedes $y(\cdot) \in C_0(\R_+,\R^n;\nu)$ die Gleichung
           $$x(t) - \int_0^t A(s)x(s) \, ds =y(t), \qquad t \in \R_+,$$
           l"osbar ist.            
           Dies besagt unter der Bedingung (\ref{WMPBedingung}) in Theorem \ref{SatzWMP} gerade Lemma \ref{LemmaDGL2}.
           Somit besitzt der Operator $\mathscr{F}_x\big(x_*(\cdot),u_*(\cdot)\big)$ eine endliche Kodimension.
\end{enumerate}

Zur Extremalaufgabe (\ref{ExtremalaufgabeWMP}) definieren wir auf
$$C_0(\R_+,\R^n;\nu) \times L_\infty(\R_+,\R^m) \times \R \times C_0^*(\R_+,\R^n;\nu) \times \R^n$$
die Lagrange-Funktion $\mathscr{L}=\mathscr{L}\big(x(\cdot),u(\cdot),\lambda_0,y^*,l_0\big)$,
$$\mathscr{L}= \lambda_0 J\big(x(\cdot),u(\cdot)\big)+ \big\langle y^*, F\big(x(\cdot),u(\cdot)\big) \big\rangle
                         +l_0^T H_0\big(x(\cdot)\big).$$
Ist $\big(x_*(\cdot),u_*(\cdot)\big)$ eine schwache lokale Minimalstelle der Aufgabe (\ref{ExtremalaufgabeWMP}),
dann existieren nach Theorem \ref{SatzExtremalprinzipSchwach}
nicht gleichzeitig verschwindende Lagrangesche Multiplikatoren $\lambda_0 \geq 0$, $y^* \in C_0^*(\R_+,\R^n;\nu)$ und $l_0 \in \R^n$ derart,
dass gelten:
\begin{enumerate}
\item[(a)] Die Lagrange-Funktion besitzt bez"uglich $x(\cdot)$ in $x_*(\cdot)$ einen station"aren Punkt, d.\,h.
          \begin{equation} \label{SatzWMPLMR1}
          \mathscr{L}_x\big(x_*(\cdot),u_*(\cdot),\lambda_0,y^*,l_0\big)=0;
          \end{equation}         
\item[(b)] Die Lagrange-Funktion erf"ullt bez"uglich $u(\cdot)$ in $u_*(\cdot)$ die Variationsungleichung
          \begin{equation} \label{SatzWMPLMR2}
           \big\langle \mathscr{L}_u\big(x_*(\cdot),u_*(\cdot),\lambda_0,y^*,l_0\big), u(\cdot)-u_*(\cdot) \big\rangle \geq 0
          \end{equation}
          f"ur alle $u(\cdot) \in L_\infty(\R_+,U)$.
\end{enumerate}

Aufgrund (\ref{SatzWMPLMR1}) ist folgende Variationsgleichung f"ur alle $x(\cdot) \in C_0(\R_+,\R^n;\nu)$ erf"ullt: 
\begin{eqnarray}
0 &=& \lambda_0 \cdot \int_0^\infty \omega(t) \big\langle f_x\big(t,x_*(t),u_*(t)\big),x(t) \big\rangle\, dt + l_0^T x(0)  \nonumber \\
  & & \label{BeweisschlussWMP1} 
      + \int_0^\infty \nu(t)\bigg[ x(t)-x(0) - \int_0^t \varphi_x\big(s,x_*(s),u_*(s)\big) x(s) \,ds \bigg]^T  d\mu(t).
\end{eqnarray}
Dabei ist $\mu$ nach Folgerung \ref{FolgerungRieszC0nu} ein signiertes regul"ares Borelsches Vektorma"s "uber $\R_+$. \\[2mm]
In der Gleichung (\ref{BeweisschlussWMP1}) ist jeder Integralterm absolut integrierbar; insbesondere gilt
\begin{eqnarray*}
\lefteqn{\int_0^\infty \nu(t)\bigg[ \int_0^t \big\|\varphi_x\big(s,x_*(s),u_*(s)\big) x(s)\big\| \,ds \bigg]  d|\mu|(t)} \\
&\leq& \int_0^\infty \bigg[ \int_0^t \frac{\nu(t)}{\nu(s)}\big\|\varphi_x\big(s,x_*(s),u_*(s)\big) [\nu(s)x(s)]\big\| \,ds \bigg]  d|\mu|(t) \\
&\leq& \sup_{t \in \R_+} \int_0^t \frac{\nu(t)}{\nu(s)}\big\|\varphi_x\big(s,x_*(s),u_*(s)\big)\big\| \,ds \cdot
       \int_0^\infty d|\mu|(t) \cdot \|x(\cdot) \|_{\infty,\nu} \leq \|\mu\| \cdot \|x(\cdot) \|_{\infty,\nu}.
\end{eqnarray*}
Damit stellt die rechte Seite in (\ref{BeweisschlussWMP1}) ein stetiges lineares Funktional im Raum $C_0(\R_+,\R^n;\nu)$ dar.
Durch vertauschen der Integrationsreihenfolge im letzten Term in (\ref{BeweisschlussWMP1}) bringen wir diese Gleichung in die Form
\begin{eqnarray}
0 &=& \int_0^\infty \bigg[ \lambda_0 \omega(t)f_x\big(t,x_*(t),u_*(t)\big) - \varphi^T_x\big(t,x_*(t),u_*(t)\big) \int_t^\infty
      \nu(s) d\mu(s) \bigg]^T  x(t) dt \nonumber \\
  & & \label{BeweisschlussWMP2} + \int_0^\infty \nu(t) [x(t)]^T d\mu(t)
            + \bigg[l_0^T x(0) - \int_0^\infty \nu(t) [x(0)]^T d\mu(t)\bigg].
\end{eqnarray}
Setzen wir $p(t)=\displaystyle \int_t^\infty \nu(s) \, d\mu(s)$,
so erhalten wir aus der eindeutigen Darstellung eines stetigen linearen Funktionals im Raum $C_0(\R_+,\R^n;\nu)$ (Folgerung \ref{FolgerungRieszC0nu}):
$$p(t) = \int_t^\infty \big[ \varphi^T_x\big(s,x_*(s),u_*(s)\big)p(s)-\lambda_0\omega(s) f_x\big(s,x_*(s),u_*(s)\big)\big] ds, \quad
  p(0) = l_0.$$
Also besitzt $p(\cdot)$ auf $\R_+$ die verallgemeinerte Ableitung
$$\dot{p}(t) = -\varphi_x^T\big(t,x_*(t),u_*(t)\big) p(t) + \lambda_0\omega(t) f_x\big(t,x_*(t),u_*(t)\big).$$
Aus der Festlegung der Funktion $p(\cdot)$ erhalten wir die Beziehung
$$\|p(t)\|^2\nu^{-1}(t) = \Big\langle \int_t^\infty \nu(s)d\mu(s),\int_t^\infty \nu(s)d\mu(s) \Big\rangle \cdot \nu^{-1}(t)
  \leq \|\mu\|^2 \nu(t).$$
Dies zeigt $p(\cdot) \in L_2(\R_+,\R^n;\nu^{-1})$ und ferner die erste ``nat"urliche'' Transversalit"atsbedingung in (\ref{WMP2}).
Au"serdem gilt
$$\big| \langle p(t), x(t) \rangle \big|^2 \leq \|p(t)\|^2 \nu^{-1}(t) \cdot \|x(t)\|^2 \nu(t).$$
Mit der eben nachgewiesenen Transversalit"atsbedingung und mit Lemma \ref{LemmaGrenzwert2} ergibt sich die zweite Bedingung in (\ref{WMP2}).
Damit sind $p(\cdot) \in L_2(\R_+,\R^n;\nu^{-1})$, sowie (\ref{WMP1}) und (\ref{WMP2}) gezeigt. \\[2mm]
Gem"a"s (\ref{SatzWMPLMR2}) gilt f"ur alle $v(\cdot) \in L_\infty(\R_+,U)$ die Ungleichung
$$\int_0^\infty \big\langle H_u\big(t,x_*(t),u_*(t),p(t),\lambda_0\big),v(t)-u_*(t) \big\rangle \, dt \leq 0.$$
Daraus folgt abschlie"send via Standardtechniken f"ur Lebesguesche Punkte die Variationsungleichung (\ref{WMP3}).
Der Beweis von Theorem \ref{SatzWMP} ist abgeschlossen. \hfill $\blacksquare$
       \subsection{Zur normalen Form der notwendigen Optimalit"atsbedingungen} \label{SectionNormalform}
Theorem \ref{SatzWMP} kann auf die folgende modifizierte Version des Beispiels von Halkin \cite{Halkin} angewendet werden.
Dieses Beispiel zeigt,
dass f"ur einen optimalen Steuerungsprozess die notwendigen Bedingungen in Theorem \ref{SatzWMP} nur f"ur nichttriviale
Multiplikatoren $\big(\lambda_0,p(\cdot)\big)$ mit $\lambda_0=0$ erf"ullt sind.

\begin{beispiel} \label{ExampleHalkin} {\rm Wir betrachten mit einem zus"atzlichen Diskontierungsfaktor die Aufgabe
\begin{eqnarray*}
&& J\big(x(\cdot),u(\cdot)\big) = \int_0^\infty \omega(t) \big(u(t)-x(t)\big) \, dt \to \sup, \\
&& \dot{x}(t) = u^2(t)+x(t), \quad x(0)=0, \quad  u(t) \in [0,1], \quad \varrho \in (0,1).
\end{eqnarray*}
In diesem Beispiel geh"ort $\big(x_*(t),u_*(t)\big) \equiv (0,0)$ zu $W^1_2(\R_+,\R;\nu) \times L_\infty(\R_+,U)$
f"ur jedes Gewicht $\nu(t)=e^{-at}$ mit $a>0$ und liefert das globale Maximum:
Denn angenommen, es gibt einen zul"assigen Steuerungsprozess $\big(x(\cdot),u(\cdot)\big)$ mit $x(t) \not\equiv 0$.
Insbesondere muss dann $\big|J\big(x(\cdot),u(\cdot)\big)\big| < \infty$ gelten.
Dann existiert ein $\tau >0$ mit $x(\tau)>0$.
Durch direktes Nachrechnen erhalten wir $x(t) \geq x(\tau)e^{t-\tau}$ f"ur alle $t \geq \tau$, und es folgt
$$J\big(x(\cdot),u(\cdot)\big) \leq \int_0^\tau e^{-\varrho t} \big(u(t)-x(t)\big) dt
              + \int_\tau^\infty e^{-\varrho t} \big(1-x(\tau)e^{t-\tau}\big) dt = - \infty.$$
Damit besteht die Menge der zul"assigen Steuerungsprozesse $\mathscr{A}_{\rm adm}$ nur aus $\big(x_*(\cdot),u_*(\cdot)\big)$. \\
In diesem Beispiel sind bei Wahl des Gewichtes $\nu(t)=e^{-at}$ mit $0<a<\varrho$ die Voraussetzungen (A$_1$) und (A$_2$) auf $U_{\gamma,\nu}$ erf"ullt.
Ferner ist $\varphi\big(t,x_*(t),u_*(t)\big) \equiv 0$.
Demzufolge gilt die Kontraktionsbedingung (\ref{WMPBedingung}).
Wenden wir Theorem \ref{SatzWMP} auf das modifizierte Halkin-Beispiel an, so erhalten wir die Adjungierte
$$p(t)=\bigg(p(0)-\frac{\lambda_0}{1-\varrho}\bigg)e^{-t}+\frac{\lambda_0}{1-\varrho}e^{(1-\varrho)t}.$$
Da nun $\varrho \in (0,1)$ gilt,
ist die erste Transversalit"atsbedingung in (\ref{WMP2}) nur dann erf"ullt,
wenn $\lambda_0=0$ ist.
Mit $p(0)<0$ gelten dann f"ur $\big(x_*(t),u_*(t)\big) \equiv (0,0)$ alle Optimalit"atsbedingungen in Theorem \ref{SatzWMP}. \hfill $\square$}
\end{beispiel}

\begin{bemerkung}{\rm
Im Beispiel \ref{ExampleHalkin} ist $\big(x_*(\cdot),u_*(\cdot)\big)$ der einzige zul"assige Steuerungsprozess.
In diesem Sonderfall kann man $\big(x_*(\cdot),u_*(\cdot)\big)$ auf der Menge der zul"assigen Steuerungsprozesse nicht variieren.
Es stellt sich damit die Frage,
ob die Existenz der nichttrivialen Multiplikatoren in Theorem \ref{SatzWMP} tats"achlich vorliegt. \\
Im Beweis des Extremalprinzips (Theorem \ref{SatzExtremalprinzipSchwach}) bezieht sich die Anwendung des Trennungssatzes auf die konvexe Menge
$\mathscr{C}$,
die nicht ausschlie"slich durch zul"assige Elemente des Steuerungsproblems (\ref{WMPAufgabe1})--(\ref{WMPAufgabe4}) gebildet werden.
Da weiterhin das Paar $\big(x_*(\cdot),u_*(\cdot)\big)$ der Menge $\mathscr{A}_{\rm Lip}$ angeh"ort und die Kontraktionsbedingung (\ref{WMPBedingung}) erf"ullt,
ergibt sich im Beweis des Extremalprinzips,
dass die Menge $\mathscr{C}$ ein nichtleeres Inneres besitzt.
Es folgt ferner die Existenz der nichttrivialen trennenden Hyperebene und demnach insbesondere die Existenz der
nichttrivialen Multiplikatoren in Theorem \ref{SatzWMP}. \hfill $\square$}
\end{bemerkung}

In diesem Abschnitt zeigen wir,
dass folgende ``Stetigkeitsbedingung'' (S) hinreichend f"ur die Normalform des Theorems \ref{SatzWMP} in der Aufgabe (\ref{WMPAufgabe1})--(\ref{WMPAufgabe3})
mit freiem rechten Endpunkt ist.

\begin{enumerate}
\item[(S)] Es existieren ein $T\geq 0$, eine Zahl $\varrho(T)>0$ und ein $\mu_T(\cdot) \in L_2(\R_+,\R_+;\nu)$ derart,
           dass f"ur alle $\zeta_T$ mit $\|\zeta_T-x_*(T)\| \leq \varrho(T)$ das System
           $\dot{x}(t) = \varphi\big(t,x(t),u_*(t)\big)$
           mit Anfangsbedingung $x(T)=\zeta_T$ eine L"osung $x(t;\zeta_T)$ auf $[T,\infty)$ besitzt und folgende Ungleichung gilt
           $$\| x(t;\zeta_T)-x_*(t)\| \leq \|\zeta_T-x_*(T)\|\mu_T(t) \quad \mbox{ f"ur alle } t \geq T \geq 0.$$
\end{enumerate}

\begin{theorem} \label{SatzNormalWMP} \index{Schwaches Optimalit"atsprinzip!normal@--, normale Form}
Sei $\big(x_*(\cdot),u_*(\cdot)\big) \in \mathscr{A}_{\rm adm} \cap \mathscr{A}_{\rm Lip}$ und seien (\ref{WMPBedingung}), (S) erf"ullt. 
Ist $\big(x_*(\cdot),u_*(\cdot)\big)$ ein schwaches lokales Minimum der Aufgabe (\ref{WMPAufgabe1})--(\ref{WMPAufgabe3}),
dann ist Theorem \ref{SatzWMP} mit $\lambda_0=1$ erf"ullt.
Ferner besitzt die Adjungierte $p(\cdot)$ die Darstellung\index{Adjungierte!eindeutig@--, eindeutige Darstellung}
\begin{equation} \label{NormalWMP1}
p(t)= -Z_*(t) \int_t^\infty \omega(s) Z^{-1}_*(s) f_x\big(s,x_*(s),u_*(s)\big) \, ds.
\end{equation}
Dabei ist $Z_*(t)$ die in $t=0$ normalisierte Fundamentalmatrix des linearen Systems
$$\dot{z}(t)=-\varphi^T_x\big(t,x_*(t),u_*(t)\big) z(t).$$
\end{theorem}
Die Darstellung (\ref{NormalWMP1}) stimmt
(bis auf das Vorzeichen, das sich durch die Minimierung statt einer Maximierung des Zielfunktionals ergibt)
mit der Formel in den Arbeiten von Aseev \& Kryazhimskii und Aseev \& Veliov \cite{AseKry,AseVel,AseVel2,AseVel3} "uberein.
Im Gegensatz zu diesen Arbeiten ist Theorem \ref{SatzNormalWMP} unter den Voraussetzungen
(A$_0$)--(A$_2$), (\ref{WMPBedingung}), (S) erf"ullt und charakterisiert schwache lokale Minimalstellen. \\[2mm]
Bevor wir Theorem \ref{SatzNormalWMP} beweisen,
wollen wir die Voraussetzung (S) und die Darstellungsformel \ref{NormalWMP1} im Differentialspiel \ref{BeispielDockner} demonstrieren.

\begin{beispiel} \label{BeispielFischereispiel}
{\rm Im Differentialspiel\index{Differentialspiel}\index{Nash-Gleichgewicht}\index{Fischereimodell}
im Beispiel \ref{BeispielDockner} mit unendlichem Zeithorizont "uberf"uhrten wir den Ansatz eines Nash-Gleichgewichtes in die gekoppelten Aufgaben
\begin{eqnarray*}
&& J_i\big(z(\cdot),u_1(\cdot),u_2(\cdot)\big)
        =\int_0^\infty e^{-\varrho t}\bigg(\frac{1}{u_1(t)+u_2(t)}-c_i\bigg)u_i(t) \, dt \to \sup, \\
&& \dot{z}(t)=-r z(t)+\alpha-u_1(t)-u_2(t), \quad z(0)=\ln x_0>0, \\
&& u_i>0, \quad \alpha,c_i,p,r,\varrho >0, \quad \alpha> \frac{1}{c_1+c_2}, \quad i=1,2.
\end{eqnarray*}
Die lineare Dynamik mit Wachstumskoeffizient $-r$ f"uhrt zur Stetigkeitsbedingung
$$\| z(t;\zeta_T)-z_*(t)\| = \|\zeta_T-z_*(T)\|e^{-r(t-T)} = \|\zeta_T-z_*(T)\|\mu_T(t) \quad \mbox{ f"ur alle } t \geq T \geq 0.$$
Dabei geh"ort die Funktion $\mu_T(t)=e^{rT}e^{-rt}$ stets dem Raum $L_2(\R_+,\R_+;\nu)$ mit Gewicht $\nu(t)=e^{-at}$, $a>0$, an.
In beiden Aufgaben ist stets $f_x\big(t,x,u_1,u_2\big) = 0$ und es liefert (\ref{NormalWMP1}) unmittelbar $p_i(t) \equiv 0$ f"ur $i=1,2$. \hfill $\square$}
\end{beispiel}

{\bf Beweis von Theorem \ref{SatzNormalWMP}}
Es seien $Y_*(t)$ bzw. $Z_*(t)$ die in $t=0$ normalisierten Fundamentalmatrizen der homogenen Systeme
$$\dot{y}(t)=\varphi_x\big(t,x_*(t),u_*(t)\big) y(t), \qquad \dot{z}(t)=-\varphi_x^T\big(t,x_*(t),u_*(t)\big) z(t).$$
Nach Voraussetzung (S) existiert die L"osung $x_\alpha(\cdot)$ der Gleichung
$$x(t)= x(T) + \int_T^t \varphi\big(s,x(s),u_*(s)\big) \, ds, \quad x(T)=x_*(T) + \alpha \xi, \quad \|\xi\|=1,$$
f"ur alle $\alpha \in [0,\varrho(T)]$.
Gem"a"s dem Satz "uber die Abh"angigkeit einer L"osung von den Anfangsdaten erhalten wir daher
$$\frac{x_\alpha(t)-x_*(t)}{\alpha}
   = \xi + \int_T^t \frac{\varphi_x\big(s,x_*(s),u_*(s)\big)\big(x_\alpha(s)-x_*(s)\big)+ o(\alpha,s)}{\alpha} \, ds.$$
Dabei ist $o(\alpha,t)/\alpha \to 0$ f"ur $\alpha \to 0^+$, und das gleichm"a"sig auf jedem endlichen Intervall $[T,K]$, $K>T$.
Der Grenz"ubergang $\alpha \to 0^+$ liefert f"ur jedes feste $t\geq T$:
$$y(t):=\lim_{\alpha \to 0^+} \frac{x_\alpha(t)-x_*(t)}{\alpha}
   = \xi + \int_T^t \varphi_x\big(s,x_*(s),u_*(s)\big) y(s) \, ds = Y_*(t) Y^{-1}_*(T) \xi.$$
Wir setzen $y(t)$ durch $y(t)=Y_*(t) Y^{-1}_*(T) \xi$ auf $\R_+$ fort.
Wegen Voraussetzung (A$_2$) folgt mit der Gronwallschen Ungleichung $\|y(t)\| \leq C \cdot e^{C_0 t}$ auf $\R_+$.
Ferner erhalten wir mit (S) f"ur alle $t \geq T$:
$$\|y(t)\| = \lim_{\alpha \to 0^+} \frac{\|x_\alpha(t)-x_*(t)\|}{\alpha} \leq \|\xi\| \cdot \mu_T(t).$$
Wegen $\mu_T(\cdot) \in L_2(\R_+,\R_+;\nu)$ zeigt dies $y(\cdot) \in W^1_2(\R_+,\R^n;\nu)$. \\[2mm]
Angenommen, es ist $\lambda_0=0$. 
Dann erf"ullt die Adjungierte $p(\cdot)$ nach Theorem \ref{SatzWMP} die Gleichung (\ref{WMP1}).
Wegen $\frac{d}{dt} \langle p(t),y(t) \rangle =0$ folgt auf $\R_+$:
$$\langle p(t),y(t) \rangle = \langle p(0),Y^{-1}_*(T)\xi \rangle.$$
Aufgrund der Transversalit"atsbedingung (\ref{WMP3}) gilt $\lim\limits_{t \to \infty} \langle p(t),y(t) \rangle = 0$.
Da $\xi \in \R^n$, $\|\xi\|=1$, beliebig war,
erhalten wir $\big(\lambda_0,p(\cdot)\big) =0$ im Widerspruch zu Theorem \ref{SatzWMP}. \\[2mm]
Es sei $\lambda_0=1$.
Dann gilt f"ur jedes $T \in \R_+$ (vgl. Aseev \& Kryazhimskii \cite{AseKry}):
$$p(t)=Z_*(t)\bigg(Z^{-1}_*(T) p(T) + \int_T^t \omega(s) Z_*^{-1}(s)f_x\big(s,x_*(s),u_*(s)\big) ds \bigg).$$
Betrachten wir nun $\langle p(t), y(t) \rangle$ und verwenden $Z^{-1}_*(t)=Y^T_*(t)$, so ergibt sich
$$\langle p(t), y(t) \rangle = 
  \Big\langle  p(T) + Z_*(T)\int_T^t \omega(s) Z_*^{-1}(s)f_x\big(s,x_*(s),u_*(s)\big) ds , \xi \Big\rangle.$$
Mit (\ref{WMP3}) liefert der Grenz"ubergang $t \to \infty$ die Darstellungsformel (\ref{NormalWMP1}). \hfill $\blacksquare$ \\[2mm]
Eine L"osung der adjungierten Gleichung (\ref{WMP1}) zu den ``nat"urlichen'' Transversalit"atsbedingungen (\ref{WMP3})
in Theorem \ref{SatzWMP} muss keineswegs eindeutig bestimmt sein.
Sondern es kann durch die notwendigen Optimalit"atsbedingungen eine ganze Schar von Kandidaten gefunden werden.
Zieht man aber in die Betrachtung die eindeutige Darstellung der Adjungierten nach (\ref{NormalWMP1}) hinzu,
so l"asst sich diese Schar h"aufig auf wenige oder einen einzigen Kandidaten reduzieren.
       \subsection{Hinreichende Bedingungen nach Mangasarian} \label{AbschnittMangasarianWMP}
Die\index{hinreichende Bedingungen!Mangasarian@-- nach Mangasarian} Darstellung der hinreichenden Bedingungen f"ur die Aufgabe
(\ref{WMPAufgabe1})--(\ref{WMPAufgabe3}) ist Seierstad \& Syds\ae ter \cite{Seierstad} entnommen.
Der wesentliche Unterschied zur Formulierung notwendiger Bedingungen in Form von Theorem \ref{SatzWMP} ist,
dass nicht die Existenz der Multiplikatoren,
sondern dass die G"ultigkeit der Bedingungen (\ref{WMP1})--(\ref{WMP3}) in Kombination mit Konkavit"atseigenschaften der Pontrjagin-Funktion $H$
bereits hinreichende Bedingungen liefern. \\
Dementsprechend weisen wir explizit darauf hin,
dass in der Herleitung hinreichender Bedingungen f"ur ein schwaches lokales Minimum die Eigenschaften der Aufgabe (\ref{WMPAufgabe1})--(\ref{WMPAufgabe3}) 
lediglich auf der gleichm"a"sigen Umgebung
$$U_\gamma= \{ (t,x,u) \in \R_+ \times \R^n \times \R^m\,|\, \|x-x(t)\| \leq \gamma, \|u-u(t)\| \leq \gamma\}$$
von Interesse sind und au"serdem die Kontraktionsbedingung (\ref{WMPBedingung}) nicht einflie"st.

\begin{theorem} \label{SatzHBWMP}
In der Aufgabe (\ref{WMPAufgabe1})--(\ref{WMPAufgabe3}) sei
$\big(x_*(\cdot),u_*(\cdot)\big)$ ein zul"assiger Steuerungsprozess,
f"ur den die Abbildungen $f(t,x,u)$ und $\varphi(t,x,u)$ auf der Menge
$$U_\gamma= \{ (t,x,u) \in \R_+ \times \R^n \times \R^m\,|\, \|x-x_*(t)\| \leq \gamma, \|u-u_*(t)\| \leq \gamma\}$$
stetig differenzierbar sind.
Ferne gelte:
\begin{enumerate}
\item[(a)] Das Tripel $\big(x_*(\cdot),u_*(\cdot),p(\cdot)\big)$
           erf"ullt (\ref{WMP1})--(\ref{WMP3}) mit $\lambda_0=1$ in Theorem \ref{SatzWMP}.        
\item[(b)] F"ur jedes $t \in \R_+$ ist die Funktion $H\big(t,x,u,p(t),1\big)$
           konkav in $(x,u)$ auf $U_\gamma$.
\end{enumerate}
Dann ist $\big(x_*(\cdot),u_*(\cdot)\big)$ ein schwaches lokales Minimum der Aufgabe (\ref{WMPAufgabe1})--(\ref{WMPAufgabe3}).
\end{theorem}

{\bf Beweis} Wir erhalten zu $T \in \R_+$ die Beziehung
\begin{eqnarray*}
    \Delta(T)
&=& \int_0^T \omega(t)\big[f\big(t,x(t),u(t)\big)-f\big(t,x_*(t),u_*(t)\big)\big] \, dt \\
&=& \int_0^T \big[H\big(t,x_*(t),u_*(t),p(t),1\big)-H\big(t,x(t),u(t),p(t),1\big)\big] \, dt \\
& & + \int_0^T \langle p(t), \dot{x}(t)-\dot{x}_*(t) \rangle \, dt.
\end{eqnarray*}
Aus den elementaren Eigenschaften konkaver Funktionen folgt
\begin{eqnarray*}
\lefteqn{H\big(t,x_*(t),u_*(t),p(t),1\big)-H\big(t,x(t),u(t),p(t),1\big)\geq} \\
&& -H_x\big(t,x_*(t),u_*(t),p(t),1\big)\big(x(t)-x_*(t)\big)
       -H_u\big(t,x_*(t),u_*(t),p(t),1\big)\big(u(t)-u_*(t)\big).
\end{eqnarray*}
Damit, sowie mit der adjungierten Gleichung (\ref{WMP1}) und der Variationsungleichung (\ref{WMP3}) gilt
$$\Delta(T) \geq  \int_0^T \langle \dot{p}(t),x(t)-x_*(t)\rangle + \langle p(t), \dot{x}(t)-\dot{x}_*(t) \rangle dt
  = \langle p(T),x(T)-x_*(T)\rangle.$$
Es folgt abschlie"send mit den ``nat"urlichen'' Transversalit"atsbedingungen in (\ref{WMP2}) die Beziehung
$$\lim_{T \to \infty} \Delta(T) \geq \lim_{T \to \infty} \langle p(T),x(T)-x_*(T)\rangle=0$$
f"ur alle zul"assigen $\big(x(\cdot),u(\cdot)\big)$ mit 
$\|x(\cdot)-x_*(\cdot)\|_\infty, \|u(\cdot)-u_*(\cdot)\|_{L_\infty} \leq \gamma$. \hfill $\blacksquare$

\begin{beispiel}\label{BeispielRegler2} {\rm Im linear-quadratischen Regler nach Pickenhain \& Wenzke \cite{PickiWenzke}
\begin{eqnarray*}
&& J\big(x(\cdot),u(\cdot)\big) = \int_0^\infty e^{-2t} \cdot \frac{1}{2}\big( x^2(t)+u^2(t)\big) \, dt \to \inf, \\
&& \dot{x}(t) = 2 x(t)+u(t), \qquad x(0)=2, \qquad u(t) \in \R
\end{eqnarray*}
gen"ugt der Steuerungsprozess und die Adjungierte
$$x_*(t)=2e^{(1-\sqrt{2})t}, \quad u_*(t)=-2(1+\sqrt{2})e^{(1-\sqrt{2})t}, \quad p(t)=e^{-2t}u_*(t)$$
den notwendigen Bedingungen (\ref{WMP1})--(\ref{WMP3}) in Theorem \ref{SatzWMP}.
Weiterhin ist $H\big(t,x,u,p(t),1\big)$ konkav bez"uglich $(x,u)$ auf $U_\gamma$ f"ur jedes $\gamma>0$.
Damit ist $\big(x_*(\cdot),u_*(\cdot)\big)$ ein schwaches lokales Minimum dieser Aufgabe. \hfill $\square$}
\end{beispiel}    

       \rhead[]{Aufgaben mit Zustandsbeschr"ankungen \hspace*{1mm} \thepage}
       \section{Die Aufgabe mit Zustandsbeschr"ankungen} \label{AbschnittZustandWMP}
Wir werden unser Vorgehen auf die Aufgabe (\ref{WMPAufgabe1})--(\ref{WMPAufgabe4}) mit Zustandsbeschr"ankungen anwenden.
Im Weiteren bezeichne $H: \R \times \R^n \times \R^m \times \R^n \times \R \to \R$ wieder die Pontrjagin-Funktion
$$H(t,x,u,p,\lambda_0) = -\lambda_0 \omega(t)f(t,x,u) + \langle p, \varphi(t,x,u) \rangle.$$

\begin{theorem}\label{SatzWMPAufgabeZB} \index{Schwaches Optimalit"atsprinzip}
Sei $\big(x_*(\cdot),u_*(\cdot)\big) \in \mathscr{A}_{\rm adm} \cap \mathscr{A}_{\rm Lip}$.
Weiterhin sei die Kontraktionsbedingung
\begin{equation} \label{WMPBedingungZB}
\sup_{t \in \R_+} \int_0^t \frac{\nu(t)}{\nu(s)} \big\|\varphi_x\big(s,x_*(s),u_*(s)\big)\big\| \, ds < 1
\end{equation}
erf"ullt und es m"ogen f"ur alle $x(\cdot) \in C_0(\R_+,\R^n;\nu)$ mit $\|x(\cdot)-x_*(\cdot)\|_{\infty,\nu} < \gamma$ die folgenden Grenzwerte existieren:
\begin{equation} \label{WMPBedingungZB2}
\lim_{t \to \infty} g_{jx}\big(t,x(t)\big), \qquad j=1,...,l.
\end{equation}
Ist $\big(x_*(\cdot),u_*(\cdot)\big)$ ein schwaches lokales Minimum der Aufgabe (\ref{WMPAufgabe1})--(\ref{WMPAufgabe4}),
dann existieren eine Zahl $\lambda_0 \geq 0$, ein Vektor $l_0 \in \R^n$, eine Funktion $p(\cdot):\R_+ \to \R^n$
und auf den Mengen
$$T_j=\big\{t \in \overline{\R}_+ =\R_+ \cup \{\infty\} \,\big|\, g_j\big(t,x_*(t)\big)=0\big\}, \quad j=1,...,l,$$
konzentrierte nichtnegative regul"are Borelsche Ma"se $\mu_j$ endlicher Totalvariation
(wobei s"amtliche Gr"o"sen nicht gleichzeitig verschwinden) derart, dass
\begin{enumerate}
\item[(a)] die Vektorfunktion $p(\cdot)$ von beschr"ankter Variation ist, der adjungierten Gleichung\index{adjungierte Gleichung}
           \begin{equation}\label{SatzWMPAufgabeZB1}
           p(t)= \int_t^\infty H_x\big(s,x_*(s),u_*(s),p(s),\lambda_0\big) \, ds
                  -\sum_{j=1}^l \int_t^\infty \nu(s) g_{jx}\big(s,x_*(s)\big) \, d\mu_j(s)
           \end{equation}
           gen"ugt und die ``nat"urlichen'' Transversalit"atsbedingungen\index{Transversalit"atsbedingungen!nat@--, nat"urliche}
           \begin{equation}\label{SatzWMPAufgabeZB2}
           \lim_{t \to \infty} \|p(t)\|^2\nu^{-1}(t) = 0, \qquad
           \lim_{t \to \infty} \langle p(t),x(t) \rangle = 0 \quad \forall \; x(\cdot) \in W^1_2(\R_+,\R^n;\nu)
           \end{equation}
           erf"ullt;
\item[(b)] in fast allen Punkten $t \in \R_+$ und f"ur alle $u \in U$ die Variationsungleichung
           \begin{equation}\label{SatzWMPAufgabeZB3}
           \big\langle H_u\big(t,x_*(t),u_*(t),p(t),\lambda_0\big),\big(u-u_*(t)\big) \big\rangle\leq 0
           \end{equation}
           gilt.
\end{enumerate}
\end{theorem}

\begin{beispiel}[Abbau einer nicht erneuerbaren Ressource] \label{ExampleRessource}\index{Ressourcenabbau}
{\rm Wir betrachten die Aufgabe
\begin{eqnarray}
&& \label{Ressource1} J\big(x(\cdot),y(\cdot),u(\cdot)\big)
   =\int_0^\infty e^{-\varrho t}\big[pf\big(u(t)\big)-ry(t)-qu(t)\big] \, dt \to \sup, \\
&& \label{Ressource2} \dot{x}(t) = -u(t),\quad \dot{y}(t)=cf\big(u(t)\big), \quad x(0)=x_0>0,\quad y(0)=y_0\geq 0, \\
&& \label{Ressource3} x(t) \geq 0, \qquad u(t) \geq 0, \qquad b,c,\varrho, q,r >0, \qquad \varrho - rc>0.
\end{eqnarray}
Die Funktion $f$ sei zweimal stetig differenzierbar, $f'>0$, $f'(0)<\infty$, $f''<0$ und es sei
$f'(u) \to 0$ f"ur $u \to \infty$.
In der vorliegenden Formulierung der Aufgabe wurde im Vergleich zu Seierstad \& Syds\ae ter \cite{Seierstad} die
Restriktion $\liminf\limits_{t \to \infty} x(t) \geq 0$ durch die Zustandsbeschr"ankung $x(t) \geq 0$ in (\ref{Ressource3})
ersetzt. \\[2mm]
"Okonomische Interpretation:
$x(t)$ bezeichnet die Menge einer nat"urlichen Ressource und $u(t)$ ist die industrielle Abbaurate dieser Ressource.
Auf Basis der Ressource werden G"uter mit der Produktionsrate $f\big(u(t)\big)$ hergestellt.
Die Kosten der Herstellung einer Produktionseinheit ist $q$ und der Preis einer G"utereinheit am Markt betr"agt $p$.
Bei der Herstellung der G"uter entstehen proportional zur Produktion Abf"alle,
deren Gesamtmenge durch $y(t)$ beschrieben wird.
Die Kosten der Beseitigung der negativen Auswirkungen der Abfallmenge sind $ry(t)$.
Im Weiteren gehen wir von dem Preis $p=1$ aus. \\[2mm]
Wir pr"ufen die Voraussetzungen an die Aufgabe:
Mit der Festlegung $L(t)=C_0e^{at}$ und $\nu(t)=e^{-at}$ mit $0<a<\varrho$ gen"ugt die Aufgabe den Voraussetzungen (A$_1$)--(A$_3$).
Weiterhin sind die Dynamiken
$$\varphi_1(t,x,y,u)=-u, \qquad \varphi_2(t,x,y,u)=cf(u)$$
unabh"angig von den Zustandsvariablen $x,y$ und die Bedingungen (\ref{WMPBedingungZB}), (\ref{WMPBedingungZB2}) sind erf"ullt.
Wir stellen die Optimalit"atsbedingungen von Theorem \ref{SatzWMPAufgabeZB} mit $\lambda_0=1$ auf:
\begin{enumerate}
\item[(a)] Die Pontrjagin-Funktion der Aufgabe (\ref{Ressource1})--(\ref{Ressource3}) lautet
           $$H(t,x,y,u,p_1,p_2,1) = p_1 (-u)+p_2 cf(u) + e^{-\varrho t}[f(u)-ry-qu].$$
\item[(b)] Die Adjungierten gen"ugen den Gleichungen
           $$p_1(t)=\int_t^\infty e^{-as} \, d\mu(s), \qquad \dot{p}_2(t)=r e^{-\varrho t} \Rightarrow p_2(t)=-\frac{r}{\varrho}e^{-\varrho t} + K.$$
           Das auf der Menge $T=\{t \in \overline{\R}_+ \,|\, x_*(t)=0\}$ konzentrierte regul"are Ma"s $\mu$ ist nichtnegativ.
           Daher ist $p_1(t) \geq 0$ "uber $\R_+$ und monoton fallend.
           Ferner erhalten wir $K=0$ aus den Transversalit"atsbedingungen bez"uglich dem Zustand $y$.
\item[(c)] Die Maximumbedingung k"onnen wir auf folgende Aufgabe reduzieren
           $$\max_{u \geq 0} \Big[ -p_1(t) u +c p_2(t) f(u) + e^{-\varrho t}[f(u)-qu]\Big].$$
           Das Einsetzen der Darstellung f"ur $p_2(t)$ liefert weiterhin mit $d=(\varrho -rc)/\varrho$:
           $$\max_{u \geq 0} \Big( d f(u)e^{-\varrho t}-u\big(p_1(t)+ qe^{-\varrho t}\big)\Big).$$
\end{enumerate}
Die Reduktion der Maximumbedingung f"uhrt f"ur festes $t$ zu der Funktion
$$g(u)= d f(u)e^{-\varrho t}-u\big(p_1(t)+ qe^{-\varrho t}\big).$$
Diese Funktion ist zweimal stetig differenzierbar und es gilt
$$g'(u)= \big(d f'(u)-q\big)e^{-\varrho t}-p_1(t), \quad g''(u)=df''(u)e^{-\varrho t}, \quad d=\frac{\varrho -rc}{\varrho}>0.$$
Daher ist $g$ streng konkav und besitzt auf der Menge $U=\{u\geq 0\}$ ein Maximum,
da $f'(u)>0$ und $f'(u) \to 0$ f"ur $u \to \infty$ gelten.
Wir diskutieren drei F"alle:
\begin{enumerate}
\item[(A)] $df'(0)\leq q$: In diesem Fall ist $g'(0) \leq 0$ und man erh"alt
           $$u_*(t) \equiv 0, \quad x_*(t) \equiv x_0, \quad
             y_*(t)=y_0 + cf(0)t, \quad p_1(t) \equiv 0, \quad p_2(t)=-\frac{r}{\varrho}e^{-\varrho t}.$$
           Da die Zustandsbeschr"ankung nichtaktiv ist, gelten die Voraussetzungen und Optimalit"atsbedingungen aus Theorem \ref{SatzWMPAufgabeZB} 
           f"ur $\nu(t)=e^{-at}$ mit $0<a<\varrho$.
\item[(B)] $df'(0)> q$ und $p_1(0)=0$:
           Aus $g'(u)=\big(d f'(u)-q\big)e^{-\varrho t}=0$ erhalten wir die optimale Strategie $u_*(t)=u_0>0$ f"ur alle $t \in \R_+$.
           Also gilt $x_*(t)=x_0-u_0t$ auf $\R_+$, was der Zustandsbeschr"ankung widerspricht.
\item[(C)] $df'(0)> q$ und $p_1(0)>0$:
           Wegen $p_1(0)>0$ wird die Ressource vollst"andig abgebaut.
           Andernfalls w"are $p_1(t)=p_1(0)>0$ "uber $\R_+$, was (\ref{SatzWMPAufgabeZB2}) widerspricht.
           Da die Ressource vollst"andig abgebaut wird,
           gibt es ein $t'>0$ mit $x_*(t)>0$ f"ur $t \in [0,t')$ und $x_*(t)=0$ f"ur $t\geq t'$.
           Demnach folgt unmittelbar $u_*(t)=0$ f"ur $t\geq t'$. \\
           F"ur $t\geq t'$ ist $p_1(\cdot)$ monoton fallend.
           Ferner erhalten wir f"ur $t \in \R_+$ die Beziehung
           $$g'(u)=0 \qquad\Rightarrow\qquad f'\big(u(t)\big)=\frac{1}{d}(q+p_1(t)e^{\varrho t}).$$
           W"urde demnach die Adjungierte $p_1(\cdot)$ f"ur $t \geq t'$ eine Unstetigkeitstelle besitzen,
           dann folgt aus der Monotonie von $p_1(\cdot)$, dass die Abbaurate sich wieder sprunghaft vergr"o"sert.
           Diese Steuerung f"uhrt zu einem erneuten Abbau der Ressource, obwohl diese bereits vollst"andig aufgebraucht ist.
           Daher ist die Adjungierte stetig. \\
           F"ur die Adjungierte erhalten wir damit
           $$p_1(t) = \big(df'(0)-q\big)e^{-\varrho t'} \mbox{ f"ur } t \leq t', \qquad
             p_1(t) = \big(df'(0)-q\big)e^{-\varrho t} \mbox{ f"ur } t \geq t'.$$
           Wir zeigen noch, dass der Zeitpunkt $t'$ existiert und eindeutig ist:
           Durch
           $$f'\big(u_\tau(t)\big)=\frac{1}{d}(q+p_1(0)e^{\varrho t})=\frac{1}{d}(q+[df'(0)-q]e^{\varrho(t-\tau)}), \quad t \in [0,\tau],$$
           und $u_\tau(t)=0$ f"ur $t \geq \tau$ wird wegen $f'\big(u_\tau(\tau)\big)=f'(0)$ eine Familie $u_\tau(\cdot)$ stetiger Funktionen definiert.
           Dabei gilt $f'\big(u_\tau(t)\big) < f'\big(u_s(t)\big)$, d.\,h. $u_\tau(t) > u_s(t)$, f"ur alle $t \in [0,\tau]$ und $\tau>s$.
           Damit ist die Familie
           $$U(\tau):= \int_0^\infty u_\tau(t) \, dt$$
           streng monoton wachsend und es gelten $U(0)=0$, $U(\tau) \to \infty$ f"ur $\tau \to \infty$.
           Der Parameter $t'$ ergibt sich dann aus der Bedingung $U(t')=x_0$. \hfill $\square$
\end{enumerate}}
\end{beispiel}
       \subsection{Der Nachweis der notwendigen Optimalit"atsbedingungen} \label{AbschnittBeweisWMPZB}
Sei $\big(x_*(\cdot),u_*(\cdot)\big) \in \mathscr{A}_{\rm Lip}$.
Die Menge $U_{\gamma,\nu}$ ist wie folgt definiert:
$$U_{\gamma,\nu}= \{ (t,x,u) \in \R_+ \times \R^n \times \R^m\,|\, e^{-at}\|x-x(t)\| \leq \gamma, \|u-u(t)\| \leq \gamma\}.$$
Wir betrachten f"ur $\big(x(\cdot),u(\cdot)\big) \in C_0(\R_+,\R^n;\nu) \times L_\infty(\R_+,\R^m)$ die Abbildungen
\begin{eqnarray*}
J\big(x(\cdot),u(\cdot)\big) &=& \int_{t_0}^{t_1} \omega(t) f\big(t,x(t),u(t)\big) \, dt, \\
F\big(x(\cdot),u(\cdot)\big)(t) &=& x(t) -x(t_0) -\int_{t_0}^t \varphi\big(s,x(s),u(s)\big) \, ds, \quad t \in \R_+,\\
H_0\big(x(\cdot)\big) &=& x(0), \\
G_j\big(x(\cdot)\big)(t) &=& \nu(t) g_j\big(t,x(t)\big), \quad t \in \R_+, \quad j=1,...,l.
\end{eqnarray*}
Dabei fassen wir sie als Abbildungen zwischen folgenden Funktionenr"aumen auf:
\begin{eqnarray*}
J &:& C_0(\R_+,\R^n;\nu) \times L_\infty(\R_+,\R^m) \to \R, \\
F &:& C_0(\R_+,\R^n;\nu) \times L_\infty(\R_+,\R^m) \to C_0(\R_+,\R^n;\nu), \\
H_0 &:& C_0(\R_+,\R^n;\nu) \to \R^n, \\
G_j &:& C_0(\R_+,\R^n;\nu) \to C_{\lim}(\R_+,\R^n), \quad j=1,...,l.
\end{eqnarray*}
Wir setzen $\mathscr{F}=(F,H_0)$ und $G=(G_1,...,G_l)$.
Au"serdem f"uhren wir f"ur die Elemente $x(\cdot)= \big(x_1(\cdot),...,x_l(\cdot)\big)$ des Raumes $C_{\lim}(\R_+,\R^l)$ folgende Halbordnung ``$\preceq$'' ein:
$$x(\cdot) \preceq y(\cdot) \qquad\Leftrightarrow\qquad x_j(t) \leq y_j(t) \mbox{ f"ur alle } t \in \R_+, \; j=1,...,l.$$
Im Raum $C_{\lim}(\R_+,\R^l)$ bezeichnen wir mit $\mathscr{K}$ den nachstehenden konvexen, abgeschlossenen Kegel mit Spitze in Null:
$$\mathscr{K} =\{ x(\cdot) \in C_{\lim}(\R_+,\R^l) \,|\, x(\cdot) \preceq 0\}.$$
F"ur den Kegel $\mathscr{K}$ ist ${\rm int\,}\mathscr{K} \not= \emptyset$. \\[2mm]
Mit diesen Setzungen pr"ufen wir f"ur die Extremalaufgabe
\begin{equation} \label{ExtremalaufgabeWMPZB}
J\big(x(\cdot),u(\cdot)\big) \to \inf, \quad \mathscr{F}\big(x(\cdot),u(\cdot)\big)=0, \quad G\big(x(\cdot)\big) \in \mathscr{K},
\quad  u(\cdot) \in L_\infty(\R_+,U)
\end{equation}
die Voraussetzungen von Theorem \ref{SatzExtremalprinzipSchwach} im Punkt $\big(x_*(\cdot),u_*(\cdot)\big)$:

\begin{enumerate}
\item[(A$_2$)] Mit Verweis auf Abschnitt \ref{AbschnittBeweisWMP} ist nur noch die Abbildung $G$ zu diskutieren,
               deren stetige Fr\'echet-Differenzierbarkeit im Beispiel \ref{DiffAbbildung2} nachgewiesen ist.
               Dabei sind die Anforderungen an die Abbildungen $g_j$ nach Voraussetzung (A$_3$) erf"ullt und es ergeben sich aus (A$_3$)
               $$\lim_{t \to \infty} \nu(t) \big\|g\big(t,x(t)\big)\big\| \leq C_0 \lim_{t \to \infty} \nu(t)(1+\|x(t)\|) = 0$$
               und mit (\ref{WMPBedingungZB2}) die Existenz der Grenzwerte
               $$\lim_{t \to \infty} g_x\big(t,x(t)\big)$$
               f"ur alle $x(\cdot) \in C_0(\R_+,\R^n;\nu)$ mit $\|x(\cdot)-x_*(\cdot)\|_{\infty,\nu} < \gamma$.
\end{enumerate}

Zur Extremalaufgabe (\ref{ExtremalaufgabeWMPZB}) definieren wir auf
$$C_0(\R_+,\R^n;\nu) \times L_\infty(\R_+,\R^m) \times \R \times C_0^*(\R_+,\R^n;\nu) \times \R^n \times \big(C_{\lim}^*(\R_+,\R^n)\big)^l$$
die Lagrange-Funktion $\mathscr{L}=\mathscr{L}\big(x(\cdot),u(\cdot),\lambda_0,y^*,l_0,z_1^*,...,z_l^*\big)$,
$$\mathscr{L}= \lambda_0 J\big(x(\cdot),u(\cdot)\big)+ \big\langle y^*, F\big(x(\cdot),u(\cdot)\big) \big\rangle
                         +l_0^T H\big(x(\cdot)\big) + \sum_{j=1}^l \big\langle z_j^*, G_j\big(x(\cdot)\big) \big\rangle.$$
Ist $\big(x_*(\cdot),u_*(\cdot)\big)$ eine schwache lokale Minimalstelle der Aufgabe (\ref{ExtremalaufgabeWMPZB}),
dann existieren nach Theorem \ref{SatzExtremalprinzipSchwach}
nicht gleichzeitig verschwindende Lagrangesche Multiplikatoren $\lambda_0 \geq 0$, $y^* \in C_0^*(\R_+,\R^n;\nu)$, $l_0 \in \R^n$
und $z_j^* \in C_{\lim}^*(\R_+,\R^n)$ derart,
dass gelten:
\begin{enumerate}
\item[(a)] Die Lagrange-Funktion besitzt bez"uglich $x(\cdot)$ in $x_*(\cdot)$ einen station"aren Punkt, d.\,h.
          \begin{equation} \label{SatzWMPZBLMR1}
          \mathscr{L}_x\big(x_*(\cdot),u_*(\cdot),\lambda_0,y^*,l_0,z_1^*,...,z_l^*\big)=0;
          \end{equation}         
\item[(b)] Die Lagrange-Funktion erf"ullt bez"uglich $u(\cdot)$ in $u_*(\cdot)$ die Variationsungleichung
           \begin{equation} \label{SatzWMPZBLMR2}
           \big\langle \mathscr{L}_u\big(x_*(\cdot),u_*(\cdot),\lambda_0,y^*,l_0,z_1^*,...,z_l^*\big), u(\cdot)-u_*(\cdot) \big\rangle \geq 0
           \end{equation}
           f"ur alle $u(\cdot) \in L_\infty(\R_+,U)$;
\item[(c)] Die komplement"aren Schlupfbedingungen gelten, d.\,h.
           \begin{equation}\label{SatzWMPZBLMR3}
           0 = \big\langle z_i^*, G_i\big(x_*(\cdot)\big) \big\rangle, \quad
           \langle z_i^*,z(\cdot) \rangle \leq 0 \quad\mbox{f"ur alle } z(\cdot) \in \mathscr{K}, \quad i=1,...,l.
           \end{equation}
\end{enumerate}

Aufgrund (\ref{SatzWMPZBLMR1}) ist folgende Variationsgleichung f"ur alle $x(\cdot) \in C_0(\R_+,\R^n;\nu)$ erf"ullt: 
\begin{eqnarray}
0 &=& \lambda_0 \cdot \int_0^\infty \omega(t) \big\langle f_x\big(t,x_*(t),u_*(t)\big),x(t) \big\rangle\, dt + l_0^T x(0) \nonumber \\
  & & + \int_0^\infty \nu(t)\bigg[ x(t)-x(0) - \int_0^t \varphi_x\big(s,x_*(s),u_*(s)\big) x(s) \,ds \bigg]^T  d\mu(t) \nonumber \\
  & & \label{BeweisschlussWMPZB1}
      + \sum_{j=1}^l \int_0^\infty \big\langle \nu(t) g_{jx}\big(t,x_*(t)\big),x(t) \big\rangle \,d\mu_j(t).
\end{eqnarray}
Dabei ist nach Folgerung \ref{FolgerungRieszC0nu} $\mu$ ein signiertes regul"ares Borelsches Vektorma"s "uber $\R_+$ und
es sind nach Satz \ref{SatzRieszClimnu} $\mu_j$ positive regul"are Borelsche Ma"se "uber $\overline{\R}_+$.
Da die Ausdr"ucke
$$\lim_{t \to \infty} \nu(t) g_{jx}\big(t,x_*(t)\big)$$
f"ur alle $j=1,...,l$ im Unendlichen verschwinden,
gehen s"amtliche atomaren Anteile in $t=\infty$ der Ma"se $\mu_j$ in der Gleichung (\ref{BeweisschlussWMPZB1}) verloren.
Wir bemerken an dieser Stelle,
dass die Betrachtung $x(\cdot) \in C_{\lim}(\R_+,\R^n;\nu)$ diesen Umstand nicht behebt. \\[2mm]
Durch vertauschen der Integrationsreihenfolge im letzten Summanden in (\ref{BeweisschlussWMPZB1}) bringen wir diese Gleichung in die Form
\begin{eqnarray}
0 &=& \int_0^\infty \bigg[ \lambda_0 \omega(t)f_x\big(t,x_*(t),u_*(t)\big) - \varphi^T_x\big(t,x_*(t),u_*(t)\big) \int_t^\infty
      \nu(s) d\mu(s) \bigg]^T  x(t) dt \nonumber \\
  & & + \int_0^\infty \nu(t) [x(t)]^T d\mu(t) + \bigg[l_0^T x(0) - \int_0^\infty \nu(t) [x(0)]^T d\mu(t)\bigg] \nonumber \\
  & & \label{BeweisschlussWMPZB2}
      + \sum_{j=1}^l \int_0^\infty \big\langle \nu(t) g_{jx}\big(t,x_*(t)\big),x(t) \big\rangle \,d\mu_j(t).
\end{eqnarray}
Da s"amtliche Integralterme absolut integrierbar sind,
definiert die rechte Seite in (\ref{BeweisschlussWMPZB2}) ein stetiges lineares Funktional auf $C_0(\R_+,\R^n;\nu)$.
Wenden wir den Darstellungssatz \ref{SatzRieszClimnu} an
und setzen $p(t)=\displaystyle \int_t^\infty \nu(s) \, d\mu(s)$, so erhalten wir
\begin{eqnarray*}
p(t) &=& \int_t^\infty \big[ \varphi^T_x\big(s,x_*(s),u_*(s)\big)p(s)-\lambda_0 f_x\big(s,x_*(s),u_*(s)\big)\big]\, ds \\
     & & - \sum_{j=1}^l \int_t^\infty \nu(t) g_{jx}\big(s,x_*(s)\big) \,d\mu_j(s), \\
p(t_0) &=& l_0.
\end{eqnarray*}
Aus der Darstellung der Adjungierten $p(\cdot)$, $p(t)=\displaystyle \int_t^\infty \nu(s) \, d\mu(s)$,
ergeben sich ebenso wie in Abschnitt \ref{AbschnittBeweisWMP} die ``nat"urlichen'' Transversalit"atsbedingungen.
Damit sind (\ref{SatzWMPAufgabeZB1}) und (\ref{SatzWMPAufgabeZB2}) gezeigt. \\[2mm]
Gem"a"s (\ref{SatzWMPZBLMR2}) gilt f"ur alle $v(\cdot) \in L_\infty(\R_+,U)$ die Ungleichung
$$\int_0^\infty \big\langle H_u\big(t,x_*(t),u_*(t),p(t),\lambda_0\big),v(t)-u_*(t) \big\rangle \, dt \leq 0.$$
Daraus folgt abschlie"send via Standardtechniken f"ur Lebesguesche Punkte die Variationsungleichung (\ref{SatzWMPAufgabeZB3}).
Der Beweis von Theorem \ref{SatzWMPAufgabeZB} ist abgeschlossen. \hfill $\blacksquare$
       \subsection{Hinreichende Bedingungen nach Mangasarian} \label{AbschnittMangasarianWMPZB}
Mit $U_\gamma$ bezeichnen wir im Weiteren die Menge
$$U_\gamma= \{ (t,x,u) \in \R_+ \times \R^n \times \R^m\,|\, \|x-x(t)\| \leq \gamma, \|u-u(t)\| \leq \gamma\}.$$

\begin{theorem} \label{SatzHBWMPZB}
In der Aufgabe (\ref{WMPAufgabe1})--(\ref{WMPAufgabe4}) sei
$\big(x_*(\cdot),u_*(\cdot)\big)$ ein zul"assiger Steuerungsprozess,
f"ur den die Abbildungen $f(t,x,u)$ und $\varphi(t,x,u)$ auf der Menge
$$U_\gamma= \{ (t,x,u) \in \R_+ \times \R^n \times \R^m\,|\, \|x-x_*(t)\| \leq \gamma, \|u-u_*(t)\| \leq \gamma\}$$
stetig differenzierbar sind.
Au"serdem sei die Vektorfunktion $p(\cdot) \in L_2(\R_+,\R^n;\nu^{-1})$ st"uckweise stetig,
besitze h"ochstens abz"ahlbar viele Sprungstellen $s_k \in (0,\infty)$,
die sich nirgends im Endlichen h"aufen,
und $p(\cdot)$ sei zwischen diesen Spr"ungen stetig differenzierbar. \\
Ferne gelte:
\begin{enumerate}
\item[(a)] Das Tripel $\big(x_*(\cdot),u_*(\cdot),p(\cdot)\big)$
           erf"ullt (\ref{SatzWMPAufgabeZB1})--(\ref{SatzWMPAufgabeZB3}) mit $\lambda_0=1$ in Theorem \ref{SatzWMPAufgabeZB}.        
\item[(b)] F"ur jedes $t \in \R_+$ ist die Funktion $H\big(t,x,u,p(t),1\big)$ konkav in $(x,u)$
           und es sind die Funktionen $g_j(t,x)$, $j=1,...,l$, konvex bez"uglich $x$ auf $U_\gamma$.
\end{enumerate}
Dann ist $\big(x_*(\cdot),u_*(\cdot)\big)$ ein schwaches lokales Minimum der Aufgabe (\ref{WMPAufgabe1})--(\ref{WMPAufgabe4}).
\end{theorem}

\begin{bemerkung}{\rm
Der Teil (a) in Theorem \ref{SatzHBWMPZB} bedarf einer detaillierteren Diskussion.
Da wir von einer st"uckweise stetigen und zwischen den Sprungstellen stetig differenzierbaren Adjungierten $p(\cdot)$ ausgehen,
k"onnen wir die adjungierte Gleichung in Integraldarstellung in die Form einer st"uckweise definierten Differentialgleichung mit
Sprungbedingungen "uberf"uhren.
Es bezeichnen $0<s_1<...<s_d<T$ die Unstetigkeitsstellen der Adjungierten $p(\cdot)$ im Intervall $(0,T)$.
Dann gelten die Sprungbedingungen
$$p(s_k^-)-p(s_k^+)= -\sum_{j=1}^l \beta_j^k g_{jx}\big(s_k,x_*(s_k)\big),\qquad \beta_j^k = \mu_j(\{s_k\}) \geq 0, \quad k=1,...,d.$$
Ferner gibt es eine st"uckweise stetige Vektorfunktion $\lambda(\cdot):[0,T] \to \R^l$ derart,
dass die Differentialgleichung
$$\dot{p}(t)=- H_x\big(t,x_*(t),u_*(t),p(t),1\big) + \sum_{j=1}^l \lambda_j(t) g_{jx}\big(t,x_*(t)\big)$$
st"uckweise auf $(s_k,s_{k+1})$, $k=0,...,d$, gilt. Dabei haben wir $s_0=0$, $s_{d+1}=T$ gesetzt. \\
Abschließend halten wir fest,
dass wegen der Positivit"at der Ma"se $\mu_j$ und der Konzentration dieser Ma"se auf den Mengen
$$T_j=\big\{t \in \overline{\R}_+ \,\big|\, g_j\big(t,x_*(t)\big)=0\big\}, \quad j=1,...,l,$$
neben $\beta_j g_j\big(T,x_*(T)\big) =0$ die Bedingungen
$$\lambda_j(t) \geq 0, \qquad \lambda_j(t)g_j\big(t,x_*(t)\big)=0$$
auf $[0,T]$ und f"ur $k=1,...,d$ in den Sprungstellen 
$$\beta_j^k \geq 0, \qquad \beta_j^k g_j\big(s_k,x_*(s_k)\big)=0$$
f"ur $j=1,...,l$ und gelten. \hfill $\square$}
\end{bemerkung}

{\bf Beweis} Im Folgenden beachte man,
dass im Fall unendlich vieler Sprungstellen $s_k$ der Grenzwert
\begin{equation} \label{BeweisHBWMPZBGrenzwert}
\lim_{k \to \infty} \|p(s_k^-)-p(s_k^+)\|=0
\end{equation}
gilt, da $p(\cdot)$ eine Funktion beschr"ankter Variation ist. \\
Aus den elementaren Eigenschaften konkaver Funktionen ergibt sich
\begin{eqnarray*}
\lefteqn{H\big(t,x_*(t),u_*(t),p(t),1\big)-H\big(t,x(t),u(t),p(t),1\big)\geq} \\
&& -H_x\big(t,x_*(t),u_*(t),p(t),1\big)\big(x(t)-x_*(t)\big)
       -H_u\big(t,x_*(t),u_*(t),p(t),1\big)\big(u(t)-u_*(t)\big).
\end{eqnarray*}
Damit, sowie mit (\ref{SatzWMPAufgabeZB1}) und (\ref{SatzWMPAufgabeZB3}) gilt f"ur $T \not=s_k$:

\begin{eqnarray*}
    \Delta(T)
&=& \int_0^T \omega(t)\big[f\big(t,x(t),u(t)\big)-f\big(t,x_*(t),u_*(t)\big)\big] \, dt \\
&\geq& \int_0^T \big[\mathscr{H}\big(t,x_*(t),p(t)\big)-\mathscr{H}\big(t,x(t),p(t)\big)\big] \, dt 
    + \int_0^T \langle p(t), \dot{x}(t)-\dot{x}_*(t) \rangle dt \\
&\geq& \int_0^T \langle \dot{p}(t),x(t)-x_*(t)\rangle + \langle p(t), \dot{x}(t)-\dot{x}_*(t) \rangle \, dt \\
&&    \hspace*{10mm} - \int_0^T \sum_{j=1}^l \lambda_j(t) \big\langle g_{jx}\big(t,x_*(t)\big) , x(t)-x_*(t) \big\rangle \, dt \\
&&    \hspace*{10mm} + \sum_{s_k <T} \langle p(s_k^-)-p(s_k^+),x(s_k)-x_*(s_k)\rangle \\
&\geq& \int_0^T \langle \dot{p}(t),x(t)-x_*(t)\rangle + \langle p(t), \dot{x}(t)-\dot{x}_*(t) \rangle \, dt
       = \langle p(T),x(T)-x_*(T)\rangle.
\end{eqnarray*}
Unter Beachtung des Grenzwertes (\ref{BeweisHBWMPZBGrenzwert}) folgt abschlie"send mit den ``nat"urlichen'' Transversalit"atsbedingungen
in (\ref{SatzWMPAufgabeZB2}) die Beziehung
$$\lim_{T \to \infty} \Delta(T) \geq \lim_{T \to \infty} \langle p(T),x(T)-x_*(T)\rangle=0$$
f"ur alle zul"assigen $\big(x(\cdot),u(\cdot)\big)$ mit 
$\|x(\cdot)-x_*(\cdot)\|_\infty, \|u(\cdot)-u_*(\cdot)\|_{L_\infty} \leq \gamma$. \hfill $\blacksquare$

\begin{beispiel} \index{Ressourcenabbau}
{\rm Im Beispiel \ref{ExampleRessource} zum Abbau einer nicht erneuerbaren Ressource ist die Funktion
$$H\big(t,x,u,p_1(t),p_2(t),1\big) = p_1(t) (-u)+p_2(t) cf(u) + e^{-\varrho t}[f(u)-ry-qu]$$
nach den Voraussetzungen an die Funktion $f(u)$ konkav.
Damit liefern die ermittelten Kandidaten in den F"allen (A) und (C) schwache lokale Minimalstellen der Aufgabe. \hfill $\square$}
\end{beispiel}   
       
       \rhead[]{Bemerkungen \hspace*{1mm} \thepage}
       \section{Bemerkungen}
In diesem Kapitel widmeten wir unsere Aufmerksamkeit der Aufgabe (\ref{WMPAufgabe1})--(\ref{WMPAufgabe4}) mit unendlichem Zeithorizont.
Dabei konzentrierten wir uns auf den Zugang im Rahmen gewichteter Funktionenr"aume.
Den Ansto"s dazu lieferte die Frage nach der L"osbarkeit der linearen Integralgleichung
$$x(t)+ \int_0^t A(s) x(s) \, ds = z(t), \qquad t \in \R_+,$$
in einem m"oglichst umfassenden Rahmen.
Eine Antwort lieferten die gewichteten R"aume stetiger Funktionen.
Als Konsequenz ergab sich f"ur die Wahl der Gewichtsfunktion $\nu(\cdot)$ die Kontraktionsbedingung
$$\sup_{t \in \R_+} \int_0^t \frac{\nu(t)}{\nu(s)} \big\|\varphi_x\big(s,x_*(s),u_*(s)\big)\big\| \, ds < 1.$$
Bezogen auf die Aufgabe (\ref{WMPAufgabe1})--(\ref{WMPAufgabe4}) flie"st die Gewichtsfunktion unmittelbar in die Annahmen (A$_1$)--(A$_3$) durch die Menge 
$$U_{\gamma,\nu}= \{ (t,x,u) \in \R_+ \times \R^n \times \R^m\,|\, \nu(t)\|x-x(t)\| \leq \gamma, \|u-u(t)\| \leq \gamma\}$$
ein.
Die Menge $U_{\gamma,\nu}$ ist bez"uglich der Zustandsvariablen verbunden mit dem Umgebungsbegriff im gewichteten Raum stetiger Funktionen,
denn es gilt:
$$\|x(\cdot) - x_*(\cdot)\|_{\infty,\nu} \leq \varepsilon \quad\Leftrightarrow\quad \nu(t)\|x(t) - x_*(t)\| \leq \varepsilon \mbox{ f"ur alle } t \in \R_+.$$
D.\,h., je gr"o"ser der Parameter $a>0$ gew"ahlt wurde,
desto schneller w"achst der Durchmesser der Menge $\{x \in \R^n | \|x\| \leq \gamma e^{at}\}$.
Deswegen kann die expandierende Umgebung in vielen Aufgaben zu sehr starken bzw. sogar zu nicht erf"ullbaren Einschr"ankungen f"uhren.
Um dies zu verdeutlichen, geben wir ein einfaches Beispiel an.

\begin{beispiel}\label{BeispielRegler} {\rm Wir betrachten den linear-quadratischen Regler:
\begin{eqnarray*}
&& J\big(x(\cdot),u(\cdot)\big) = \int_0^\infty e^{-2t} \cdot \frac{1}{2}\big( x^2(t)+u^2(t)\big) \, dt \to \inf, \\
&& \dot{x}(t) = 2 x(t)+u(t), \qquad x(0)=2, \qquad u(t) \in \R.
\end{eqnarray*}
Wegen dem Term $e^{-2t} x^2(t)$ im Integranden ist das Zielfunktional ausschlie"slich f"ur Trajektorien $x(\cdot) \in C_0(\R_+,\R^n;\nu)$ 
zu Gewichten $\nu(t)=e^{-at}$ mit $a<1$ wohldefiniert.
Demgegen"uber gilt in Thereom \ref{SatzWMP} die Kontraktionsbedingung (\ref{WMPBedingung}) genau dann, wenn $a>2$ ist.
Daher darf das Theorem \ref{SatzWMP} nicht in diesem Beispiel angewendet werden. \hfill $\square$}
\end{beispiel}

Dieses Beispiel demonstriert,
dass die Daten der Aufgabe auf der expandierenden Menge $U_{\gamma,\nu}$ zu stark wachsen und somit die Annahmen an die Aufgabe nicht erf"ullt sein k"onnen.
Wiederum das Beispiel \ref{BeispielDockner} beinhaltet eine Dynamik der Form
$$\dot{x}(t)=x(t) \big(\alpha-r\ln x(t) \big) -u(t)x(t),$$
die f"ur eine beschr"ankte Trajektorie $x(\cdot)$ auf keiner Menge $U_{\gamma,\nu}$ wohldefiniert ist.
Treten demnach in der Aufgabe Terme auf,
die nur gewisse Bereiche abbilden k"onnen,
so ist der Rahmen gewichteter R"aume meist ungeeignet. \\[2mm]
In Pickenhain \cite{Pickenhain} geh"oren bei einem zul"assigen Steuerungsprozess sowohl der Zustand als auch die Steuerung einem gewichteten Raum an.
Im Gegensatz zur vorliegenden Herangehensweise sind damit unbeschr"ankte Steuerungen in \cite{Pickenhain} inbegriffen. \\
Durch die Anpassung der Voraussetzungen (A$_0$)--(A$_3$) auf eine Menge,
die auch bez"uglich der Steuerungsvariable expandiert,
lassen sich die Schwachen Optimalit"atsprinzipien dieses Kapitels nachweisen.
Wir betrachten unbeschr"ankte Steuerungen in einer Aufgabe der Neoklassischen Wachstumstheorie (vgl. z.\,B. \cite{Arnold,Barro}):

\begin{beispiel}\label{BeispielNoPonzi} \index{Kapitalakkumulation}
{\rm Mit der isoelastischen Nutzenfunktion
$$U(C)=\frac{C^{1-\sigma}-1}{1-\sigma} \quad\mbox{f"ur } \sigma >0,\sigma\not=1
  \qquad\mbox{bzw.}\qquad U(C)=\ln(C) \quad\mbox{f"ur } \sigma=1,$$
diskutieren wir die Aufgabe
\begin{eqnarray*}
&& J\big(K(\cdot),C(\cdot)\big) = \int_0^\infty e^{-\varrho t} U\big(C(t)\big) \, dt \to \sup, \\
&& \dot{K}(t) = rK(t)+W(t)-C(t), \qquad K(0)=K_0>0, \qquad C(t) > 0.
\end{eqnarray*}
Die adjungierte Gleichung (\ref{WMP1}) und die Variationsungleichung (\ref{WMP3}) liefern
$$\dot{p}(t)=-rp(t), \qquad C_*^{-\sigma}(t)=p(t)e^{\varrho t}, \qquad C_*^{-\sigma}(t)=p_0 e^{-(r-\varrho)t}.$$
Es folgt daraus die bekannte Ramsey-Regel der konstanten Wachstumsrate:
$$\frac{\dot{C}_*(t)}{C_*(t)}=\frac{r-\varrho}{\sigma}.$$
Da aber keine Beschr"ankungen an die Entwicklung des Kapitalstocks vorliegen, liefert die Maximierung des Zielfunktionals
$$J\big(K(\cdot)\big) = \int_0^\infty e^{-\varrho t} U\big(C(t)\big) \, dt$$
die wenig brauchbare L"osung $C_*(t)\equiv \infty$. \\
Oft wird an dieser Stelle ausgeschlossen, dass sich ein Haushalt f"ur immer verschuldet,
indem er alte Schulden durch Aufnahme immer neuer Kredite finanziert.
Das bedeutet, dass der Barwert der Ausgaben den Barwert der Einnahmen nicht "ubersteigen darf:
$$K(0)+\int_0^\infty W(t) e^{-rt} \, ds \geq \int_0^\infty C(t) e^{-rt} \, dt.$$
Bei Vergleich dieser Ungleichung mit der expliziten L"osung der linearen Dynamik,
$$K(t) = e^{rt} \bigg[K(0)+\int_0^t \big(W(s)-C(s)\big) e^{-rs} \, ds \bigg],$$
erhalten wir die ``No-Ponzi'' Bedingung \index{No-Ponzi Bedingung}
$$\lim_{t \to \infty} K(t) e^{-rt} = 0.$$
Speziell im Fall $\sigma=1$ und $W(t) \equiv W$ erhalten wir aus den konstanten Wachstumsraten
$$C_*(t)=c_0 e^{(r-\varrho)t}$$
und ferner für die Entwicklung des Kapitalbestandes
$$K(t) = e^{rt}\bigg[K(0)+\int_0^t\big(W-C_*(s)\big) e^{-rs} \, ds \bigg]
       = e^{rt}\bigg[K(0) - \frac{W}{r} (e^{-rt}-1) + \frac{c_0}{\varrho}(e^{-\varrho t}-1) \bigg].$$
Die zus"atzliche ``No-Ponzi'' Randbedingung ist also nur dann erf"ullt,
wenn
$$K(0) + \frac{W}{r} - \frac{c_0}{\varrho}=0$$
gilt.
Unter der Bedingung im Uendlichen ergibt sich also eine eindeutige L"osung. \\
Die Anforderung $K(\cdot) \in W^1_2(\R_+,\R^n;\nu)$ an eine zul"assige Trajektorie f"uhrt unmittelbar zu der Bedingung
$$\lim_{t \to \infty} K(t) e^{-at} = 0,$$
wobei wegen der Kontraktionsbedingung $a>r$ gelten muss.
Aus der expliziten Darstellung der Funktionen $K(\cdot)$ folgt offensichtlich,
dass stets dieser Grenzwert erf"ullt ist.
Dementsprechend stellt die ``No-Ponzi'' Bedingung eine zus"atzliche Restriktion an eine zul"assige Trajektorie dar. \hfill $\square$}
\end{beispiel}

In der Aufgabe mit Zustandsbeschr"ankungen ergibt sich die Frage nach den aktiven Ungleichungen.
Denn die Ungleichungen (\ref{WMPAufgabe4}) besitzen "uber dem unbeschr"ankten Intervall $\R_+$ einen anderen Charakter
als in den Aufgaben mit endlichem Zeithorizont.
Dies liegt darin begr"undet,
dass im Gegensatz zum unbeschr"ankten Intervall jede stetige Funktion "uber einer kompakten Menge stets ein Maximum
und ein Minimum besitzt.
Diese Eigenschaft geht "uber $\R_+$ verloren.
Eine Art von Ungleichungen, die sich "uber dem unendlichen Zeithorizont ergeben,
sind nun diejenigen, die nur im Unendlichen aktiv sind:
 
\begin{beispiel}
{\rm Wir f"ugen im linear-quadratischen Regler eine zus"atzliche Zustandsbeschr"ankung ein:
\begin{eqnarray*}
&& J\big(x(\cdot),u(\cdot)\big)=\int_0^\infty e^{-2t}\frac{1}{2}\big(x^2(t)+u^2(t)\big) \, dt \to \inf, \\
&& \dot{x}(t)=2x(t)+u(t), \quad x(0)=2, \quad x(t) \geq 0, \quad u(t) \in \R.
\end{eqnarray*}
Die Beschr"ankung $-x(t) \leq 0$ hat offenbar keinen Einfluss auf die globale L"osung
$$x_*(t)=2e^{(1-\sqrt{2})t}, \quad u_*(t)= -2(1+\sqrt{2})e^{(1-\sqrt{2})t}.$$
Jedoch hat die Trajektorie $x_*(\cdot)$ die Eigenschaften 
$$-x_*(t) < 0 \quad\mbox{f"ur alle } t \in \R_+ \quad\mbox{ und }\quad \sup_{t \in \R_+} \big(-x_*(t)\big) = 0.$$
Daher ist f"ur $x_*(\cdot)$ die Zustandsbeschr"ankung im Unendlichen aktiv.} \hfill $\square$
\end{beispiel}

Im letzten Beispiel hat die Zustandsbeschr"ankung,
die nur im Unendlichen aktiv wird,
keinen Einfluss auf die L"osung der Aufgabe.
Intuitiv k"onnte man daher meinen,
dass es sich bei einer Zustandsbeschr"ankung, die im Unendlichen aktiv ist, um eine nichtaktive Beschr"ankung handelt.
Allerdings stellt das einen fehlerhaften Schluss dar:

\begin{beispiel} \label{BeipielZBUnendlich}
{\rm Wir betrachten die Aufgabe
\begin{eqnarray*}
&& J\big(x(\cdot),u(\cdot)\big)=\int_0^\infty e^{-\varrho t}\big(1-u(t)\big)x(t) \, dt \to \sup, \\
&& \dot{x}(t)=u(t)x(t), \qquad x(0)=1, \qquad u \in [0,1], \qquad \varrho \in (0,1), \\
&& \dot{z}(t)=e^{-\varrho t}x(t), \qquad z(0)=0, \qquad z(t) \leq Z \quad \mbox{f"ur alle } t \in \R_+, \qquad Z> \frac{1}{\varrho}.
\end{eqnarray*}
Die Zustandsgleichung und -beschr"ankung f"ur die Trajektorie $z(\cdot)$ ergibt sich aus der isoperimetrischen Nebenbedingung in Form einer Budgetbeschr"ankung
$$\int_0^\infty e^{-\varrho t}x(t) \, dt \leq Z.$$
Offensichtlich ist $\dot{z}(t) >0$ auf $\R_+$ und damit $z(t)$ streng monoton wachsend.
Demzufolge kann die Beschr"ankung $z(t) \leq Z$ erst im Unendlichen aktiv werden und greift nur durch das Verhalten in $t=\infty$ in die
gestellte Aufgabe ein. \\
Da stets $\dot{x}(t)\geq 0$ ist,
muss f"ur jede zul"assige Trajektorie $e^{-\varrho t}x(t) \to 0$ f"ur $t \to \infty$ gelten,
denn nur dann ist $z(t) \leq Z$ erf"ullt.
Damit erhalten wir f"ur zul"assige Steuerungsprozesse zun"achst im Zielfunktional
\begin{eqnarray*}
    J\big(x(\cdot),z(\cdot),u(\cdot)\big)
&=& \int_0^\infty e^{-\varrho t}\big(1-u(t)\big)x(t) \, dt = \int_0^\infty \dot{z}(t) \, dt - \int_0^\infty e^{-\varrho t} \dot{x}(t) \, dt \\
&=& \int_0^\infty \dot{z}(t) \, dt + 1 - \varrho \int_0^\infty e^{-\varrho t} x(t) \, dt \leq 1 + (1-\varrho) Z.
\end{eqnarray*}
Es ergibt sich also f"ur das Zielfunktional die obere Schranke $1 + (1-\varrho) Z$.
D.\,h., dass jeder Steuerungsprozess $\big(x(\cdot),z(\cdot),u(\cdot)\big)$,
f"ur den die Zustandsbeschr"ankung $z(t) \leq Z$ im Unendlichen aktiv wird,
global optimal ist und $J\big(x(\cdot),z(\cdot),u(\cdot)\big) = 1 + (1-\varrho) Z$ gilt.
Demnach liefert
$$x_*(t) = e^{\alpha t},\qquad z_*(t)= \frac{1}{\alpha - \varrho} (e^{(\alpha-\varrho)t}-1),\qquad
  u_*(t)=\alpha,\qquad \alpha=\varrho-\frac{1}{Z} \in (0, \varrho)$$
ein globales Maximum.
Wir wenden formal die notwendigen Optimalit"atsbedingungen an:
Die Pontrjagin-Funktion lautet im normalen Fall
$$H(t,x,z,u,p,q,1)= pux+qe^{-\varrho t}x + e^{-\varrho t}(1-u)x.$$
Daraus ergeben sich die adjungierten Gleichungen
$$\dot{p}(t) = -u_*(t) p(t) -q(t) e^{-\varrho t} - e^{-\varrho t}\big(1-u_*(t)\big), \qquad \dot{q}(t) = 0$$
und f"ur alle $u \in [0,1]$ die Variationsungleichungen
$$\big\langle H_u\big(t,x_*(t),u_*(t),p(t),\lambda_0\big),\big(u-u_*(t)\big) \big\rangle\leq 0
  \quad\Leftrightarrow\quad \big(p(t)- e^{-\varrho t}\big) \cdot \big(u-u_*(t)\big) \leq 0.$$
Wegen $u_*(t) = \alpha \in (0,1)$ erhalten wir f"ur die Adjungierten
$$p(t)= e^{-\varrho t}, \qquad q(t)= \varrho-1.$$
D.\,h., dass bez"uglich der Zustandsbeschr"ankung $z(t) \leq Z$ das korrespondierende positive Ma"s in $t=\infty$
mit dem Atom $\mu(\{\infty\})=1-\varrho$ konzentriert sein muss.} \hfill $\square$
\end{beispiel}

Das Beispiel \ref{BeipielZBUnendlich} zeigt,
dass die im Unendlichen aktiven Ungleichungen wesentlich in die Aufgabe eingreifen k"onnen.
Dar"uberhinaus liefern die formal gebildeten notwendigen Bedingungen,
dass in diesem Fall die vollst"andige Information "uber die Multiplikatoren "uber $\overline{\R}_+$ n"otig sind.
Die Behandlung des unendlichen Zeithorizontes "uber $\R_+$ bzw. $\overline{\R}_+$ birgt bez"uglichdes theoretischen Hintergrundes gravierende Unterschiede in sich,
auf die wir nun eingehen werden. \\
Auff"allig im Beweis im Abschnitt \ref{AbschnittBeweisWMPZB} ist, 
dass wir die Funktionen $G_j$ als Abbildungen in den Raum $C_{\lim}(\R_+,\R^n)$ auffassen obwohl unter den getroffenen Annahmen stets
$$\lim_{t \to \infty} G_j\big(x(\cdot)\big)(t) = 0$$
gilt.
Diesbez"uglich bemerken wir zun"achst, dass der Kegel
$$\mathscr{K} =\{ x(\cdot) \in C_0(\R_+,\R^l) \,|\, x(\cdot) \preceq 0\}$$
als Teilmenge des Raumes $C_0(\R_+,\R^l)$ ein leeres Inneres besitzt.
Denn zu jeder Funktion $x(\cdot) \in C_0(\R_+,\R)$ mit $x(\cdot) \preceq 0$ und zu jedem $\varepsilon >0$ gibt es eine Zahl $T>0$ mit 
$$|x(t)| \leq \varepsilon \qquad\mbox{f"ur alle } t \geq T.$$
Ferner k"onnen wir eine Zahl $T' \in [T,T+2\varepsilon]$ derart angeben,
dass die Funktion
$$y(t)= \left\{ \begin{array}{ll} x(t), & t \in [0,T), \\ x(T)+(t-T),& t \in [T,T'), \\ -x(t),& t \in [T',\infty) \end{array} \right.$$
dem Raum $C_0(\R_+,\R)$ aber nicht dem Kegel $\mathscr{K}$ angeh"ort und $\|x(\cdot)-y(\cdot)\|_\infty \leq \varepsilon$ gilt. \\
Demgegen"uber gilt aber ${\rm int\,}\mathscr{K} \not= \emptyset$ f"ur
$\mathscr{K} =\{ x(\cdot) \in C_{\lim}(\R_+,\R^l) \,|\, x(\cdot) \preceq 0\}$ als Kegel im Raum $C_{\lim}(\R_+,\R^l)$
und die Voraussetzungen des Schwachen Extremalprinzips \ref{SatzExtremalprinzipSchwach} an den Kegel $\mathscr{K}$ sind erf"ullt. \\
Ferner zeigen die verschiedenen Versionen des Satzes Riesz im Anhang \ref{AnhangRiesz},
dass die Darstellung der stetigen linearen Funktionale im Raum $C_{\lim}(\R_+,\R^n)$ gegen"uber dem Raum $C_0(\R_+,\R^n)$ auf den Abschluss $\overline{\R}_+$
fortgesetzt werden.
Auf diesem Weg ist es m"oglich eine vollst"andige Information,
insbesondere das in $t=\infty$ konzentrierte Borelsche Ma"s,
zu erhalten.
Ein weiterer wesentlicher Punkt den Raum $C_0(\R_+,\R^n)$ durch den Raum $C_{\lim}(\R_+,\R^n)$ zu ersetzen ergibt sich im Anhang \ref{AbschnittKonvexeAnalysis}
zur Konvexen Analysis.
In den Beispeilen \ref{SubdifferentialMaximum5} und \ref{SubdifferentialMaximum6} wird gezeigt,
dass die Darstellung des Subdifferentials der Supremumsfunktion im Raum $C_0(\R_+,\R^n)$ unscharf ist.
Im Gegensatz dazu ergibt sich im Raum $C_{\lim}(\R_+,\R^n)$ bez"uglich der Maximumfunktion die bekannte Darstellung
(Beispiele \ref{SubdifferentialMaximum3}, \ref{SubdifferentialMaximum4}).  \\[2mm]
Die "Uberf"uhrung der Aufgabe (\ref{WMPAufgabe1})--(\ref{WMPAufgabe4}) in die Extremalaufgabe \ref{ExtremalaufgabeWMPZB} und die Einbettung in den Rahmen
gewichteter Funktionenr"aume hat weiterhin zur Folge,
dass die Abbildungen $G_j\big(x(\cdot)\big)(t) = \nu(t) g_j\big(t,x(t)\big)$ stets im Unendlichen verschwinden m"ussen.
Dementsprechend sind zwar stets die ``nat"urlichen'' Transversalit"atsbedingungen erf"ullt,
aber der Zugang birgt die in $t=\infty$ konzentrierten Ma"se in sich.
Bezogen auf das Beispiel \ref{BeipielZBUnendlich} bedeutet dies,
dass f"ur das globale Maximum
$$x_*(t) = e^{\alpha t},\qquad z_*(t)= \frac{1}{\alpha - \varrho} (e^{(\alpha-\varrho)t}-1),\qquad
  u_*(t)=\alpha,\qquad \alpha=\varrho-\frac{1}{Z} \in (0, \varrho)$$
die notwendigen Bedingungen genau dann erf"ullt sind, wenn 
$$\lambda_0=0, \qquad p(t) \equiv 0, \qquad q(t) \equiv 0$$
und das regul"are Borelsche Ma"s in $t=\infty$ konzentriert sind. \\[2mm]
In der Einleitung dieses Kapitels haben wir bereits erw"ahnt,
dass wir im Rahmen der gewichteten Funktionenr"aume Aufgaben mit Randwerten im Unendlichen nicht behandeln werden.
Eine Randbedingung im Unendlichen ist im Beispiel \ref{BeispielNoPonzi} in Form der ``No-Ponzi'' Bedingung aufgetreten.
Der Verzicht auf diese Restriktionen ist darin begr"undet,
dass der Operator
$$H_1:C_{\lim}(\R_+,\R^n;\nu) \to \R^n, \qquad H_1\big(x(\cdot)\big)= \lim_{t \to \infty} x(t),$$
im Raum $C_{\lim}(\R_+,\R^n;\nu)$ nirgends stetig ist.
Denn f"ur jede Funktion $x(\cdot)$, f"ur die der Grenzwert im Unendlichen existiert, liefert
$$y(t)=x(t) +\varepsilon \sin(t)$$
eine Funktion, f"ur die der Operator $H_1$ nicht erkl"art ist und f"ur die $\|x(\cdot)-y(\cdot)\|_{\infty,\nu} \leq \varepsilon$ gilt.
Da der Operator $H_1$ nirgends stetig ist,
kann demzufolge das Extremalprinzip in Anhang \ref{AbschnittEPschwach} nicht angewendet werden.

\begin{appendix}
\lhead[\thepage \hspace*{1mm} Funktionalanalytische Hilfsmittel]{}
\rhead[]{Funktionalanalytische Hilfsmittel \hspace*{1mm} \thepage}
\section{Funktionalanalytische Hilfsmittel}
\subsection{Grundprinzipien der Funktionalanalysis}
\begin{theorem}[Satz von Hahn-Banach; Fortsetzungsversion]  \index{Satz, Darstellungssatz von Riesz!von HahnBanach@-- von Hahn-Banach}
Sei $X$ ein normierter Raum und $U$ ein Untervektorraum.
Zu jedem stetigen linearen Funktional $u^*:U \to \R$ existiert dann ein stetiges lineares Funktional Funktional
$x^*: X \to \R$ mit
$$x^*\big|_U=u^*, \qquad \|x^*\|=\|u^*\|.$$
\end{theorem}

\begin{folgerung} \label{FolgerungAnnulator}
Seien $X$ ein normierter Raum, $U$ ein abgeschlossener Unterraum und $x \in X$, $x \not\in U$.
Dann existiert ein $x^* \in X^*$ mit
$$x^*\big|_U=0,\qquad \langle x^*,x \rangle \not=0.$$
\end{folgerung}

\begin{theorem}[Satz von Hahn-Banach; Trennungsversion]  \index{Satz, Darstellungssatz von Riesz!Trennungssatz@--, Trennungssatz}
Seien $X$ ein normierter Raum, $V_1,V_2 \subseteq X$ konvex und $V_1$ offen.
Es gelte $V_1 \cap V_2 = \emptyset$.
Dann existiert ein $x^* \in X^*$ mit
$$\langle x^*, v_1 \rangle < \langle x^*, v_2 \rangle \qquad \mbox{ f"ur alle } v_1 \in V_1, v_2 \in V_2.$$
\end{theorem}

Eine Abbildung $T$ hei"st offen\index{Abbildung, beschr"ankte!offene@--, offene}, wenn $T$ offene Mengen auf offene Mengen abbildet.

\begin{theorem}[Satz von der offenen Abbildung] \label{SatzOffeneAbbildung}  \index{Satz, Darstellungssatz von Riesz!von der offen@-- von der offenen Abbildung}
Es seien $X$, $Y$ Banachr"aume und $T \in L(X,Y)$ surjektiv.
Dann ist $T$ offen.
\end{theorem}

Seien $X,Y$ normierte R"aume und $T \in L(X,Y)$.
Der adjungierte Operator\index{adjungierter Operator} $T^*:Y^* \to X^*$ ist durch
$\langle T^*y^*, x \rangle = \langle y^*, Tx \rangle$
definiert. Offensichtlich folgt daraus $T^* \in L(Y^*,X^*)$.
Seien nun $U \subseteq X$ und $V \subseteq X^*$.
Wir definieren die Mengen
\begin{eqnarray*}
U^\perp &=& \{ x^* \in X^* \,|\, \langle x^*, x \rangle =0 \mbox{ f"ur alle } x \in U\}, \\
V_\perp &=& \{ x \in X \,|\, \langle x^*, x \rangle =0 \mbox{ f"ur alle } x^* \in V\}.
\end{eqnarray*}

\begin{lemma}[Satz vom abgeschlossenen Bild] \label{SatzAbgeschlossenesBild} \index{Satz, Darstellungssatz von Riesz!vom ab@-- vom abgeschlossenen Bild}
Seien $X$, $Y$ Banachr"aume, und es sei $T \in L(X,Y)$.
Dann gelten die "Aquivalenzen:
\begin{eqnarray*}
{\rm Im\,}T \mbox{ ist abgeschlossen } &\Leftrightarrow& {\rm Im\,}T=({\rm Ker\,}T^*)_\perp \\
\Leftrightarrow \; {\rm Im\,}T^* \mbox{ ist abgeschlossen} &\Leftrightarrow&  {\rm Im\,}T^*=({\rm Ker\,}T)^\perp.
\end{eqnarray*}
\end{lemma}

\begin{satz}[Fixpunktsatz von Weissinger] \label{SatzWeissinger} \index{Satz, Darstellungssatz von Riesz!Fixpunkt@--, Fixpunktsatz von Weissinger}
Es sei $U$ eine nichtleere abgeschlossene Teilmenge des Banachraumes $X$,
ferner $\displaystyle \sum_{n=1}^\infty a_n$ eine konvergente Reihe positiver Zahlen und
$A:U \to U$ eine Selbstabbildung von $U$ mit
$$\|A^n u-A^n v\| \leq a_n\|u-v\| \qquad\mbox{ f"ur alle } u,v \in U,\; n \in \N.$$
Dann besitzt $A$ genau einen Fixpunkt, d.\,h. es gibt genau ein $u \in U$ mit $Au=u$. \\
Dieser Fixpunkt ist Grenzwert der Iterationsfolge $u_n=Au_{n-1}, n=1,2,...,$ bei beliebigem Startwert $u_0 \in U$.
Schlie"slich gilt die Fehlerabsch"atzung
$$\|u-u_n\| \leq \|u_1-u_0\| \cdot \sum_{k=n}^\infty a_k.$$
\end{satz}

\subsection{Der Darstellungssatz von Riesz} \label{AnhangRiesz}
Es sei $I \subseteq \R_+$ und es bezeichne $\mathscr{M}(I)$ die Menge der signierten regul"aren Borelschen Ma"se auf der
Borelschen $\sigma$-Algebra auf $I$.
Au"serdem bezeichnet $C_{\lim}(\R_+,\R^n)$ den Raum der stetigen Funktionen $x(\cdot)$, die im Unendlichen einen Grenzwert besitzen.
Als abgeschlossener Unterraum des Raumes $C_b(\R_+,\R^n)$ ist $C_{\lim}(\R_+,\R^n)$ vollst"andig. \\
Es sei $\nu(\cdot)$ eine positive stetige Funktion.
Mit $C_{\lim}(\R_+,\R^n;\nu)$ bezeichnen wir den Raum stetiger Funktionen,
die bez"uglich der Gewichtsfunktion $\nu(\cdot)$ im Unendlichen einen Grenzwert besitzen:
$$C_{\lim}(\R_+,\R^n;\nu) = \big\{ x(\cdot) \in C(\R_+,\R^n) \,\big|\, \lim_{t \to \infty}\nu(t) x(t) \mbox { existiert} \big\}.$$
Wir versehen den Raum  $C_{\lim}(\R_+,\R^n;\nu)$ mit der gewichteten Supremumsnorm
$$\|x(\cdot)\|_{\infty,\nu} = \sup_{t \in \R_+} \nu(t) \|x(t)\|.$$
Offensichtlich gelten damit die "Aquivalenzen
$$x(\cdot) = \nu^{-1}(\cdot) y(\cdot) \in  C_{\lim}(\R_+,\R^n;\nu) \Leftrightarrow \nu(\cdot)x(\cdot) =  y(\cdot) \in C_{\lim}(\R_+,\R^n;\nu)$$
und $\|x(\cdot)\|_{\infty,\nu} = \|y(\cdot)\|_\infty.$
Damit ist $C_{\lim}(\R_+,\R^n;\nu)$ bez"uglich $\|\cdot\|_{\infty,\nu}$ vollst"andig.
Der Raum $C_0(\R_+,\R^n;\nu)$ der stetigen Funktionen,
die bez"uglich der positiven und stetigen Gewichtsfunktion $\nu(\cdot)$ im Unendlichen verschwinden,
ist versehen mit der gewichteten Supremumsnorm ein abgeschlossener Unterraum von $C_{\lim}(\R_+,\R^n;\nu)$ und daher vollst"andig.

\begin{satz}[Rieszscher Darstellungssatz] \label{SatzRieszC0}\index{Satz, Darstellungssatz von Riesz}\index{Darstellungssatz von Riesz}
Der Dualraum $C_0^*(\R_+,\R)$ ist isometrisch isomorph zu $\mathscr{M}(\R_+)$ unter der Abbildung
$$\Lambda:\mathscr{M}(\R_+) \to C_0^*(\R_+,\R), \qquad \Lambda(\mu)x(\cdot)= \int_{\R_+} x(t) \, d\mu(t).$$
\end{satz}

\begin{satz}[Rieszscher Darstellungssatz] \label{SatzRieszClimnu}\index{Satz, Darstellungssatz von Riesz}\index{Darstellungssatz von Riesz}
Der Dualraum $C_{\lim}^*(\R_+,\R^n;\nu)$ ist unter der Abbildung
$$\Lambda(\mu)x(\cdot)= \int_{\R_+} \nu(t) \langle x(t) , d\mu_0(t) \rangle + \lim_{t \to \infty} \langle \nu(t) x(t) , \mu_\infty \rangle$$
isometrisch isomorph zu den signierten Vektorma"sen $\mu \in \mathscr{M}(\overline{\R}_+)$.
Dabei besitzt $\mu$ die Darstellung $\mu=\mu_0+\mu_\infty$
mit einem $\mu_0 \in \mathscr{M}(\R_+)$ und einem in $t =\infty$ konzentrierten signierten Ma"s $\mu_\infty$.
\end{satz}

{\bf Beweis} Wir betrachten die stetige lineare Abbildung $T:C_{\lim}(\R_+,\R^n;\nu) \to C_0(\R_+,\R^n)$,
$$Tx(\cdot)(t)=\nu(t)x(t) - \lim_{t \to \infty} \nu(t)x(t).$$
Der Setzung nach bildet $T$ auf den gesamten Raum $C_0(\R_+,\R^n)$ ab.
Ferner gilt
$${\rm Ker\,}T = \{ x(\cdot) \in C_{\lim}(\R_+,\R^n;\nu) \,|\, \nu(t) x(t)= \mbox{konstant}\}.$$
Sei $x^* \in C_{\lim}^*(\R_+,\R^n;\nu)$ und sei
$a \in \R^n$ mit den Komponenten $a_i = \langle x^*,e_i(\cdot) \rangle$,
wobei die $i$-te Komponente der Funktion $e_i(\cdot)$ identisch Eins ist und alle weiteren Komponenten identisch gleich Null sind.
Wir definieren das Funktional $x_1^*$ durch
$$\langle x_1^*, x(\cdot) \rangle = \langle x^*, x(\cdot) \rangle - \lim_{t \to \infty} \langle a,\nu(t) x(t) \rangle.$$
Dann gilt $x_1^* \in ({\rm Ker\,}T)^\perp$ und es existiert nach dem Satz vom abgeschlossenen Bild ein $y^* \in C_0^*(\R_+,\R^n)$ mit $x_1^*= T^*y^*$.
Daraus folgt mit dem Rieszschen Darstellungssatz
\begin{eqnarray*}
\langle x^*, x(\cdot) \rangle &=& \langle x_1^*, x(\cdot) \rangle + \lim_{t \to \infty} \langle a,\nu(t) x(t) \rangle
                                  = \langle y^*, Tx(\cdot) \rangle +\lim_{t \to \infty} \langle a,\nu(t) x(t) \rangle \\
              &=& \int_{\R_+}  \big\langle \nu(t)x(t)-\lim_{t \to \infty}\nu(t) x(t) , d\mu_0(t) \big\rangle + \lim_{t \to \infty} \langle a,\nu(t) x(t) \rangle \\
              &=& \int_{\R_+}  \nu(t) \langle x(t) , d\mu_0(t) \rangle + \lim_{t \to \infty} \langle \nu(t) x(t) , \mu_\infty \rangle.
\end{eqnarray*}
Der Darstellungssatz ist damit nachgewiesen. \hfill $\blacksquare$

\begin{folgerung} \label{FolgerungRieszC0nu}
Der Dualraum $C_0^*(\R_+,\R^n;\nu)$ ist unter der Abbildung
$$\Lambda(\mu)x(\cdot)= \int_{\R_+} \nu(t) \langle x(t) , d\mu(t) \rangle$$
isometrisch isomorph zu den signierten Vektorma"sen $\mu \in \mathscr{M}(\R_+)$.
\end{folgerung}
       

\subsection{Der Satz von Ljusternik}
In diesem Abschnitt befassen wir uns mit dem fundamentalen Satz von Ljusternik (Ljusternik \cite{Ljusternik}).
Die vorliegende zusammenfassende Darstellung und die Verallgemeinerung ist Ioffe \& Tichomirov \cite{Ioffe} entnommen.
Eine vollst"andige Beweisf"uhrung ist wiederum in Ioffe \& Tichomirov \cite{Ioffe} zu finden.

\begin{definition}[Lokaler Tangentialkegel]
Sei $X$ ein Banachraum und $x_0 \in M \subseteq X$.
Mit $\mathscr{C}(M,x_0)$ bezeichnen wir die Menge aller Elemente $x \in X$,
zu denen ein $\varepsilon_0 > 0$ und eine Abbildung $r(\varepsilon): [0,\varepsilon_0] \to X$ mit den Eigenschaften
$$\lim_{\varepsilon \to 0^+} \frac{\| r(\varepsilon) \|}{\varepsilon} = 0
  \qquad\mbox{und}\qquad x_0 + \varepsilon x + r(\varepsilon) \in M \quad \mbox{ f"ur alle } \varepsilon \in [0,\varepsilon_0]$$
existieren. 
$\mathscr{C}(M,x_0)$ hei"st der lokale Tangentialkegel an die Menge $M$ im Punkt $x_0$.
\end{definition}

\begin{theorem}[Satz von Ljusternik] \label{SatzLjusternik}\index{Satz, Darstellungssatz von Riesz!von Ljusternik@-- von Ljusternik}
Es seien $X$ und $Y$ Banachr"aume,
$V$ eine Umgebung des Punktes $x_* \in X$ und $F$ eine Fr\'echet-differenzierbare Abbildungen der Menge $V$ in $Y$.
Wir setzen voraus,
die Abbildung $F$ sei regul"ar im Punkt $x_*$, d.\,h., es gelte
$${\rm Im\,} F'(x_*) = Y,$$
au"serdem sei ihre Ableitung in diesem Punkt in der gleichm"a"sigen Operatorentopologie des Raumes $L(X,Y)$ stetig. \\
Unter diesen Voraussetzungen stimmt dann der lokale Tangentialkegel an die Menge
$$M= \big\{ x \in V \,\big|\, F(x) = F(x_*) \big\}$$
im Punkt $x_*$ mit dem Kern des Operators $F'(x_*)$ "uberein:
$$\mathscr{C}(M,x_*) = {\rm Ker\,} F'(x_*).$$
\end{theorem}


\subsection{Differenzierbarkeit konkreter Abbildungen}
\begin{beispiel} \label{DiffAbbildung2}
{\rm Sei $\nu(t)=e^{-at}$.
Zu $y_*(\cdot) \in C_{\lim}(\R_+,\R^n;\nu)$ definieren wir die Menge
$$V_{\gamma,\nu} = \{ (t,y) \in \R_+ \times \R^n \,|\, e^{-at}\|y-y_*(t)\| \leq \gamma\}.$$
Wir nehmen an,
dass die Abbildung $g(t,y) : \R \times \R^n \to \R^m$ auf der Menge $V_{\gamma,\nu}$
stetig, stetig differenzierbar bez"uglich $y$ und 
$$\|g_y(t,y)-g_y(t,y')\| \leq C_0 e^{-at} \|y-y'\|$$
mit einem $C_0>0$ f"ur alle $(t,y),(t,y') \in V_{\gamma,\nu}$ ist.
Au"serdem setzen wir voraus,
dass f"ur alle $y(\cdot) \in C_{\lim}(\R_+,\R^n;\nu)$ mit $\|y(\cdot)-y_*(\cdot)\|_{\infty,\nu} \leq \gamma$
die Grenzwerte
$$\lim_{t \to \infty} \nu(t) g\big(t,y(t)\big), \qquad \lim_{t \to \infty} g_y\big(t,y(t)\big)$$
existieren.
Dann ist die Abbildung $G: C_{\lim}(\R_+,\R^n;\nu) \to C_{\lim}(\R_+,\R^m;\nu)$,
$$G\big(y(\cdot)\big)(t) = g\big(t,y(t)\big),$$
im Punkt $y_*(\cdot)$ stetig Fr\'echet-differenzierbar und es gilt
$$\big[G'\big(y_*(\cdot)\big) y(\cdot)\big](t) = g_y\big(t,y_*(t)\big) y(t), \quad t \in \R_+.$$
Denn: Aufgrund unserer Annahmen "uber das Verhalten im Unendlichen stellen $G\big(y(\cdot)\big)$ und $G'\big(y(\cdot)\big)$ Abbildungen des Raumes
$C_{\lim}(\R_+,\R^n;\nu)$ in sich dar. \\
F"ur $t \in \R_+$, $\|y(\cdot)-y_*(\cdot)\|_{\infty,\nu} < \gamma$ und $0 < \lambda <\lambda_0$ gilt
$$\bigg[\frac{G\big(y(\cdot) + \lambda x(\cdot)\big) - G\big(y(\cdot)\big)}{\lambda} - G'\big(y(\cdot)\big)x(\cdot)\bigg](t)
  = \hspace*{-1mm}\int_0^1  \hspace*{-1mm} \big[g_y\big(t,y(t)+s\lambda x(t)\big) - g_y\big(t,y(t)\big) \big] x(t) ds.$$
Nach Voraussetzung an die Abbildung $g_y(t,y)$ existiert auf der Menge $V_{\gamma,\nu}$ ein $C_0>0$ mit 
\begin{eqnarray*}
&& \hspace*{-5mm}\bigg\| \frac{G\big(y(\cdot) + \lambda x(\cdot)\big) - G\big(y(\cdot)\big)}{\lambda}
   - G'\big(y(\cdot)\big)x(\cdot) \bigg\|_{\infty,\nu}
   \leq \sup_{t \in \R_+} \nu(t) \int_0^1 C_0\nu(t)  \|\lambda x(t)\| \|x(t)\| \, ds \\
&&\hspace*{20mm} \leq \sup_{t \in \R_+} \int_0^1 C_0 \lambda \|\nu(t) x(t)\| \, ds \cdot \|x(\cdot)\|_{\infty,\nu} = \lambda C_0,
\end{eqnarray*}
d.\,h. der Grenzwert $\lambda \to 0^+$ konvergiert gleichm"a"sig bez"uglich $\|x(\cdot)\|_{\infty,\nu} = 1$.
Also ist die Abbildung $G$ auf einer Umgebung von $y_*(\cdot)$ Fr\'echet-differenzierbar. 
Weiterhin ergibt sich f"ur die Abbildung $y(\cdot) \to G'\big(y(\cdot)\big)$ bez"uglich der Operatornorm in $y_*(\cdot)$:
\begin{eqnarray*}
    \big\|G'\big(y(\cdot)\big) - G'\big(y_*(\cdot)\big) \big\|
&=& \sup_{\|x(\cdot)\|_{\infty,\nu} =1} \big\| \big[G'\big(y(\cdot)\big) - G'\big(y_*(\cdot)\big)\big] x(\cdot) \big\|_{\infty,\nu} \\
&\leq& \sup_{t \in \R_+} \big\|  g_y\big(t,y(t)\big) - g_y\big(t,y_*(t)\big) \big\|
       \leq C_0 \|y(\cdot)-y_*(\cdot) \|_{\infty,\nu}.
\end{eqnarray*}
Somit ist die stetige Fr\'echet-Differenzierbarkeit nachgewiesen. \hfill $\square$}
\end{beispiel}

\begin{beispiel} \label{DiffZielfunktionalW}
{\rm Es sei $\nu(t)=e^{-at}$ ein Gewicht mit $a>0$ und seien $x_*(\cdot) \in C_{\lim}(\R_+,\R^n;\nu)$, $u_*(\cdot) \in L_\infty(\R_+,\R^m)$.
Ferner nehmen wir an,
die Funktion $f(t,x,u): \R_+ \times \R^n \to \R$ besitzt auf dem Abschluss der Menge
$$U_{\gamma,\nu} = \{ (t,x,u) \in \R_+ \times \R^n \times \R^m \,|\, \nu(t)\|x-x_*(t)\| \leq \gamma, \; \|u-u_*(t)\| \leq \gamma\}$$
folgende Eigenschaften:
\begin{enumerate}
\item[(a)] $f$ ist stetig und stetig differenzierbar bez"uglich $x$ und $u$;
\item[(b)] es existiert f"ur alle $(t,x,u) \in U_{\gamma,\nu}$ eine Funktion $L(\cdot) \in L_1(\R_+,\R;\omega)$ mit
           $$| f(t,x,u)| \leq L(t), \quad \|f_x(t,x,u) \| \leq L(t) \nu(t), \quad \|f_u(t,x,u) \| \leq L(t).$$
\end{enumerate}
Dann ist $J: C_{\lim}(\R_+,\R^n;\nu) \times L_\infty(\R_+,\R^m) \to \R$,
$$J\big(x(\cdot),u(\cdot)\big) = \int_0^\infty \omega(t)f\big(t,x(t),u(t)\big) \, dt,$$
im Punkt $\big(x_*(\cdot),u_*(\cdot)\big)$ Fr\'echet-differenzierbar und es gilt
\begin{eqnarray*}
J'\big(x_*(\cdot),u_*(\cdot)\big) \big(x(\cdot),u(\cdot)\big)
  &=& \int_0^\infty \omega(t)\big\langle f_x\big(t,x_*(t),u_*(t)\big), x(t) \big\rangle \, dt \\
  & & + \int_0^\infty \omega(t)\big\langle f_u\big(t,x_*(t),u_*(t)\big), u(t) \big\rangle \, dt.
\end{eqnarray*}
Nach Voraussetzung (b) erhalten wir auf einer Umgebung des Punktes $\big(x_*(\cdot),u_*(\cdot)\big)$ die Wohldefiniertheit des Zielfunktionals:
$$\big|J\big(x(\cdot),u(\cdot)\big)\big|   \leq \int_0^\infty \omega(t)L(t) \, dt < \infty.$$
Weiterhin ist die lineare Abbildung $J'\big(x_*(\cdot),u_*(\cdot)\big)$ stetig:
\begin{eqnarray*}
       \big|J'\big(x_*(\cdot),u_*(\cdot)\big) \big(x(\cdot),u(\cdot)\big)\big|
&\leq& \int_0^\infty \omega(t) L(t) \big(\nu(t) \|x(t)\| + \|u(t)\|\big) \, dt \\
&\leq& C \big(\|x(\cdot)\|_{\infty,\nu} + \|u(\cdot)\|_{L_\infty}\big).
\end{eqnarray*}
Sei $\varepsilon >0$ gegeben.
Dann k"onnen wir eine Zahl $T>0$ mit der Eigenschaft
$$\int_T^\infty 2 \omega(t) L(t) \,dt \leq \frac{\varepsilon}{3}$$
w"ahlen.
Nach dem Satz von Lusin existiert eine kompakte Teilmenge $K$ von $[0,T]$ derart,
dass $t \to \omega(t) f_x\big(t,x_*(t),u_*(t)\big)$ und $t \to \omega(t) f_u\big(t,x_*(t),u_*(t)\big)$ stetig auf $K$ sind und zudem
die Relation
$$\int_{[0,T] \setminus K} 2 \omega(t) L(t) \, dt \leq \frac{\varepsilon}{3}$$
erf"ullt ist.
Aufgrund der Stetigkeit von $f_x$ und $f_u$ gibt es eine Zahl $\lambda_0>0$ mit
\begin{eqnarray*}
&& \omega(t) \cdot \big\|f_x\big(t,x_*(t)+\lambda x,u_*(t)+\lambda u\big)-f_x\big(t,x_*(t),u_*(t)\big)\big\| 
   \leq \frac{\varepsilon}{3T} \nu(t), \\[1mm]
&& \omega(t) \cdot \big\|f_u\big(t,x_*(t)+\lambda x,u_*(t)+\lambda u\big)-f_u\big(t,x_*(t),u_*(t)\big)\big\|
   \leq \frac{\varepsilon}{3T}
\end{eqnarray*}
f"ur alle $t \in K$, f"ur alle $\nu(t) \|x\| \leq 1$, $\|u \| \leq 1$ und alle $0 < \lambda \leq \lambda_0$.
Zusammen erhalten wir
\begin{eqnarray*}
\lefteqn{\bigg|\frac{J\big(x_*(\cdot)+\lambda x(\cdot),u_*(\cdot)+\lambda u(\cdot)\big)
                - J\big(x_*(\cdot),u_*(\cdot)\big)}{\lambda}
        - J'\big(x_*(\cdot),u_*(\cdot)\big)\big(x(\cdot),u(\cdot)\big)\bigg|} \\
&\leq & \int_0^\infty \omega(t)\bigg| \int_0^1
        \big\langle f_x\big(t,x_*(t)+\lambda s x(t),u_*(t)+\lambda s u(t)\big) - f_x\big(t,x_*(t),u_*(t)\big), x(t) \big\rangle
        ds \bigg|  \, dt \\
&    &  + \int_0^\infty \omega(t) \bigg| \int_0^1 
        \big\langle f_u\big(t,x_*(t)+\lambda s x(t),u_*(t)+\lambda s u(t)\big) - f_u\big(t,x_*(t),u_*(t)\big), u(t) \big\rangle
         ds \bigg| \, dt \\
&\leq& \int_K \frac{\varepsilon}{3T} \big(\nu(t)\|x(t)\|+\|u(t)\|\big) \, dt
      + \int_{[0,T] \setminus K} 2 \omega(t)L(t) \big(\nu(t)\|x(t)\|+\|u(t)\|\big) \, dt \\
&    & + \int_T^\infty 2 \omega(t)L(t) \big(\nu(t)\|x(t)\|+\|u(t)\|\big) \, dt \leq \varepsilon
\end{eqnarray*}
f"ur alle $\|x(\cdot) \|_{\infty,\nu} \leq 1$, $\|u(\cdot) \|_{L_\infty} \leq 1$ und alle
$0 < \lambda \leq \lambda_0$. \hfill $\blacksquare$}
\end{beispiel}

\begin{beispiel} \label{DiffDynamikW1}
{\rm Es sei $\nu(t)=e^{-at}$ mit $a>0$.
Weiterhin seien $x_*(\cdot) \in C_0(\R_+,\R^n;\nu)$ und $u_*(\cdot) \in L_\infty(\R_+,\R^m)$.
Ferner nehmen wir an,
die Abbildung $\varphi(t,x,u): \R_+ \times \R^n \to \R^n$ besitzt auf dem Abschluss der Menge
$$U_{\gamma,\nu} = \{ (t,x,u) \in \R_+ \times \R^n \times \R^m \,|\, \nu(t)\|x-x_*(t)\| \leq \gamma, \; \|u-u_*(t)\| \leq \gamma\}$$
folgende Eigenschaften:
\begin{enumerate}
\item[(a)] $\varphi$ ist stetig und stetig differenzierbar bez"uglich $x$ und $u$;
\item[(b)] es existiert f"ur alle $(t,x,u) \in U_{\gamma,\nu}$ eine Zahl $C_0>0$ mit
           $$\| \varphi(t,x,u) \| \leq C_0(1+\|x\|+\|u\|), \quad \|\varphi_x(t,x,u) \| \leq C_0, \quad \|\varphi_u(t,x,u) \| \leq C_0,$$
\item[(c)] f"ur alle $(t,x,u),(t,x',u') \in U_{\gamma,\nu}$ gelten
           \begin{eqnarray*}
           && \|\varphi_x(t,x',u') - \varphi_x(t,x,u) \| \leq C_0 \big(e^{-at}\|x'-x\| + \|u'-u\|\big), \\
           && \|\varphi_u(t,x',u') - \varphi_u(t,x,u) \| \leq C_0 \big(\|x'-x\| + e^{at} \|u'-u\|\big).
           \end{eqnarray*}
\end{enumerate}
Dann ist $F: C_0(\R_+,\R^n;\nu) \times L_\infty(\R_+,\R^m) \to C_0(\R_+,\R^n;\nu)$,
$$F\big(x(\cdot),u(\cdot)\big)(t) = \int_0^t \varphi\big(s,x(s),u(s)\big) \, ds, \quad t \in \R_+,$$
im Punkt $\big(x_*(\cdot),u_*(\cdot)\big)$ stetig Fr\'echet-differenzierbar und es gilt
\begin{eqnarray*}
    \big[F'\big(x_*(\cdot),u_*(\cdot)\big) \big(x(\cdot),u(\cdot)\big)\big](t)
&=& \int_0^t \varphi_x\big(s,x_*(s),u_*(s)\big) x(s) \, ds \\
& & + \int_0^t \varphi_u\big(s,x_*(s),u_*(s)\big) u(s) \, ds, \quad t \in \R_+.
\end{eqnarray*}
Wir zeigen, dass der Operator $F$ in den Raum $C_0(\R_+,\R^n;\nu)$ abbildet.
Sei $\varepsilon >0$ gegeben.
Wegen $\nu(t)=e^{-at}$, $x(\cdot) \in C_0(\R_+,\R^n;\nu)$, $u(\cdot) \in L_\infty(\R_+,U)$ gilt nach Voraussetzung
$$\lim_{t \to \infty} \nu(t) \big\|\varphi\big(t,x(t),u(t)\big)\big\| = \lim_{t \to \infty} \nu(t) C_0(1+\|x(t)\|+\|u(t)\|) =0.$$
Daher lassen sich Zahlen $0<T<T'$ angeben mit
$$\nu(t) \int_0^T \big\|\varphi\big(t,x(t),u(t)\big)\big\| \, dt \leq \varepsilon$$
f"ur alle $t \geq T'$ und
$$\nu(s) \big\|\varphi\big(s,x(s),u(s)\big)\big\|  \, dt \leq \varepsilon$$
f"ur alle $s \geq T$.
Damit erhalten wir f"ur alle $t \geq T'$
\begin{eqnarray*}
       \nu(t) \big\|F\big(x(\cdot),u(\cdot)\big)(t)\big\|
&\leq& \nu(t) \int_0^T \big\|\varphi\big(s,x(s),u(s)\big)\big\| \, ds + \nu(t) \int_T^t \big\|\varphi\big(s,x(s),u(s)\big)\big\| \, ds \\
&\leq& \varepsilon + \int_T^t \frac{\nu(t)}{\nu(s)} \nu(s) \big\|\varphi\big(s,x(s),u(s)\big)\big\| \, ds \\
&\leq& \varepsilon + \varepsilon \int_T^t e^{-a(t-s)} \, ds \leq \varepsilon + \varepsilon \cdot \frac{1}{a} = K \varepsilon
\end{eqnarray*}
mit einer von $\varepsilon>0$ unabh"angigen Zahl $K$.
Also ergibt sich $\nu(t) \big\|F\big(x(\cdot),u(\cdot)\big)(t)\big\| \to 0$ f"ur $t \to \infty$. 
Auf die gleiche Weise l"asst sich zeigen,
dass $F'\big(x_*(\cdot),u_*(\cdot)\big)$ in den Raum $C_0(\R_+,\R^n;\nu)$ abbildet.
Ferner ist $F'\big(x_*(\cdot),u_*(\cdot)\big)$ stetig, denn 
\begin{eqnarray*}
       \big\| F'\big(x_*(\cdot),u_*(\cdot)\big) \big(x(\cdot),u(\cdot)\big)\big\|_{\infty,\nu}
&\leq& \sup_{t \in \R_+} \int_0^t \frac{\nu(t)}{\nu(s)} C_0 \nu(s)\big(\|x(s)\| + \|u(s)\|\big)\, ds \\
&\leq& C_0 \sup_{t \in \R_+} \int_0^t \frac{\nu(t)}{\nu(s)} \,ds \cdot \big(\|x(\cdot)\|_{\infty,\nu} + \|u(\cdot)\|_{L_\infty}\big).
\end{eqnarray*}
Wir zeigen die  Fr\'echet-Differenzierbarkeit der Abbildung $F$ im Punkt $\big(x_*(\cdot),u_*(\cdot)\big)$.
Dazu bringen wir die Differenz
$$\bigg[\frac{F\big(x_*(\cdot)+\lambda x(\cdot),u_*(\cdot)+\lambda u(\cdot)\big) - F\big(x_*(\cdot),u_*(\cdot)\big)}{\lambda}
   - F'\big(x_*(\cdot),u_*(\cdot)\big) \big(x(\cdot),u(\cdot)\big)\bigg](t)$$
in die Gestalt
\begin{eqnarray*}
&& \int_0^t \bigg[ \int_0^1 \big[
   \varphi_x(\tau,x_*(\tau)+\lambda s x(\tau),u_*(\tau)+\lambda s u(\tau)) - \varphi_x(\tau,x_*(\tau),u_*(\tau))\big]x(\tau) ds \\
&& \quad + \int_0^1 \big[
          \varphi_u(\tau,x_*(\tau)+\lambda s x(\tau),u_*(\tau)+\lambda s u(\tau)) - \varphi_u(\tau,x_*(\tau),u_*(\tau))\big] u(\tau) ds
            \bigg]  \, d\tau.
\end{eqnarray*}
Nach Voraussetzung (c) ergibt sich damit die Ungleichung
\begin{eqnarray*}
&& \bigg\|\bigg[\frac{F\big(x_*(\cdot)+\lambda x(\cdot),u_*(\cdot)+\lambda u(\cdot)\big) - F\big(x_*(\cdot),u_*(\cdot)\big)}{\lambda}
   - F'\big(x_*(\cdot),u_*(\cdot)\big) \big(x(\cdot),u(\cdot)\big)\bigg](t)\bigg\| \\
&& \hspace*{10mm} \leq \int_0^t C_0 \big(\nu(\tau) \|\lambda x(\tau)\| + \| \lambda u(\tau)\|\big) \cdot \|x(\tau)\| \, d\tau \\
&& \hspace*{20mm} + \int_0^t C_0 \big( \|\lambda x(\tau)\| + \nu^{-1}(\tau)\| \lambda u(\tau)\|\big) \cdot \|u(\tau)\| \, d\tau.
\end{eqnarray*}
Damit erhalten wir
\begin{eqnarray*}
&&\hspace*{-10mm} \bigg\|\frac{F\big(x_*(\cdot)+\lambda x(\cdot),u_*(\cdot)+\lambda u(\cdot)\big) - F\big(x_*(\cdot),u_*(\cdot)\big)}{\lambda}
            - F'\big(x_*(\cdot),u_*(\cdot)\big) \big(x(\cdot),u(\cdot)\big)\bigg]\bigg\|_{\infty, \nu} \\
&\leq& \sup_{t \in \R_+} \int_0^t \frac{\nu(t)}{\nu(s)} C_0 \big(\nu(s) \|\lambda x(s)\| + \| \lambda u(s)\|\big) \cdot \nu(s)\|x(s)\| \, ds \\
&&   + \sup_{t \in \R_+} \int_0^t \frac{\nu(t)}{\nu(s)} C_0 \big(\nu(s) \|\lambda x(s)\| + \| \lambda u(s)\|\big) \cdot \|u(s)\| \, ds
   \leq 2C_0 \sup_{t \in \R_+} \int_0^t \frac{\nu(t)}{\nu(s)} \, ds \cdot \lambda\leq \varepsilon
\end{eqnarray*}
f"ur alle $\|x(\cdot)\|_{\infty,\nu} \leq 1$, $\|u(\cdot)\|_{L_\infty} \leq 1$ und alle $0 < \lambda \leq \lambda_0$. \\
Mit den gleichen Argumenten,
mit denen eben die Fr\'echet-Differenzierbarkeit der Abbildung $F$ im Punkt $\big(x_*(\cdot),u_*(\cdot)\big)$ gezeigt wurde,
lassen sich die Differenzierbarkeit auf einer Umgebung des Punktes $\big(x_*(\cdot),u_*(\cdot)\big)$ im Raum
$C_0(\R_+,\R^n;\nu) \times L_\infty(\R_+,\R^m)$ und die Stetigkeit der Fr\'echet-Ableitung im Punkt $\big(x_*(\cdot),u_*(\cdot)\big)$ zeigen. \hfill $\square$}
\end{beispiel}


\subsection{Zu absolutstetigen Funktionen "uber $\R_+$}
\begin{lemma} \label{LemmaGrenzwert1} 
Sei $x(\cdot) \in W^1_1(\R_+,\R^n)$.
Dann gilt nach Magill \cite{Magill}:
$$\lim_{t \to \infty} \|x(t)\|=0.$$
\end{lemma}

{\bf Beweis} Da $x(\cdot) \in W^1_1(\R_+,\R^n)$ ist, gilt die Darstellung
$$x(t)=x(0)+\int_0^t \dot{x}(s)\,ds, \quad t \in \R_+.$$
Darin ist der Integralterm "uber $\R_+$ absolut integrierbar.
Also besitzt $x(t)$ einen Grenzwert f"ur $t \to \infty$.
Dieser muss gleich Null sein, da $x(\cdot)$ "uber $\R_+$ integrierbar ist. \hfill $\blacksquare$ \\[2mm]
Es sei $\mu$ ein $\sigma$-endliche Ma"s, dass bez"uglich dem Lebesgue-Ma"s die Dichte $\nu(\cdot)$ besitzt.
Dann ergeben sich aus den Lebesgue- und Sobolevr"aume bez"uglich dem Ma"s $\mu$ gewichteten Lebesgue- und Sobolevr"aume zum Gewicht $\nu$,
die wir mit $L_p(I,\R^n;\nu)$ bzw. $W^1_p(I,\R^n;\nu)$ bezeichnen.
F"ur $p=2$ lautet die gewichtete Norm des Raumes $L_2(I,\R^n;\nu)$
$$\|x(\cdot)\|^2_{L_2(\nu)}= \int_0^\infty \|x(t)\|^2 \nu(t) \,dt$$
und die gewichtete Norm des Raumes $W^1_2(I,\R^n;\nu)$
$$\|x(\cdot)\|_{W^1_2(\nu)} = \|x(\cdot)\|_{L_2(\nu)} + \|\dot{x}(\cdot)\|_{L_2(\nu)}.$$
Es sei im Weiteren $\nu(\cdot) \in W^1_1(\R_+,\R)$ eine Gewichtsfunktion,
zu der eine Konstante $K>0$ mit $|\dot{\nu}(t)| \leq K \nu(t)$ f"ur alle $t \in \R_+$ existiert.

\begin{lemma} \label{LemmaGrenzwert2}
Es sei $x(\cdot) \in W^1_2(\R_+,\R^n;\nu)$.
Wir betrachten die Funktionen
$$f(t)=\|x(t)\|^2 \nu(t), \qquad \psi(t)=\nu(t)x(t).$$
Dann gelten $f(\cdot) \in W^1_1(\R_+,\R)$, $\psi(\cdot) \in W^1_1(\R_+,\R^n)$ und damit nach Lemma \ref{LemmaGrenzwert1}:
$$\lim_{t \to \infty} f(t)=0, \qquad \lim_{t \to \infty} \|\psi(t)\|=0.$$
\end{lemma}

{\bf Beweis} Wir schreiben $f(\cdot)$ in der Form $f(t)=\langle x(t), x(t) \rangle \nu(t)$.
Dann erhalten wir unmittelbar die verallgemeinerte Ableitung
$$\dot{f}(t)=2\langle x(t), \dot{x}(t) \rangle \nu(t) + \|x(t)\|^2 \dot{\nu}(t).$$
Bei Anwendung der Beziehungen
$$\big|\langle x(t), \dot{x}(t) \rangle \nu(t)\big|^2 \leq \|x(t)\|^2 \nu(t) \cdot \|\dot{x}(t)\|^2 \nu(t)$$
ergibt sich mit der Cauchy-Schwarzschen Ungleichung:
$$\|\dot{f}(\cdot)\|_{L_1}+\|f(\cdot)\|_{L_1}
  \leq 2 \cdot \|\dot{x}(\cdot)\|_{L_2(\nu)} \cdot \|x(\cdot)\|_{L_2(\nu)} + (1+K) \cdot \|x(\cdot)\|^2_{L_2(\nu)}.$$
Die Funktion $\psi(\cdot)$ besitzt eine verallgemeinerte Ableitung und wir erhalten wegen der Eigenschaften der Gewichtsfunktion $\nu(\cdot)$ und
mit der Cauchy-Schwarzschen Ungleichung
$$\|\psi(\cdot)\|_{L_1} \leq \|\nu(\cdot)\|^{1/2}_{L_1} \cdot \|x(\cdot)\|_{L_2(\nu)}, \quad
  \|\dot{\psi}(\cdot)\|_{L_1} \leq (1+K) \|\nu(\cdot)\|^{1/2}_{L_1} \cdot \|x(\cdot)\|_{W^1_2(\nu)}.$$
Damit sind alle Behauptungen gezeigt. \hfill $\blacksquare$

\begin{lemma} \label{LemmaGrenzwert3}
Es seien $x(\cdot) \in W^1_2(\R_+,\R^n;\nu)$ und $y(\cdot) \in W^1_2(\R_+,\R^n;\nu^{-1})$.
Wir setzen:
$$g(t)=\|y(t)\|^2 \nu^{-1}(t), \qquad h(t)= \langle x(t),y(t) \rangle.$$
Dann gelten $g(\cdot), h(\cdot) \in W^1_1(\R_+,\R)$ und nach Lemma \ref{LemmaGrenzwert1}:
$$\lim_{t \to \infty} g(t)=0, \qquad \lim_{t \to \infty} h(t)=0.$$
\end{lemma}

{\bf Beweis} Wir bringen $g(\cdot)$ in die Gestalt $g(t)=\langle y(t), y(t) \rangle \nu^{-1}(t)$.
Dann erhalten wir die verallgemeinerte Ableitung
$$\dot{g}(t)=2\langle y(t), \dot{y}(t) \rangle \nu^{-1}(t) + \|y(t)\|^2 \cdot \frac{d}{dt}\nu^{-1}(t).$$
Bei Anwendung der Ungleichung
$$\big|\langle y(t), \dot{y}(t) \rangle \nu^{-1}(t)\big|^2 \leq \|y(t)\|^2 \nu^{-1}(t) \cdot \|\dot{y}(t)\|^2 \nu^{-1}(t)$$
ergibt sich mit der Cauchy-Schwarzschen Ungleichung:
$$\|\dot{g}(\cdot)\|_{L_1}+\|g(\cdot)\|_{L_1} \leq 2 \cdot \|\dot{y}(\cdot)\|_{L_2(\nu^{-1})} \cdot \|y(\cdot)\|_{L_2(\nu^{-1})}
         + (1+K) \cdot \|y(\cdot)\|^2_{L_2(\nu^{-1})}.$$
Die Funktion $h(\cdot)$ besitzt eine verallgemeinerte Ableitung und es gilt
\begin{eqnarray*}
\|h(\cdot)\|_{L_1} &\leq& \|y(\cdot)\|_{L_2(\nu^{-1})} \cdot \|x(\cdot)\|_{L_2(\nu)}, \\
\|\dot{h}(\cdot)\|_{L_1} &\leq& \|\dot{y}(\cdot)\|_{L_2(\nu^{-1})} \cdot \|x(\cdot)\|_{L_2(\nu)} 
                                  + \|y(\cdot)\|_{L_2(\nu^{-1})} \cdot \|\dot{x}(\cdot)\|_{L_2(\nu)}.
\end{eqnarray*}
Damit ist Lemma \ref{LemmaGrenzwert3} gezeigt. \hfill $\blacksquare$
\lhead[\thepage \hspace*{1mm} Differentialgleichungen]{}
\rhead[]{Differentialgleichungen \hspace*{1mm} \thepage}
\section{Lineare Differentialgleichungen}
Wir betrachten f"ur $A(t) \in \R^{n \times n}$ und $a(t) \in \R^n$ das lineare Differentialgleichungssystem
$$\dot{x}(t)=A(t)x(t)+a(t).$$

\begin{lemma} \label{LemmaDGL2}
Es seien die Abbildung $t \to A(t)$ und die Vektorfunktion $a(\cdot)$ "uber $\R_+$ me"sbar und beschr"ankt.
Ferner gelte mit der Funktion $\nu(t)=e^{-at}$, $a>0$, die Bedingung
$$\sup_{t \in \R_+} \int_0^t \frac{\nu(t)}{\nu(s)} \|A(s)\| \, ds = q < 1,$$
die z.\,B. f"ur $a > \|A(\cdot)\|_{L_\infty}$ erf"ullt ist.
Dann existiert zu jedem $\zeta(\cdot) \in C_0(\R_+,\R^n;\nu)$ und jedem $\tau \in \R_+$
eine eindeutig bestimmte Vektorfunktion $x(\cdot) \in C_0(\R_+,\R^n;\nu)$ derart,
dass
$$\zeta(t)=x(t) - \int_\tau^t  \big[A(s)x(s)+a(s)\big] \, ds$$
f"ur alle $t \in \R_+$ erf"ullt ist.
\end{lemma}

{\bf Beweis} Wir werden im Folgenden zeigen,
dass die Fixpunktgleichung $x(\cdot) = T\big(x(\cdot)\big)$, wobei der Operator $T$ durch
$$x(\cdot) \to T\big(x(\cdot)\big), \quad
  T\big(x(\cdot)\big)(t) = \zeta(t) + \int_\tau^t \big[A(s)x(s)+a(s)\big] \, ds, \quad t \in \R_+,$$
gegeben wird, stets eine eindeutige L"osung besitzt.
WZun"achst bildet der Operator $T$ den Raum $C_0(\R_+,\R^n;\nu)$ in sich ab:
Wegen $\zeta(\cdot), x(\cdot) \in C_0(\R_+,\R^n;\nu)$ und da $a(\cdot)$ beschr"ankt ist, gelten die Grenzwerte
$$\lim_{t \to \infty} \nu(t) \zeta(t)=0, \quad \lim_{t \to \infty} \nu(t) x(t)=0, \quad
  \lim_{t \to \infty} \nu(t) \int_0^t \|a(s)\| \, ds =0.$$
Sei weiterhin $\varepsilon >0$ gegeben.
Dann lassen sich Zahlen $0<T<T'$ angeben mit
$$\nu(t) \|x(t)\| \leq \varepsilon \mbox{ f"ur alle } t \geq T, \qquad \nu(t) \int_0^T \|A(s) x(s)\| \, ds  \leq \varepsilon \mbox{ f"ur alle } t \geq T'.$$
Damit erhalten wir f"ur alle $t \geq T'$
\begin{eqnarray*}
&& \nu(t) \int_0^t \|A(s) x(s)\| \, ds
   = \nu(t) \int_0^T \|A(s) x(s)\| \, ds + \nu(t) \int_T^t \|A(s) x(s)\| \, ds \\
&& \hspace*{10mm} \leq \varepsilon + \|A(\cdot)\|_{L_\infty} \int_T^t \frac{\nu(t)}{\nu(s)} [\nu(s) \|x(s)\|] \, ds
       \leq \varepsilon + \varepsilon \cdot \frac{1}{a} \|A(\cdot)\|_{L_\infty}.
\end{eqnarray*}
Zusammen ergibt sich damit $T\big(x(\cdot)\big)(t) \to 0$ f"ur $t \to \infty$. \\
Wir zeigen nun, dass der Operator $T$ kontraktiv ist:
Bei mehrfacher Anwendung des Operators $T$ ergeben sich f"ur $x_1(\cdot),x_2(\cdot) \in C_0(\R_+,\R^n;\nu)$ die Beziehungen
\begin{eqnarray*}
\lefteqn{\nu(t)\big\| \big[T\big(x_1(\cdot) - x_2(\cdot)\big)\big](t) \big\|
         \leq \nu(t)\int_\tau^t \|A(s)\| \| x_1(s) - x_2(s) \|\, ds} \\
&& \hspace*{10mm} \leq \int_0^t \frac{\nu(t)}{\nu(s)} \big\|A(s)\big\| \cdot \nu(s)\| x_1(s) - x_2(s) \|\, ds
       \leq q \cdot \big\| x_1(\cdot) - x_2(\cdot) \big\|_{\infty,\nu}, \\
\lefteqn{\nu(t)\big\| \big[T^2\big(x_1(\cdot) - x_2(\cdot)\big)\big](t) \big\|
         \leq \nu(t)\int_\tau^t \|A(s)\| \big\| \big[T\big(x_1(\cdot) - x_2(\cdot)\big)\big](s) \big\| \, ds} \\
&& \hspace*{10mm} \leq \int_0^t \frac{\nu(t)}{\nu(s)} \|A(s)\| \cdot \nu(s) \big\|\big[ T\big(x_1(\cdot) - x_2(\cdot)\big)\big](s) \big\| \, ds
       \leq q^2 \cdot \big\| x_1(\cdot) - x_2(\cdot) \big\|_{\infty,\nu}.
\end{eqnarray*}
Sukzessive erhalten wir f"ur $m \in \N$:
\begin{eqnarray*}
\lefteqn{\nu(t) \big\| \big[T^m \big(x_1(\cdot) - x_2(\cdot)\big)\big](t) \big\|
         \leq \nu(t) \int_\tau^t \|A(s)\| \big\| \big[T^{m-1}\big(x_1(\cdot) - x_2(\cdot)\big)\big](s) \big\| \, ds} \\
&& \hspace*{10mm} \leq \int_0^t \frac{\nu(t)}{\nu(s)} \|A(s)\| \cdot \nu(s) \big\| \big[T^{m-1}\big(x_1(\cdot) - x_2(\cdot)\big)\big](s) \big\| \, ds
       \leq q^m \| x_1(\cdot) - x_2(\cdot) \|_{\infty,\nu}.
\end{eqnarray*}
In der Topologie des Raumes $C_0(\R_+,\R^n;\nu)$ gilt daher
$$\big\|T^m \big(x_1(\cdot) - x_2(\cdot)\big) \big\|_{\infty,\nu} \leq q^m \cdot \| x_1(\cdot) - x_2(\cdot) \|_{\infty,\nu}.$$
Nach Voraussetzung ist $q<1$.
Daher existiert nach dem Banachschen Fixpunktsatz genau ein $x(\cdot)$ mit $x(\cdot) = T x(\cdot)$.
\hfill $\blacksquare$
\lhead[\thepage \hspace*{1mm} Konvexe Analysis]{}
\rhead[]{Konvexe Analysis \hspace*{1mm} \thepage}
\section{Elemente der Konvexen Analysis} \label{AbschnittKonvexeAnalysis}
Bei der Zusammenstellung der grundlegenden Ergebnisse beschr"anken wir uns auf die Eigenschaften konvexer und lokalkonvexer Funktionen
nach Clarke \cite{Clarke}, Ioffe \& Tichomirov \cite{Ioffe} und Rockafellar \cite{Rockafellar}.
Im vorliegenden Rahmen stimmen die klassische Richtungsableitung und der Clarkesche Gradient "uberein.
Deswegen verweisen wir bez"uglich Lemma \ref{LemmaRichtungsableitung} und bez"uglich der Kettenregel \ref{SatzKettenregel}
auf \cite{Ioffe}.

\subsection{Das Subdifferential konvexer Funktionen}
Es seien $X,Y$ Banachr"aume.
Eine Funktion $f$ auf $X$ ist in der Konvexen Analysis eine Abbildung in die erweiterte reelle Zahlengerade, d.\,h.
$f:X \to \overline{\R} = [-\infty,\infty]$.
Der effektive Definitionsbereich
der Abbildung $f$ ist die Menge ${\rm dom\,}f = \{ x \in X | f(x) < \infty\}$. 
Die Funktion $f$ hei"st eigentlich\index{Funktion, absolutstetige!eigentliche@--, eigentliche},
falls ${\rm dom\,}f \not= \emptyset$ und $f(x) > -\infty$ f"ur alle $x \in X$ gelten. \\
Die eigentliche Funktion $f$ hei"st konvex\index{Funktion, absolutstetige!konvexe@--, konvexe},
wenn f"ur alle $x_1,x_2 \in X$ und alle $0 \leq \alpha \leq 1$ folgende Ungleichung gilt:
$$f\big(\alpha x_1 + (1-\alpha) x_2 \big) \leq \alpha f(x_1) + (1-\alpha) f(x_2).$$
Die Funktion $f$ hei"st homogen\index{Funktion, absolutstetige!homogene@--, homogene},
falls $f(0) = 0$ und $f(\lambda x) = \lambda f(x)$ f"ur alle $x \in X,\, \lambda > 0$ ist.
Eine eigentliche konvexe Funktion ist genau dann in einem Punkt stetig,
wenn sie auf einer Umgebung dieses Punktes nach oben beschr"ankt ist.
In diesem Fall ist das Innere des effektiven Definitionsbereichs nichtleer.
Ist andererseits eine homogene Funktion auf einer Umgebung des Nullpunktes stetig, so ist sie auf $X$ stetig. \\
Ist $f$ eine eigentliche konvexe Funktion auf $X$,
dann existiert in jedem Punkt der Menge ${\rm dom\,}f$ die klassische
Richtungsableitung\index{Ableitung, Fr\'echet-!Richtungs@--, Richtungs-}\label{Richtungsableitung},
d.\,h. f"ur alle $z \in X$ der Grenzwert
$$f'(x;z)=\lim_{\lambda \to 0^+} \frac{f(x + \lambda z) - f(x)}{\lambda}.$$
Sei $f$ eigentlich, konvex und in $x$ stetig.
Dann ist $f$ auf einer Umgebung des Punktes $x$ nach oben beschr"ankt,
in $x$ lokal Lipschitz-stetig und es existiert der Clarkesche Gradient\index{Clarkescher Gradient}\label{ClarkescherGradient}
$$f^\circ(x;z) = \limsup_{\substack{y \to x \\ \lambda \to 0^+}} \frac{f(y + \lambda z) - f(y)}{\lambda}.$$
Unter diesen Voraussetzungen ist au"serdem die Funktion $f$ ist im Punkt $x$ regul"ar\index{Funktion, absolutstetige!regul"are@--, regul"are}
im Sinn der Konvexen Analysis,
d.\,h. $f'(x;\cdot)=f^\circ(x;\cdot)$.
Das Subdifferential\index{Subdifferential} der eigentlichen konvexen Funktion $f$ besteht im Punkt $x$ aus allen Subgradienten
$x^* \in X^*$, d.\,h. \label{Subdifferential}
$$\partial f(x) = \{ x^* \in X^* | f(z) - f(x) \geq \langle x^*, z-x \rangle \mbox{ f"ur alle } z \in X \}.$$
F"ur eine eigentliche konvexe Funktion $f$ gilt $\partial f(x)=\partial f'(x;0)$ f"ur alle $x \in {\rm dom\,}f$.
Ist $f$ eine eigentliche homogene konvexe Funktion und $x \not= 0$, dann ist
$$\partial f(x) = \{ x^* \in \partial f(0) | f(x) = \langle x^*, x \rangle \}.$$
Mit $\partial_x f(x,y)$ bezeichnen wir das Subdifferential der Abbildung $x \to f(x,y)$.\label{partSubdifferential}

\subsection{Lokalkonvexe Funktionen}
Es seien $X,Y$ Banachr"aume.
Eine auf $X$ definierte Funktion $G$ hei"st im Punkt $x_0$ lokalkonvex\index{Funktion, absolutstetige!lokalkonvexe@--, lokalkonvexe},
wenn ihre Richtungsableitung in diesem Punkt existiert und $x \to G'(x_0;x)$ konvex ist.
Im Folgenden seien $g: X \to Y$ im Punkt $x_0 \in X$ Fr\'echet-differenzierbar und
$f: Y \to \overline{\R}$ eigentlich, konvex und im Punkt $g(x_0)$ stetig.

\begin{lemma} \label{LemmaRichtungsableitung}
Die Funktion $G: X \to \overline{\R}$, $G(x) = f\big( g(x)\big)$, besitzt in $x_0$ eine klassische Richtungsableitung, es gilt
$$G'(x_0;x) = f'\big(g(x_0);g'(x_0)x\big)$$
und die Richtungsableitung konvergiert bez"uglich jeder Richtung $x$
gleichm"a"sig\index{Funktion, absolutstetige!gleichmassig@--, gleichm"a"sig differenzierbare}:
$$\bigg| \frac{G(x_0+\lambda z) - G(x_0)}{\lambda} - G'(x_0;x) \bigg| < \varepsilon \quad
  \mbox{f"ur alle } z \in U(x), \; 0<\lambda<\lambda_0.$$
Insbesondere folgt aus der gleichm"a"sigen Differenzierbarkeit bez"uglich aller Richtungen,
dass die Richtungsableitung der Abbildung $G$ eine stetige Funktion ist.
\end{lemma}

\begin{lemma} \label{LemmaClarke}
Die Funktion $G(x) = f\big(g(x)\big)$ ist in $x_0$ regul"ar.
\end{lemma}
{\bf Beweis} Nach Definition des $\limsup$ existieren Folgen $x_n \to x_0$ und $\lambda_n \to 0^+$ mit
$$\lim_{n \to \infty} \frac{G(x_n + \lambda_n x) - G(x_n)}{\lambda_n} = G^\circ(x_0;x).$$
Unter den getroffenen Voraussetzungen gilt $f'\big(g(x_0);g'(x_0)x\big) = f^\circ\big(g(x_0);g'(x_0)x\big)$ und wir erhalten 
\begin{eqnarray*}
G^\circ(x_0;x) &=& \lim_{n \to \infty} \frac{G(x_n + \lambda_n x) - G(x_n)}{\lambda_n}
                   = \lim_{n \to \infty} \frac{f\big(g(x_n + \lambda_n x)\big) - f\big(g(x_n)\big)}{\lambda_n} \\
               &\leq& \limsup_{\substack{y \to x \\ \lambda \to 0^+}} \frac{f\big(g(y + \lambda x)\big) - f\big(g(y)\big)}{\lambda} 
                   = f'\big(g(x_0);g'(x_0)x\big) = G'(x_0;x).
\end{eqnarray*}
Andererseits folgt unmittelbar die Relation
\begin{eqnarray*}
G'(x_0;x) &=& \lim_{\lambda \to 0^+} \frac{G(x_0 + \lambda x) - G(x_0)}{\lambda}
            \leq \limsup_{\substack{y \to x \\ \lambda \to 0^+}} \frac{G(y + \lambda x) - G(y)}{\lambda} = G^\circ(x_0;x).
\end{eqnarray*}
Beide Ungleichungen zeigen $G'(x_0;\cdot)=G^\circ(x_0;\cdot)$. \hfill $\blacksquare$

\begin{satz}[Kettenregel] \label{SatzKettenregel}
Es seien $g: X \to Y$ im Punkt $x_0 \in X$ Fr\'echet-differenzierbar und $f: Y \to \overline{\R}$ eigentlich,
konvex und im Punkt $g(x_0)$ stetig.
Dann ist die Funktion $G(x) = f\big( g(x)\big)$ im Punkt $x_0$ regul"ar und es gilt
$$\partial G(x_0) = g'^*(x_0) \partial f\big(g(x_0)\big).$$
\end{satz}


\subsection{Das Subdifferential konkreter Funktionen}
\begin{beispiel} \label{SubdifferentialMaximum3}
{\rm Wir betrachten im Raum $C_{\lim}(\R_+,\R)$ die Funktion
$$f\big(x(\cdot)\big) = \max_{t \in \overline{\R}_+} x(t).$$
Diese Funktion ist stetig, konvex und homogen.
Mit den gleichen Argumenten wie im Fall des Raumes der stetigen Funktionen "uber einem kompakten Intervallergibt sich,
dass das Subdifferential $\partial f(0)$ aus denjenigen signierten regul"aren Borelschen Ma"sen $\mu$ auf $\overline{\R}_+$ besteht,
die nichtnegativ sind und die Totalvariation $\|\mu\| = 1$ besitzen. \\
Weiterhin besteht das Subdifferential der Funktion $f$ in einem vom Nullpunkt verschiedenen Punkt $x(\cdot)$ aus den regul"aren
Borelschen Ma"sen auf $\overline{\R}_+$, deren Norm gleich Eins ist und die auf der Menge
$T = \big\{ t \in \overline{\R}_+ \,\big|\, x(t) = f\big(x(\cdot)\big) \big\}$
konzentriert sind. \hfill $\square$}
\end{beispiel}

\begin{beispiel} \label{SubdifferentialMaximum4}
{\rm Es sei $g(t,x)$ eine Funktion auf $\R_+ \times \R^n$,
die bez"uglich beider Ver"anderlicher gleichm"a"sig stetig und f"ur jedes $t \in \R_+$ nach $x$ gleichm"a"sig stetig differenzierbar ist.
Dann ist die Abbildung
$$\tilde{g}: C_{\lim}(\R_+,\R^n) \to C_{\lim}(\R_+,\R), \qquad \big[\tilde{g}\big(x(\cdot)\big)\big] (t) = g\big(t,x(t)\big),\quad t \in \R_+,$$
Fr\'echet-differenzierbar und es gilt
$$\big[\tilde{g}'\big(x(\cdot)\big) z(\cdot)\big] (t) = \big\langle g_x \big(t,x(t)\big), z(t) \big\rangle, \qquad t \in \R_+.$$
Weiterhin ist die Funktion $f$ im Beispiel \ref{SubdifferentialMaximum3} auf $C_{\lim}(\R_+,\R)$ stetig.
Auf die Funktion
$$G\big(x(\cdot)\big) = f\Big(g\big(x(\cdot)\big)\Big) =\max_{t \in \overline{\R}_+} g\big(t,x(t)\big).$$
wir die Kettenregel an:
Das Subdifferential der Funktion $G$ besteht im Punkt $x(\cdot)$ genau aus denjenigen stetigen linearen Funktionalen $x^*$,
die die Darstellung
$$\big\langle x^*,z(\cdot) \big\rangle = \int_0^\infty \big\langle g_x \big(t,x(t)\big), z(t) \big\rangle \,d\mu(t) + 
                                         \lim_{t \to \infty} \big\langle g_x \big(t,x(t)\big), z(t) \big\rangle \, \mu(\{\infty\})$$
besitzen,
wobei das regul"are Borelsches Ma"s $\mu$ auf $T= \big\{ t \in \overline{\R}_+ \,\big|\, g\big(t,x(t)\big) = G\big(x(\cdot)\big) \big\}$
konzentriert ist und $\|\mu\| =1$ gilt. \hfill $\square$}
\end{beispiel}

\begin{beispiel} \label{SubdifferentialMaximum5}
{\rm Im Raum $C_0(\R_+,\R)$ betrachten wir die Funktion $f\big(x(\cdot)\big) = \sup\limits_{t \in \R_+} x(t)$.
Sie ist stetig, homogen und konvex.
F"ur das Subdifferential $\partial f(0)$ liefert die Ungleichung
$$\sup_{t \in \R_+} x(t) \geq \int_0^\infty x(t) \, d\mu(t) \quad \mbox{ f"ur alle }x(\cdot) \in C_0(\R_+,\R),$$
dass das "uber $\R_+$ signierte regul"are Borelsche Ma"s $\mu$ nichtnegativ ist.
Aber im Gegensatz zu Beispiel \ref{SubdifferentialMaximum3} erhalten wir aus der Ungleichungskette 
$$\sup_{t \in \R_+} x(t) \geq \int_0^\infty x(t) \, d\mu(t) \geq \inf_{t \in \R_+} x(t)$$
lediglich $\|\mu\|\leq 1$, da die Funktionen $x(\cdot)$ im Unendlichen stets verschwinden. \\
Umgekehrt gilt, wenn $\mu \geq 0$ und $\|\mu\|\leq 1$, die Ungleichung
$$\sup_{t \in \R_+} x(t) \geq \sup_{t \in \R_+} x(t) \cdot \int_0^\infty d\mu(t) = \int_0^\infty \sup_{t \in \R_+} x(t) \,d\mu(t)
  \geq \int_0^\infty x(t) \,d\mu(t)$$
f"ur alle $x(\cdot) \in C_0(\R_+,\R)$.
Damit besteht das Subdifferential der Funktion $f$ in $x(\cdot)=0$ aus allen regul"aren Borelschen Ma"sen $\mu \geq 0$ mit $\|\mu\| \leq 1$. \\
F"ur $x(\cdot) \not= 0$ m"ussen wir die Unterscheidung treffen,
ob die Funktion $x(\cdot) \in C_0(\R_+,\R)$ ein Maximum besitzt oder nicht.
Die Funktion $x(\cdot)$ besitzt genau dann kein Maximum "uber $\R_+$,
wenn $x(t)<0$ f"ur alle $t$ gilt.
In diesem Fall erhalten wir
$$0 = \sup_{t\in \R_+} x(t) = \int_{\R_+} x(t)\, d\mu(t) \qquad\Leftrightarrow\qquad \|\mu\|=0.$$
Nimmt die Funktion $x(\cdot)$ ihr Maximum "uber $\R_+$ an,
dann gilt $\|\mu\|=1$ und $\mu$ ist auf der Menge $T=\big\{ t \in \R_+ \,\big|\, x(t) = f\big(x(\cdot)\big) \big\}$ konzentriert. \hfill $\square$}
\end{beispiel}

\begin{beispiel} \label{SubdifferentialMaximum6}
{\rm Unter den Voraussetzungen des Beispiels \ref{SubdifferentialMaximum4} ist die Abbildung $\tilde{g}$
im Rahmen des Raumes $C_0(\R_+,\R^n)$ Fr\'echet-differenzierbar.
Somit k"onnen wir zur Berechnung des Subdifferentials der Funktion $G\big(x(\cdot)\big) = f\Big(g\big(x(\cdot)\big)\Big)$
die Kettenregel (Satz \ref{SatzKettenregel}) anwenden
und erhalten f"ur $x^* \in \partial G\big(x(\cdot)\big)$:
$$\langle x^*,z(\cdot) \rangle = \int_0^\infty \big\langle g_{jx}\big(t,x(t)\big),z(t) \big\rangle d\mu(t).$$
F"ur das Ma"s $\mu$ m"ussen dabei folgende F"alle unterschieden werden:
\begin{enumerate}
\item[(a)] Ist $\tilde{g}\big(x(\cdot)\big)=0$, so besitzt das Ma"s $\mu$ eine Totalvariation $\|\mu\|\leq 1$.
\item[(b)] Ist $\tilde{g}\big(x(\cdot)\big)\not=0$ und besitzt kein Maximum "uber $\R_+$, so ist $\|\mu\|=0$.
\item[(c)] Ist $\tilde{g}\big(x(\cdot)\big) \not=0$ und besitzt ein Maximum "uber $\R_+$,
           so ist $\|\mu\| = 1$ und $\mu$ ist auf der Menge
           $T = \big\{ t \in \R_+ \,\big|\, g\big(t,x(t)\big) = G\big(x(\cdot)\big) \big\}$ konzentriert. \hfill $\square$
\end{enumerate}}
\end{beispiel}
\lhead[\thepage \hspace*{1mm} Theorie der Extremalaufgaben]{}
\rhead[]{Theorie der Extremalaufgaben \hspace*{1mm} \thepage}
\section{Ein Extremalprinzip f"ur ein schwaches lokales Minimum} \label{AbschnittEPschwach}
In der folgenden Darstellung des Extremalprinzips f"ur ein schwaches lokales Minimum seien $X$, $Y$, $Z$ Banachr"aume,
$\mathscr{U}$ ein normierter Raum, $U \subseteq \mathscr{U}$,
$\mathscr{K} \subseteq Z$ ein abgeschlossener konvexer Kegel mit Spitze in Null und es sei ${\rm int\,}\mathscr{K} \not= \emptyset$.
Weiterhin seien $J$ ein Funktional auf $X \times U$,
$\mathscr{F}$ eine Abbildung des Produktes $X \times U$ in den Raum $Y$ und $G:X \to Z$. \\[2mm]
Unter diesen Angaben betrachten wir in diesem Abschnitt die Extremalaufgabe
\begin{equation}\label{EAAnhang2}
J(x,u) \to \inf; \quad \mathscr{F}(x,u)=0, \quad G(x) \in \mathscr{K}, \quad x \in X,\; u \in U,\; U \mbox{ konvex}.
\end{equation}
Der Punkt $(x,u)$ ist ein zul"assiges Element der Aufgabe (\ref{EAAnhang2}),
falls s"amtliche Nebenbedingungen erf"ullt sind.
Der Punkt $(x_*,u_*)$ hei"st ein schwaches\index{Minimum, schwaches lokales} lokales Minimum der Extremalaufgabe (\ref{EAAnhang2}),
wenn ein $\varepsilon >0$ derart existiert,
dass f"ur alle zul"assigen Paare $(x,u)$ mit $\|(x,u)-(x_*,u_*)\|_{X \times \mathscr{U}} \leq \varepsilon$
die Ungleichung $J(x_*,u_*) \leq J(x,u)$ gilt. \\[2mm]
Auf $X \times U \times \R \times Y^* \times Z^*$ definieren wir zur Aufgabe (\ref{EAAnhang2}) die Lagrange-Funktion
$$\mathscr{L}(x,u,\lambda_0,y^*,z^*)
  = \lambda_0 J(x,u) + \langle y^*, \mathscr{F}(x,u) \rangle + \langle z^*, G(x) \rangle.$$
Au"serdem bezeichnet $\mathscr{K}^*$ den dualen Kegel
$$\mathscr{K}^*=\{z^* \in Z^* \,|\, \langle z^*,z \rangle \leq 0 \mbox{ f"ur alle } z \in \mathscr{K} \}.$$

\begin{theorem}[Extremalprinzip] \label{SatzExtremalprinzipSchwach} \index{Extremalprinzip}
Sei $(x_*,u_*)$ ein zul"assiges Element der Aufgabe (\ref{EAAnhang2}).
\begin{enumerate}
\item[(A)] Wir nehmen an, dass der Punkt $(x_*,u_*)$ eine Umgebung $V$ mit folgenden Eigenschaften besitzt:
           \begin{enumerate}
           \item[(A$_1$)] Die Funktion $J(x,u)$ ist im Punkt $(x_*,u_*)$ Fr\'echet-differenzierbar;
           \item[(A$_2$)] Die Abbildung $\mathscr{F}(x,u)$ ist auf der Umgebung $V$ Fr\'echet-differenzierbar
                          und im Punkt $(x_*,u_*)$ stetig Fr\'echet-differenzierbar.
           \item[(A$_3$)] Die Abbildung $G(x)$ ist im Punkt $x_*$ Fr\'echet-differenzierbar.
           \end{enumerate}
\item[(B)] Weiterhin setzen wir voraus, dass der Operator $\mathscr{F}_x(x_*,u_*)$ eine endliche Kodimension besitzt.
\end{enumerate}
Ist dann $(x_*,u_*)$ schwache lokale Minimalstelle der Aufgabe (\ref{EAAnhang2}),
so ist f"ur die Aufgabe (\ref{EAAnhang2}) das Lagrangesche Prinzip g"ultig,
d.\,h., es existieren nicht gleichzeitig verschwindende Lagrangesche Multiplikatoren $\lambda_0 \geq 0$, $y^* \in Y^*$ und $z^* \in Z^*$
derart,
dass folgende Bedingungen gelten:
\begin{enumerate}
\item[(a)] Die Lagrange-Funktion besitzt bez"uglich $x$ in $x_*$ einen station"aren Punkt, d.\,h.
           \begin{equation}\label{SatzEPschwach1}
           0 = \mathscr{L}_x(x_*,u_*,\lambda_0,y^*,z^*);
           \end{equation}   
\item[(b)] Die Lagrange-Funktion erf"ullt bez"uglich $u$ in $u_*$ die Variationsungleichung
           \begin{equation}\label{SatzEPschwach2}
           0 \leq \langle \mathscr{L}_u(x_*,u_*,\lambda_0,y^*,z^*), u-u_* \rangle \qquad \mbox{f"ur alle } u \in U;
           \end{equation}
\item[(c)] Die komplement"aren Schlupfbedingungen gelten, d.\,h.
           \begin{equation}\label{SatzEPschwach3}
           0 = \langle z^*, G(x_*) \rangle, \qquad z^* \in \mathscr{K}^*.
           \end{equation}
\end{enumerate}
\end{theorem}

Wir werden anschlie"send folgende Bezeichnungen verwenden:
$$L_0={\rm Im\,}\mathscr{F}_x(x_*,u_*) \subseteq Y,$$
die Wertemenge des stetigen linearen Operators $\mathscr{F}_x(x_*,u_*)$;
$$B = L_0 + \mathscr{F}_u(x_*,u_*)(U-u_*),$$
die Gesamtheit derjenigen $y \in Y$, zu denen es ein $x \in X$ und ein $u \in U$ gibt mit
$$y=\mathscr{F}_x(x_*,u_*)x+\mathscr{F}_u(x_*,u_*)(u-u_*).$$
Au"serdem bezeichnet $L= {\rm lin\,} B$ die lineare H"ulle der Menge $B$.
Nach Voraussetzung hat der Teilraum $L_0$ eine endliche Kodimension,
so dass $L_0$ und $L$ abgeschlossene Teilr"aume von $Y$ sind.
(Die Annahme, ${\rm Im\,}\mathscr{F}_x(x_*,u_*)$ sei abgeschlossen, w"are an dieser Stelle nicht ausreichend.)

\begin{lemma} \label{LemmaEPschwach}
Es sei $L=Y$.
Dann ist ${\rm int\,}B \not= \emptyset$.
Ist au"serdem $0 \in {\rm int\,}B$, so existieren $u_1,...,u_m \in U$ derart, dass
f"ur $z_j=\pi\big(\mathscr{F}_u(x_*,u_*)(u_j-u_*)\big)$ folgende Beziehungen gelten:
$${\rm lin\,}\{z_1,...,z_m\}= Y/L_0, \qquad z_1+...+z_m=0.$$
Dabei bezeichnet $\pi:Y \to Y/L_0$ die kanonische Abbildung,
d.\,h. es gilt $\pi y_1=\pi y_2$ genau dann, wenn $y_1-y_2 \in L_0$ ist.
\end{lemma}

{\bf Beweis:} Da die Kodimension der Menge $L_0$ endlich ist,
ist der Quotientenraum $Y/L_0$ endlichdimensional.
Die Menge $B$ ist offensichtlich konvex.
Deshalb ist auch $\pi(B)$ konvex.
Da die lineare H"ulle der Menge $B$ mit $Y$ "ubereinstimmt,
f"allt die lineare H"ulle der Menge $\pi(B)$ mit $Y/L_0$ zusammen.
Wegen der Konvexit"at von $\pi(B)$ muss daher die Menge $\pi(B)$ ein nichtleeres Inneres besitzen.
Au"serdem ist die Beziehung $\pi^{-1}\big(\pi(B)\big) = B$ erf"ullt.
Weil $\pi$ eine stetige Abbildung ist, bedeutet dies, dass ${\rm int\,}B \not= \emptyset$ gilt. \\
Weil $Y/L_0$ endlichdimensional ist,
existieren im Fall $0 \in {\rm int\,}B$ offenbar endlich viele Punkte $z_1,...,z_m$ aus $\pi(B)$ mit den geforderten Eigenschaften.
Wir w"ahlen dazu einfach $z_j$ als die Ecken eines hinreichend kleinen Würfels entsprechender Dimension, dessen Mittelpunkt im Ursprung liegt.
Nach Definition der Menge $\pi(B)$ gibt es ferner Elemente $u_j \in U$ mit $z_j=\pi\big(\mathscr{F}_u(x_*,u_*)(u_j-u_*)\big)$. \hfill $\blacksquare$ \\[2mm]
Wir wenden uns nun dem Beweis des Extremalprinzips \ref{SatzExtremalprinzipSchwach} zu.
Der Beweis ist in drei F"alle aufgeteilt, n"amlich zwei entarteten und einem nichtentarteten. \\[2mm]
{\bf Beweis im ersten entarteten Fall:} Es sei $L \not= Y$.
Dann existiert nach Folgerung \ref{FolgerungAnnulator} ein nichttriviales Funktional $y^* \in Y^*$ mit
$$0=\langle y^*, \mathscr{F}_x(x_*,u_*)x+\mathscr{F}_u(x_*,u_*)(u-u_*) \rangle$$
f"ur alle $x \in X$ und $u \in U$.
Setzen wir au"serdem $\lambda_0=0$ und $z^*=0$,
dann gelten die Bedingungen des Extremalprinzips \ref{SatzExtremalprinzipSchwach}. \\[2mm]
{\bf Beweis im zweiten entarteten Fall:} Es sei $L = Y$ und $0 \not\in {\rm int\,}B$.
Da die konvexe Menge $B$ in diesem Fall ein nichtleeres Inneres besitzt,
existiert nach dem Trennungssatz ein nichttriviales $y^* \in Y^*$,
das die Menge $B$ und das Nullelement trennt.
Dies bedeutet,
dass f"ur alle $x \in X$ und $u \in U$ die Ungleichung
$$0\leq \langle y^*, \mathscr{F}_x(x_*,u_*)x+\mathscr{F}_u(x_*,u_*)(u-u_*) \rangle$$
erf"ullt ist.
Setzen wir hierin $u=u_*$, so erhalten wir
$$0\leq \langle y^*, \mathscr{F}_x(x_*,u_*)x \rangle$$
f"ur alle $x \in X$. Folglich gilt $\mathscr{F}_x^*(x_*,u_*)y^*=0$.
Andererseits ergibt sich f"ur $x=0$ die Ungleichung
$$0\leq \langle y^*, \mathscr{F}_u(x_*,u_*)(u-u_*) \rangle$$
f"ur alle $u \in U$.
Wie im vorhergehenden Fall sind daher $\lambda_0=0$, $z^*=0$ und $y^*$ ein System gesuchter Multiplikatoren. \\[2mm]
{\bf Beweis im regul"aren Fall:} Es sei $L = Y$ und $0 \in {\rm int\,}B$.
Wir nehmen $J(x_*,u_*)=0$ an.
Im Weiteren sei $\mathscr{C}$ die Menge derjenigen
$(\eta_0,y,z) \in \R \times Y \times Z$ mit folgender Eigenschaft:
Zu jedem Element existieren ein $x \in X$, ein $u \in U$ und ein $\eta \in {\rm int\,}\mathscr{K}$ mit
$$\eta_0 > J'(x_*,u_*)(x,u-u_*), \quad y = \mathscr{F}'(x_*,u_*)(x,u-u_*), \quad z = G(x_*) + G'(x_*)x - \eta.$$
Zum Beweis des Extremalprinzips \ref{SatzExtremalprinzipSchwach} gen"ugt es im regul"aren Fall die Beziehungen
$${\rm int\,}\mathscr{C} \not= \emptyset, \qquad 0 \not\in {\rm int\,}\mathscr{C}$$
nachzuweisen.
Die Menge $\mathscr{C}$ ist offensichtlich konvex, da die Menge $U$ nach Voraussetzung konvex ist.
Wenn also diese Beziehungen gelten,
dann folgt die Existenz eines nichttrivialen Funktionals $(\lambda_0,y^*,z^*) \in \R \times Y^* \times Z^*$,
das die Menge ${\rm int\,}\mathscr{C}$ vom Ursprung trennt,
d.\,h. ein solches Funktional, dass
$$0 \leq \lambda_0 \eta_0 + \langle y^*,y \rangle + \langle z^*,z \rangle$$
f"ur alle $(\eta_0,y,z) \in\mathscr{C}$ gilt.
Beachten wir $J(x_*,u_*)=0$ und $\mathscr{F}(x_*,u_*)=0$,
so erhalten wir daraus, dass f"ur alle $x \in X$, $u \in U$ und alle $\eta \in \mathscr{K}$ die Beziehung
\begin{eqnarray*} 
       \langle z^*,\eta \rangle
&\leq& \lambda_0 [J(x_*,u_*) + J'(x_*,u_*)(x,u-u_*)] \\
&    & + \langle y^*, \mathscr{F}(x_*,u_*) + \mathscr{F}'(x_*,u_*)(x,u-u_*) \rangle + \langle z^*,G(x_*)+G'(x_*)x \rangle
\end{eqnarray*}
erf"ullt ist.
Betrachten wir diese Ungleichung zun"achst f"ur $x=0$ und $u=u_*$,
dann folgt
$$\langle z^*,\eta \rangle \leq \langle z^*,G(x_*) \rangle$$
f"ur alle $\eta \in {\rm int\,}\mathscr{K}$.
Da $\mathscr{K}$ ein abgeschlossener Kegel mit Spitze in Null ist, ergeben sich daraus die Beziehungen
$\langle z^*,z \rangle \leq 0$ f"ur alle $z \in \mathscr{K}$ und $\langle z^*,G(x_*) \rangle = 0$.
Damit ist (\ref{SatzEPschwach3}) gezeigt. \\
Zusammen ergibt sich so f"ur alle $x \in X$ und $u \in U$:
$$0 \leq \langle \mathscr{L}_x(x_*,u_*,\lambda_0,y^*,z^*), x \rangle + \langle \mathscr{L}_u(x_*,u_*,\lambda_0,y^*,z^*),u-u_* \rangle.$$
Betrachten wir nun nacheinander $u=u_*$ und $x=0$,
so kommen wir zu den Beziehungen (\ref{SatzEPschwach1}), (\ref{SatzEPschwach2}) des Extremalprinzips \ref{SatzExtremalprinzipSchwach}. \\[2mm]
Wir zeigen nun ${\rm int\,}\mathscr{C} \not= \emptyset$:
Es seien $u_1,...,u_m$ die Elemente mit den entsprechenden Eigenschaften in Lemma \ref{LemmaEPschwach}.
Im Weiteren seien 
$$V=\{x \in X\,|\, \|x\| <1\}, \qquad U_0={\rm conv\,} \{u_1,...,u_m\} \subseteq U$$
und
$$B_0=\mathscr{F}_x(x_*,u_*)V + \mathscr{F}_u(x_*,u_*)(U_0-u_*).$$
Die Menge $B_0$ ist konvex und besitzt ein nichtleeres Inneres,
da $\pi(B_0)$ die Punkte $z_1,...,z_m$ (vgl. Lemma \ref{LemmaEPschwach}) enth"alt.
Au"serdem ist nach dem Satz von der offenen Abbildung die Menge $\mathscr{F}_x(x_*,u_*)V$ offen in $L_0$.
Ferner seien
$$c_0 = \max_{j=1,...,m} \big(J_u(x_*,u_*)(u_j-u_*) + \|J_x(x_*,u_*)\| \big), \qquad \R_{c_0} = \{ a \in \R \,|\, a > c_0\}.$$
Dann besitzt die Menge $\mathscr{C}_0= \R_{c_0} \times B_0 \times (G(x_*)-{\rm int\,}\mathscr{K})$
ein nichtleeres Inneres und es gilt $\mathscr{C}_0 \subseteq \mathscr{C}$.
Demzufolge ist ${\rm int\,}\mathscr{C} \not= \emptyset$. \\[2mm]
Wir nehmen nun an, es ist $0 \in {\rm int\,}\mathscr{C}$.
Dann existieren Vektoren $\overline{x} \in X$, $\overline{u} \in U$, $\overline{\eta} \in {\rm int\,}\mathscr{K}$ und eine Zahl $c>0$ mit
\begin{eqnarray}
\label{BeweisWEPr2} -c &>& J'(x_*,u_*)(\overline{x},\overline{u}-u_*), \\
\label{BeweisWEPr3} 0 &=& \mathscr{F}'(x_*,u_*)(\overline{x},\overline{u}-u_*), \\
\label{BeweisWEPr4} 0 &=& G(x_*)+ G'(x_*)\overline{x} - \overline{\eta}.
\end{eqnarray}
Angenommen,
die Beziehungen (\ref{BeweisWEPr2})--(\ref{BeweisWEPr4}) seien erf"ullt.
Es sei $\varrho >0$ fest gew"ahlt.
Ferner seien $u_1,...,u_m \in U$ die Elemente in Lemma \ref{LemmaEPschwach}.
Dann liefert
$$u_*+\alpha_0(\overline{u}-u_*) + \varrho\sum_{j=1}^m \alpha_j (u_j-u_*)$$
f"ur $\alpha_0 + \varrho (\alpha_1 + ... + \alpha_m) \in [0,1]$ eine Konvexkombination.
Daher ist in einer Umgebung des Punktes $(x_*,0,0) \in X \times \R \times \R^m$ durch
$$\Phi(x,\alpha_0,\alpha) = \mathscr{F}\bigg(x,u_*+\alpha_0(\overline{u}-u_*) + \varrho\sum_{j=1}^m \alpha_j (u_j-u_*)\bigg),
  \quad \alpha=(\alpha_1,...,\alpha_m),$$
eine Abbildung in den Raum $Y$ definiert.
Dabei gilt $\Phi(x_*,0,0) = 0$.
Die Bedingung (A$_2$) liefert,
dass $\Phi$ auf einer Umgebung des Punktes $(x_*,0,0)$ Fr\'echet-differenzierbar und im Punkt $(x_*,0,0)$ stetig Fr\'echet-differenzierbar mit der Ableitung
$$\Phi'(x_*,0,0)(x,\alpha_0,\alpha) = \mathscr{F}'(x_*,u_*)\big(x,\alpha_0(\overline{u}-u_*)\big) + \varrho\sum_{j=1}^m \alpha_j \mathscr{F}_u(x_*,u_*)(u_j-u_*)$$
ist.
Ferner enth"alt die Wertemenge des linearen Operators $\Phi'(x_*,0,0)$ die Menge $B_0$ und stimmt folglich mit dem ganzen Raum $Y$ "uberein.
Au"serdem gibt es nach Wahl der Elemente $u_j$ in Lemma \ref{LemmaEPschwach} ein $x' \in X$ mit
$$\mathscr{F}_x(x_*,u_*)x' + \sum_{j=1}^m \mathscr{F}_u(x_*,u_*)(u_j-u_*) =0.$$
Damit folgt zusammen mit (\ref{BeweisWEPr3})
$$\Phi'(x_*,0,0)(\overline{x} + \varrho x',1,1)=0,$$
d.\,h., dass der Punkt $(\overline{x}+ \varrho x',1,1)$ dem Kern des Operators $\Phi'(x_*,0,0)$ angeh"ort.
Nach dem Satz von Ljusternik (Theorem \ref{SatzLjusternik}) existieren eine Zahl $\varepsilon_0 >0$ und Abbildungen
$\varepsilon \to \big(x(\varepsilon),\alpha_0(\varepsilon),\alpha(\varepsilon)\big)$ des Intervalls $[0,\varepsilon_0]$ in den Raum
$X \times \R \times \R^m$ derart, dass
$$\lim_{\varepsilon \to 0^+} \big\|\big(x(\varepsilon),\alpha_0(\varepsilon),\alpha(\varepsilon)\big)\big\|=0$$
gilt und au"serdem
\begin{equation} \label{BeweisWEPr5}
\Phi\big(x_* + \varepsilon [\overline{x} + \varrho x' + x(\varepsilon)],\varepsilon[1 + \alpha_0(\varepsilon)],
         \varepsilon[1 + \alpha(\varepsilon)]\big) =\Phi(x_*,0,0)
\end{equation}
f"ur alle $\varepsilon \in [0,\varepsilon_0]$ erf"ullt ist.
Wir setzen der K"urze halber
\begin{eqnarray*}
\tilde{x}(\varepsilon) &=& x_* + \varepsilon [\overline{x} + \varrho x' + x(\varepsilon)], \\
\tilde{u}(\varepsilon) &=& u_*+[\varepsilon + \alpha_0(\varepsilon)](\overline{u}-u_*) + \varrho\sum_{j=1}^m [\varepsilon + \alpha_j(\varepsilon)] (u_j-u_*).
\end{eqnarray*}
Darin stellt $\tilde{u}(\varepsilon)$ f"ur alle $\varepsilon \in [0,\varepsilon_0]$ mit hinreichend kleinem $\varepsilon_0>0$ eine
Konvexkombination der Elemene $u_*,\overline{u},u_1,...,u_m$ dar und es folgt unmittelbar aus (\ref{BeweisWEPr5}):
\begin{equation} \label{BeweisWEPr6}
\mathscr{F}\big(\tilde{x}(\varepsilon),\tilde{u}(\varepsilon)\big) = 0.
\end{equation}
Da die Abbildung $G$ im Punkt $x_*$ Fr\'echet-differenzierbar ist,
ergibt sich
$$G\big(\tilde{x}(\varepsilon)\big) = G(x_*) + \varepsilon G'(x_*)(\overline{x}+\varrho x') + r(\varepsilon), \quad
  \lim_{\varepsilon \to 0^+} \frac{\|r(\varepsilon)\|}{\varepsilon}=0.$$
Wir beachten,
dass $\mathscr{K}$ konvex ist,
sowie $G(x_*) \in \mathscr{K}$ und $G(x_*)+ G'(x_*)\overline{x} = \overline{\eta} \in {\rm int\,}\mathscr{K}$ nach (\ref{BeweisWEPr4}) gelten.
Ferner sei $\varrho >0$ derart gew"ahlt, dass $\overline{\eta} + \varrho G'(x_*) x' \in {\rm int\,}\mathscr{K}$ erf"ullt ist. 
Damit ergibt sich
\begin{eqnarray*}
G\big(\tilde{x}(\varepsilon)\big) &=& (1-\varepsilon) G(x_*) + \varepsilon [G(x_*)+ G'(x_*)\overline{x}]
                                      +\varepsilon \varrho  G'(x_*) x' + r(\varepsilon) \\
                                  &=& (1-\varepsilon) G(x_*) + \varepsilon [\overline{\eta} + \varrho G'(x_*) x' + r(\varepsilon) / \varepsilon].
\end{eqnarray*}
Wir erhalten daraus f"ur alle $\varepsilon \in [0,\varepsilon_0]$ mit einem hinreichend kleinen $\varepsilon_0$:
\begin{equation} \label{BeweisWEPr7}
G\big(\tilde{x}(\varepsilon)\big) \in {\rm int\,}\mathscr{K}.
\end{equation}
Zus"atzlich zu der Bedingung $\overline{\eta} + \varrho G'(x_*) x' \in {\rm int\,}\mathscr{K}$ w"ahlen wir $\varrho>0$ so,
dass mit der Zahl $c >0$ in (\ref{BeweisWEPr2}) die Relation
$$\varrho \bigg[ J_x(x_*,u_*)x' + \sum_{j=1}^m J_u(x_*,u_*)(u_j-u_*) \bigg] \leq  \frac{c}{2}$$
erf"ullt ist.
Da das Funktional $J$ im Punkt $(x_*,u_*)$ Fr\'echet-differenzierbar ist, ergibt sich
\begin{eqnarray*}
    J\big(\tilde{x}(\varepsilon),\tilde{u}(\varepsilon)\big)
&=& J(x_*,u_*) + \varepsilon J'(x_*,u_*)(\overline{x},\overline{u}-u_*) \\
& &   + \varepsilon \varrho \bigg[ J_x(x_*,u_*)x' + \sum_{j=1}^m J_u(x_*,u_*)(u_j-u_*) \bigg] + o(\varepsilon).
\end{eqnarray*}
Mit (\ref{BeweisWEPr2}) und nach Wahl von $\varrho>0$ erhalten wir f"ur alle $\varepsilon \in [0,\varepsilon_0]$:
\begin{equation} \label{BeweisWEPr8}
J\big(\tilde{x}(\varepsilon),\tilde{u}(\varepsilon)\big) \leq J(x_*,u_*) - \varepsilon c + \varepsilon \frac{c}{2} + o(\varepsilon).
\end{equation}
Die Beziehungen (\ref{BeweisWEPr6}) und (\ref{BeweisWEPr7}) zeigen,
dass f"ur hinreichend kleine $\varepsilon >0$ das Paar $\big(\tilde{x}(\varepsilon),\tilde{u}(\varepsilon)\big)$ zul"assig in der Aufgabe (\ref{EAAnhang2}) ist.
Au"serdem ist $J\big(\tilde{x}(\varepsilon),\tilde{u}(\varepsilon)\big) < J(x_*,u_*)$ nach (\ref{BeweisWEPr8}) f"ur hinreichend kleine $\varepsilon >0$.
Wegen $\tilde{x}(\varepsilon) \to x_*$ und $\tilde{u}(\varepsilon) \to u_*$ f"ur $\varepsilon \to 0$ bedeutet dies,
dass der Punkt $(x_*,u_*)$ im Widerspruch zur Voraussetzung kein schwaches lokales Minimum sein k"onnte. \hfill $\blacksquare$   
\end{appendix}       
   
\addcontentsline{toc}{section}{Literatur}
\lhead[\thepage \, Literatur]{Optimale Steuerung mit unendlichem Zeithorizont}
\rhead[Optimale Steuerung mit unendlichem Zeithorizont]{Literatur \thepage}

\end{document}